\documentclass{report}
\usepackage{amsbsy,amssymb,amscd,amsfonts,latexsym,amstext,delarray,
amsmath,epsfig}
\usepackage{epsfig}
\input psfig.sty
\input{amssymb.sty}

\def\n{{\noindent}}
\def\no{{\noindent}}

\newtheorem{thm}{Theorem}[chapter]
\newtheorem{prop}[thm]{Proposition}
\newtheorem{cor}[thm]{Corollary}
\newtheorem{lem}[thm]{Lemma}

\newtheorem{defn}[thm]{Definition}
\newtheorem{rem}[thm]{Remark}

\numberwithin{equation}{chapter}

\newcommand{\ie}{{\it i.e.\/}\ }
\newcommand{\eg}{{\it e.g.\/}\ }
\newcommand{\cf}{{\it cf.\/}\ }

\def\Cb{{\mathbb C}}
\def\Qb{{\mathbb Q}}
\def\Hb{{\mathbb H}}
\def\Nb{{\mathbb N}}

\def\Rb{{\mathbb R}}

\def\Zb{{\mathbb Z}}

\def\L{\Lambda}
\def\Ac{{\mathcal A}}
\def\Bc{{\mathcal B}}
\def\Ac{{\mathcal A}}
\def\Hc{{\mathcal H}}
\def\Lc{{\mathcal L}}

\def\Oc{{\mathcal O}}

\def\Uc{{\mathcal U}}

\def\Xc{{\mathcal X}}

\def\Ec{{\mathcal E}}
\def\Hc{{\mathcal H}}

\def\Lc{{\mathcal L}}
\def\Rc{{\mathcal R}}

\def\Uc{{\mathcal U}}

\def\a{\alpha}
\def\b{\beta}
\def\d{\delta}

\def\g{\gamma}

\def\t{\theta}
\def\ve{\varepsilon}

\def\z{\zeta}

\def\D{\Delta}
\def\G{\Gamma}

\def\part{\partial}

\def\RR{R}

\def\qqq{\,,\quad \forall}

\def\text{\hbox}

\def\Aut{\mathop{\rm Aut}\nolimits}

\def\Trace{\mathop{\rm Trace}\nolimits}

\def\bP{{\mathbb P}}
\def\R{{\mathbb R}}
\def\Q{{\mathbb Q}}
\def\C{{\mathbb C}}
\def\Z{{\mathbb Z}}
\def\N{{\mathbb N}}
\def\A{{\mathbb A}}

\def\GL{{\rm GL}}
\def\SL{{\rm SL}}
\def\Gal{{\rm Gal}}

\def\Tr{{\rm Tr}}

\def\sA{{\mathcal A}}
\def\cA{{\mathcal A}}

\def\sE{{\mathcal E}}

\def\cH{{\mathcal H}}

\def\H{{\mathbb H}}

\def\n{{\noindent}}

\def\Cb{{\mathbb C}}
\def\Qb{{\mathbb Q}}
\def\Hb{{\mathbb H}}
\def\Nb{{\mathbb N}}
\def\Rb{{\mathbb R}}

\def\Zb{{\mathbb Z}}

\def\L{\Lambda}
\def\Ac{{\mathcal A}}
\def\Bc{{\mathcal B}}
\def\Ac{{\mathcal A}}
\def\Hc{{\mathcal H}}
\def\Lc{{\mathcal L}}

\def\Oc{{\mathcal O}}

\def\Uc{{\mathcal U}}

\def\Ec{{\mathcal E}}
\def\Hc{{\mathcal H}}

\def\Lc{{\mathcal L}}
\def\Rc{{\mathcal R}}

\def\Uc{{\mathcal U}}

\def\a{\alpha}
\def\b{\beta}
\def\d{\delta}

\def\g{\gamma}

\def\t{\theta}
\def\ve{\varepsilon}

\def\z{\zeta}

\def\D{\Delta}
\def\G{\Gamma}

\def\part{\partial}

\def\text{\hbox}

\def\Aut{\mathop{\rm Aut}\nolimits}

\def\Trace{\mathop{\rm Trace}\nolimits}

\def\bP{{\mathbb P}}
\def\R{{\mathbb R}}
\def\Q{{\mathbb Q}}
\def\C{{\mathbb C}}
\def\Z{{\mathbb Z}}
\def\N{{\mathbb N}}
\def\A{{\mathbb A}}

\def\GL{{\rm GL}}
\def\SL{{\rm SL}}
\def\Gal{{\rm Gal}}

\def\Tr{{\rm Tr}}

\def\sA{{\mathcal A}}
\def\cA{{\mathcal A}}

\def\sE{{\mathcal E}}

\def\cH{{\mathcal H}}

\def\H{{\mathbb H}}

\title
{From Physics to Number theory via Noncommutative Geometry}
\author{Alain Connes and Matilde Marcolli}
\date{}

\begin{document}

\maketitle

\tableofcontents

\chapter*{Introduction}

\begin{verse}
{\em e volta nostra poppa nel mattino, \\
de' remi facemmo ali al folle volo  \\
------ Dante, Inf. XXVI 124-125 }
\end{verse}

\n Several recent results reveal a surprising connection between
modular forms and noncommutative geometry. The first
occurrence came from the classification of noncommutative three
spheres, \cite{0CDV1} \cite{0CDV2}. Hard computations with the
noncommutative analog of the Jacobian involving the ninth power of
the Dedekind eta function were necessary in order to analyze the
relation between such spheres and noncommutative nilmanifolds.
Another occurrence can be seen in the computation of the explicit
cyclic cohomology Chern character of a spectral triple on
$SU_q(2)$ \cite{0CoQgr}. Another
surprise came recently from a remarkable action of the Hopf
algebra of transverse geometry of foliations of codimension one on
the space of lattices modulo Hecke correspondences, described in
the framework of noncommutative geometry, using a
modular Hecke algebra obtained as the cross product of modular forms
by the action of Hecke
correspondences  \cite{0CM1} \cite{0CM2}. This
action determines a differentiable structure on this
noncommutative space, related to the Rankin--Cohen brackets of
modular forms, and shows their compatibility with Hecke operators.
Another instance where properties of modular forms can be recast
in the context of noncommutative geometry can be found in the
theory of modular symbols and Mellin transforms of cusp forms of
weight two, which can be recovered from the geometry of the moduli
space of Morita equivalence classes of noncommutative tori viewed
as boundary of the modular curve \cite{0MM}.

\smallskip

\n In this paper we show that
the theory of modular Hecke algebras, the spectral
realization of zeros of $L$-functions, and the arithmetic properties
of KMS states in quantum statistical mechanics combine into a
unique general picture based on the noncommutative geometry
of the space of commensurability classes of $\Q$-lattices.

\smallskip

\n An n-dimensional $\Q$-lattice consists of an ordinary lattice
$\Lambda$ in $\R^n$ and a homomorphism
$$ \phi: \Q^n/\Z^n \to \Q\Lambda/\Lambda. $$
Two such $\Q$-lattices are {\em commensurable} if and only if the
corresponding lattices are commensurable and the maps agree modulo
the sum of the lattices.

\smallskip

\n The description of the spaces of commensurability classes of
$\Q$-lattices via noncommutative geometry yields two quantum
systems related by a duality. The first system is of quantum
statistical mechanical nature, with the algebra of coordinates
parameterizing commensurability classes of $\Q$-lattices modulo
scaling and with a time evolution with eigenvalues given by the
index of pairs of commensurable $\Q$-lattices. There is a symmetry
group acting on the system, in general by {\em endomorphisms}. It
is this symmetry that is spontaneously broken at low temperatures,
where the system exhibits distinct phases parameterized by
arithmetic data. We completely analyze the phase transition with
spontaneous symmetry breaking in the two-dimensional case, where a
new phenomenon appears, namely that there is a second critical
temperature, beyond which no equilibrium state survives.

\smallskip

\n In the ``dual system'', which corresponds just to
commensurability of $\Q$-lattices, the scaling group is acting. In
physics language, what emerges is that the zeros of zeta appear as
an absorption spectrum of the scaling action in the $L^2$ space of
the space of commensurability classes of $\Q$-lattices as in
\cite{0Co-zeta}. While the zeros of zeta and $L$-functions appear
at the critical temperature, the analysis of the low temperature
equilibrium states concentrates on the subspace
$$ \GL_n(\Q)\backslash \GL_n(\A) $$
of {\em invertible} $\Q$-lattices, which as is well known plays a
central role in the theory of automorphic forms.

\smallskip

\n While, at first sight, at least in the
1-dimensional case, it would seem easy to classify
commensurability classes of $\Q$-lattices, we shall see that ordinary geometric
tools fail because of the ergodic nature of the equivalence
relation.

\smallskip

\n Such quotients are fundamentally of ``quantum nature'', in that,
even though they are sets in the ordinary sense, it is impossible
to
 distinguish points
 by any finite (or countable) collection of invariants.
Noncommutative geometry is specifically designed to handle such
quantum spaces by encoding them by algebras of non-commuting
coordinates and extending the techniques of ordinary geometry
using the tools of functional analysis, noncommutative algebra,
and quantum physics.

\smallskip

\n Direct attempts to define function spaces for such quotients
lead to invariants that are of a cohomological nature. For
instance, let the fundamental group $\Gamma$ of a Riemann surface
act on the boundary $\bP^1(\R)$ of its universal cover identified
with the Poincar\'e disk. The space
$$ L^\infty(\Gamma\backslash \bP^1(\R)) :=L^\infty(\bP^1(\R))^\G $$
is in natural correspondence with global sections of the sheaf of
(real parts of) holomorphic functions on the Riemann surface, as
boundary values. More generally, the cyclic cohomology of the
noncommutative algebra of coordinates on such quotients is
obtained by applying derived functors to these naive functorial
definition of function spaces.

\smallskip

\n In the 1-dimensional case, the states at zero temperature are
related to the Kronecker--Weber construction of the maximal
abelian extension $\Q^{ab}$. In fact, in this case the quantum
statistical mechanical system is the one constructed in
\cite{0BC}, which has underlying geometric space $X_1$
parameterizing commensurability classes of 1-dimensional
$\Q$-lattices modulo scaling by $\R^\ast_+$. The corresponding algebra of
coordinates is a Hecke algebra for an almost normal pair of
solvable groups. The regular representation is of type III$_1$ and
determines the time evolution of the system, which has the set of
$\log (p)$, $p$ a prime number, as set of basic frequencies. The
system has an action of the id\`eles class group modulo the connected
 component of identity as a group of symmetries. This induces a Galois
action on the ground states of the system at zero temperature.
When raising the temperature the system has a phase transition,
with a unique equilibrium state above the critical temperature.
The Riemann zeta function appears as the partition function of the
system, as in \cite{julia0}.

\smallskip

\n Each equivalence class of $\Q$-lattices determines an
irreducible covariant representation, where the Hamiltonian is
implemented by minus the log of the covolume. For a general class,
this is not bounded below. It is so, however, in the case of
equivalence classes of {\em invertible} $\Q$-lattices, \ie where
the labelling of torsion points is one to one. These classes then
define positive energy representations and corresponding KMS
states for all temperatures below critical.

\smallskip

\n In the 2-dimensional case, as the temperature lowers, the
system settles down on these invertible $\Q$-lattices, so that the
zero temperature space is commutative and is given by the Shimura
variety
$$ \GL_2(\Q) \backslash \GL_2(\A) / \C^*. $$
The action of the symmetry group, which in this case is
 nonabelian and isomorphic to
$\Q^*\backslash \GL_2(\A_f)$, is more subtle
due to the presence of inner automorphisms
and the necessary use of the formalism
of superselection sectors.
Moreover, its effect on
 the zero
temperature states is not obtained directly
but is induced by the action at non-zero
temperature, which involves the full noncommutative system.
The quotient $\GL_2(\Q)\backslash (M_2(\A_f) \times \GL_2(\R))
/\C^*$ and the space of 2--dimensional $\Q$-lattices modulo
commensurability and scaling are the same, hence the corresponding
algebras are Morita equivalent. However, it is preferable to work
with the second description, since, by taking the classical
quotient by the action of the subgroup $\SL_2(\Z)$, it reduces the group part
in the cross product to the classical Hecke algebra.

\smallskip

\n The $\GL_2$ system has an arithmetic structure provided by a
{\em rational} subalgebra, given by a natural condition on the
coefficients of the $q$-series. We show that it is a Hecke algebra
of modular functions, closely related to the modular Hecke algebra
of \cite{0CM1}, \cite{0CM2}. The symmetry group acts on the values
of ground states on this rational subalgebra as the automorphism
group of the modular field.

\smallskip

\n Evaluation of a generic ground state $\varphi$ of the system on
the rational subalgebra generates an embedded copy of the modular
field in $\C$ and there exists a unique isomorphism of the
symmetry group of the system with the Galois group of the embedded
modular field, which intertwines the Galois action on the image
with the symmetries of the system,
$$ \theta(\sigma) \circ \varphi = \varphi \circ \sigma. $$
The relation between this $\GL_2$ system and class field theory
is being investigated in ongoing work \cite{0CMR}.

\smallskip

\n The arithmetic structure is inherited by the dual of the
$\GL_2$ system and enriches the structure of the noncommutative
space of commensurability classes of 2-dimensional $\Q$-lattices
to that of a ``noncommutative arithmetic variety''. The relation
between this dual system and the spectral realization of zeros of
$L$-functions is the central topic of Chapter 2.

\smallskip

\n The dual of the $\GL_1$ system, under the duality obtained by
taking the cross product by the time evolution, corresponds to
the space of commensurability classes of 1-dimensional
$\Q$-lattices, not considered up to scaling. This corresponds
geometrically to the total space $\Lc$ of a principal $\R^*_+$
bundle over the base $X_1$, and determines a natural scaling
action of $\R^*_+$. The space $\Lc$ is described by the quotient
$$ \Lc= \GL_1(\Q)\backslash \A^\cdot, $$
where $\A^\cdot$ denotes the set of ad\`eles with nonzero
archimedean component. The corresponding algebra of coordinates is
Morita equivalent to $C(X_1) \rtimes_{\sigma_t} \R$.

\smallskip

\n Any approach to a spectral realization of the zeros of zeta
through the quantization of a classical dynamical system faces the
problem of obtaining the leading term in the Riemann counting
function for the number of zeros of imaginary part less than $E$
as a volume in phase space. The solution \cite{0Co-zeta} of this
issue is achieved
in a remarkably simple way, by the scaling action of $\R^*_+$ on the
phase space of the real line $\R$, and will be the point of
departure for the second part of the paper.

\smallskip

\n In particular, this shows that the space $\Lc$ requires a
further compactification at the archimedean place, obtained by
replacing the quotient $\Lc= \GL_1(\Q)\backslash \A^\cdot$ by
$\overline{\Lc}= \GL_1(\Q)\backslash \A$
\ie dropping the non vanishing of the archimedean component.
This compactification
has an analog for the $\GL_2$ case, given by the noncommutative
boundary of modular curves considered in \cite{0MM}, which
corresponds to replacing $\GL_2(\R)$ by $M_2(\R)$ at the
archimedean place, and is related to class field theory for real
quadratic fields through Manin's real multiplication program.

\smallskip

\n The space $\overline{\Lc}$ appears as the configuration space
for a quantum field theory, where the degrees of freedom are
parameterized by prime numbers, including infinity. When only
finitely many degrees of freedom are considered, and in particular
only the place at infinity, the semiclassical approximation
exhibits the main terms in the asymptotic formula for the number
of zeros of the Riemann zeta function.

\smallskip

\n The zeros of zeta appear as an absorption spectrum, namely as
lacunae in a continuous spectrum, where the width of the
absorption lines depends on the presence of a cutoff. The full
id\`eles class group appears as symmetries of the system and
$L$-functions with Gr\"ossencharakter replace the Riemann zeta
function in nontrivial sectors.

\smallskip

\n From the point of view of quantum field theory, the field
configurations are given by ad\`eles, whose space $\A$ is then
divided by the action of the gauge group $\GL_1(\Q)$. As mentioned
above, the quotient space is essentially the same as the space
$\Lc$ of commensurability classes of 1-dimensional $\Q$-lattices.
The $\log(p)$ appear as periods of the orbits of the scaling
action. The Lefschetz formula for the scaling action recovers the
Riemann--Weil explicit formula as a semi-classical approximation.
The exact quantum calculation for finitely many degrees of freedom
confirms this result. The difficulty in extending this calculation
to the global case lies in the quantum field theoretic problem of
passing to infinitely many degrees of freedom.

\smallskip

\n The main features of the dual
systems in the $\GL_1$ case are summarized in the following table:
\bigskip

\begin{tabular}{|c|c|}\hline & \\
Quantum statistical mechanics & Quantum field theory \\[2mm]
\hline & \\
Commensurability classes & Commensurability classes \\
of $\Q$-lattices modulo scaling & of $\Q$-lattices \\[2mm] \hline & \\
 $A = C^*(\Q/\Z)\rtimes \N^\times$ &
$A \rtimes_{\sigma_t} \R$ \\ & \\ \hline & \\
Time evolution $\sigma_t$ & Energy scaling $U(\lambda)$,
$\lambda \in \R^*_+$ \\[2mm]
\hline & \\
$\{ \log p \}$ as frequencies & $\{ \log p \}$ as periods of
orbits \\[2mm] \hline & \\
 Arithmetic rescaling  $\mu_n$  & Renormalization group flow
$\mu
\partial_\mu$ \\[2mm]  \hline & \\
Symmetry group  $\hat \Z^*$ & Id\`eles class group \\
as Galois action on $T=0$ states &  as
gauge group \\[2mm] \hline & \\
System at zero temperature &  $\GL_n(\Q)\backslash
\GL_n(\A)$ \\[2mm]
\hline & \\
System at critical temperature & Spectral realization \\
 (Riemann's $\zeta$ as partition function) & (Zeros
of $\zeta$ as absorption spectrum)
\\[2mm] \hline & \\
Type III$_1$ & Type II$_\infty$ \\[2mm]
\hline
\end{tabular}

\bigskip
\bigskip

\n There is a similar duality (and table) in the $\GL_2$ case,
where part of the picture remains to be clarified. The relation
with the modular Hecke algebra of \cite{0CM1}, \cite{0CM2} is more
natural in the dual system where modular forms with non-zero
weight are naturally present.

\smallskip

\n The dual system $\Lc$ can be interpreted physically as a
``universal scaling system'', since it exhibits the continuous
renormalization group flow and its relation with the discrete
scaling by powers of primes. For the primes two and three, this
discrete scaling manifests itself in acoustic systems, as is well
known in western classical music, where the two scalings
correspond, respectively, to passing to the octave (frequency
ratio of 2) and transposition (the perfect fifth is the frequency
ratio 3/2), with the approximate value $\log(3)/\log(2)\sim 19/12$
responsible for the difference between the ``circulating
temperament'' of the Well Tempered Clavier and the ``equal
temperament'' of XIX century music. It is precisely the
irrationality of $\log(3)/\log(2)$ which is responsible for the
noncommutative nature of the quotient corresponding to the three
places $\{2, 3, \infty\}$.

\smallskip

\n The key role of the continuous renormalization group flow as a
symmetry of the dual system $\Lc$ and its similarity with a Galois
group at the archimedean place brings us to the last part of this
paper. In Chapter 3, we analyze the quantum statistical mechanics
of $\Q$-lattices at critical temperature.

\smallskip

\n The fact that the KMS state at critical temperature can be
expressed as a noncommutative residue (Dixmier trace) shows that
the system at critical temperature should be analyzed with tools
from quantum field theory and renormalization.

\smallskip

\n The mathematical theory of renormalization in QFT developed in
\cite{0CK1} \cite{0CK2} shows in geometric terms that the
procedure of perturbative renormalization can be described as the
Birkhoff decomposition
$$ g_{\rm eff}(\ve) =  g_{{\rm eff}_+} (\ve)\, (g_{{\rm
eff}_-}(\ve))^{-1} $$
on the projective line of
complexified dimensions  $\ve$
of the  loop
$$
g_{\rm eff}(\ve)\in G\,=\,\text{ formal diffeomorphisms of }\C
$$
given by the unrenormalized effective coupling constant. The
$g_{{\rm eff}_-}$ side of the Birkhoff decomposition yields the
counterterms and the $g_{{\rm eff}_+}$ side evaluated at the
critical dimension gives the renormalized value of the effective
coupling. This explicit knowledge of the counterterms suffices to
determine the full renormalized theory.

\smallskip

\n The principal $G$-bundle on $\bP^1(\C)$, with trivialization
given by the Birkhoff decomposition, has a flat connection with
regular singularities, coming from a Riemann--Hilbert problem
determined by the representation datum given by the
$\beta$-function of renormalization, viewed as the logarithm of
the monodromy around the singular dimension. The
problem of incorporating nonperturbative effects leads to a
more sophisticated Riemann--Hilbert problem in terms of a
representation of the wild fundamental group of Martinet--Ramis,
related to applying Borel summation techniques to the
unrenormalized effective coupling constant.

\smallskip

\n In the perturbative theory the renormalization group appears as
a natural 1-parameter subgroup of the group of ``diffeographisms''
which governs the ambiguity in the choice of the physical
solution. In fact, the renormalization group, the wild fundamental
group and the connected component of identity in the id\`eles
class group are all incarnations of a still mysterious Galois
theory at the archimedean place.

\smallskip

\n It was shown in \cite{0CoZ} that the classification of
approximately finite factors provides a nontrivial Brauer theory
for central simple algebras over $\C$, and an archimedean analog
of the module of central simple algebras over nonarchimedean
fields. The relation of Brauer theory to the Galois group is via
the construction of central simple algebras as cross products of
a field by a group of automorphisms. It remained for a long time
an elusive point to obtain in a natural manner
factors as cross product  by a group of automorphisms
of a {\it field}, which is a
transcendental extension of $\C$.
This was achieved in \cite{0CDV2} for type II$_1$, via the cross
product of the field of elliptic functions by an automorphism
given by translation on the elliptic curve. The results on the
$\GL_2$ system give an analogous construction for type
III$_1$  factors using the modular field.

\smallskip

\n The physical reason for considering such field extensions of
$\C$ lies in the fact that the coupling constants $g$ of the
fundamental interactions (electromagnetic, weak and strong)
are not really constants but
 depend on the energy scale $\mu$ and are therefore
functions $g(\mu)$. Thus, high energy physics implicitly extends
the ``field of constants'', passing from the field of scalars $\C$
to a field of functions containing all the $g(\mu)$. On this
field, the renormalization group provides the corresponding theory
of ambiguity, and acquires an interpretation as the missing Galois
group at archimedean places.

\medskip

\n The structure of the paper is organized as follows.
\begin{itemize}
\item
The first part is dedicated to the quantum statistical
mechanical system of $\Q$-lattices, in the cases of dimension one and
two, and its
behavior at zero temperature.

\item The second part deals
with the $\GL_1$ system at critical temperature and its dual
system, and their relation to the spectral realization of the
zeros of zeta of \cite{0Co-zeta}.

\item The last part is dedicated to the system at critical
temperature, the theory of renormalization of \cite{0CK1},
\cite{0CK2}, the Riemann-Hilbert problem and the missing Galois
theory at the archimedean place.

\end{itemize}

\bigskip

\bibliographystyle{alpha}

\bigskip
\chapter{ Quantum Statistical Mechanics of $\Q$-Lattices}

\section{Introduction}

\noindent In this chapter we shall start by giving a geometric
interpretation in terms of the space of commensurability classes
of $1$-dimensional $\Q$-lattices of the quantum statistical
dynamical system (BC \cite{BC}). This system exhibits the relation
between the phenomenon of spontaneous symmetry breaking and number
theory. Its dual system obtained by taking the cross product by
the time evolution is basic in the spectral interpretation of
zeros of zeta.

\smallskip

\n Since $\Q$-lattices and commensurability continue to make sense
in dimension $n$, we shall obtain an analogous system in higher
dimension and  in particular we derive a complete picture of the
system in dimension $n=2$. This shows two distinct phase
transitions with arithmetic spontaneous symmetry breaking.

\smallskip

\n In the initial model of BC (\cite{BC}) the partition function
is the Riemann zeta function. Equilibrium states are characterized
by the KMS-condition. While at large temperature there is only one
equilibrium state, when the temperature gets smaller than the
critical temperature, the equilibrium states are no longer unique
but fall in distinct phases parameterized by number theoretic
data. The pure phases are parameterized by the various embeddings
of the cyclotomic field $\Q^{ab}$ in $\Cb$.

\smallskip
\n The physical observables of the BC system form a
$C^\ast$-algebra endowed with a natural time evolution $\sigma_t$.
This algebra is interpreted here as the algebra of noncommuting
coordinates on the space of commensurability classes of
$1$-dimensional $\Q$-lattices up to scaling by $\R^\ast_+$.

\smallskip

\n What is remarkable about the ground states of this system is
that, when evaluated on the rational observables of the system,
they only affect values that are algebraic numbers. These span the
maximal abelian extension of $\Q$. Moreover, the class field
theory isomorphism intertwines the two actions of the id\`eles
class group, as symmetry group of the system, and of the Galois
group, as permutations of the expectation values of the rational
observables. That the latter action preserves positivity is a rare
property of states. We abstract this property as a definition of
``fabulous\footnote{This terminology is inspired from John
Conway's talk on ``fabulous" groups.} states'', in the more
general context of arbitrary number fields and review recent
developments in the direction of extending this result to other
number fields.

\smallskip

\n We present a new approach, based on the construction of an
analog of the BC system in the $\GL_2$ case.
 Its
relation to the complex multiplication case of the Hilbert 12th
problem will be discussed specifically in ongoing work of the two
authors with N. Ramachandran \cite{CMR}.

\smallskip

\n The $C^*$-algebra of observables in the $\GL_2$-system
describes the non-commutative space of commensurability classes of
$\Q$-lattices in $\C$ up to scaling by $\C^\ast$.

\smallskip

\n A $\Q$-lattice in $\C$ is a pair $(\L,\phi)$ where $\L \subset
\C$ is a lattice while
$$\displaystyle \phi :  \mathbb{Q}^2\slash\mathbb{Z}^2
\longrightarrow\mathbb{Q}\L \slash \L $$ is a homomorphism of
abelian groups (not necessarily invertible). Two $\Q$-lattices
$(\L_j,\phi_j)$ are commensurable iff the lattices $\L_j$ are
commensurable (\ie $\Q\Lambda_1=\Q\Lambda_2$) and the maps
$\phi_j$ are equal modulo $\L_1+\L_2$. The time evolution
corresponds to the ratio of covolumes of pairs of commensurable
$\Q$-lattices. The group
$$
S=\,\Q^\ast\backslash {\rm GL}_2 ({\A_f})
$$
quotient of  the finite ad\`elic group of $\GL_2$ by the
multiplicative group $\Q^\ast$ acts as symmetries of the system,
and the action is implemented by endomorphisms, as in the theory
of superselection sectors of Doplicher-Haag-Roberts (\cite{haag}).

\smallskip

\n It is this symmetry which is spontaneously broken below the
critical temperature $T= \frac{1}{2}$. The partition function of
the $\GL_2$ system is $\zeta(\beta)\zeta(\beta-1)$, for
$\beta=1/T$, and the system exhibits three distinct phases, with
two phase transitions at $T=\frac{1}{2}$ and at $T=1$. At low
temperatures ($T<\frac{1}{2}$) the pure phases are parameterized
by the set
$$ \GL_2(\Q)\backslash \GL_2(\A)/ \C^* $$
of classes of invertible $\Q$-lattices (up to scaling). The
equilibrium states of the ``crystalline phase'' merge as $T\to
1/2$ from below, as the system passes to a ``liquid phase'', while
at higher temperatures ($T\geq 1$) there are no  KMS states.

\smallskip

\n The subalgebra of {\em rational} observables turns out to be
intimately related to the modular Hecke algebra introduced in
(Connes-Moscovici \cite{CM10}) where its surprising relation with
transverse geometry of foliations is analyzed. We show that the
KMS states at zero temperature when evaluated on the rational
observables generate a specialization of the modular function
field $F$. Moreover, as in the BC system the state intertwines the
two actions of the group $S$,
 as symmetry group of the system, and
as permutations of the expectation values of the rational
observables by the Galois group of the modular field, identified
with $S$ by Shimura's theorem (\cite{Shimura}).

\smallskip

\n We shall first explain the general framework of quantum
statistical mechanics, in terms of $C^*$-algebras and KMS states.
Noncommutative algebras concretely represented in Hilbert space
inherit a canonical time evolution, which allows for phenomena of
phase transition and spontaneous symmetry breaking for KMS states
at different temperatures.

\smallskip

\n There are a number of important nuances between the abelian BC
case and the higher dimensional non-abelian cases. For instance,
in the abelian case, the subfield of $\C$ generated by the image
of the rational subalgebra under an extremal KMS$_\infty$ state
does not depend on the choice of the state and the intertwining
between the symmetry and the Galois actions is also independent of
the state. This no longer holds in the non-abelian case, because
of the presence of inner automorphisms of the symmetry group $S$.

\smallskip

\n Moreover, in the $\GL_2$ case, the action of $S$ on the
extremal KMS$_\infty$ states is not transitive, and the
corresponding invariant of the orbit of a state $\varphi$ under
$S$ is the subfield $F_\varphi\subset \C$, which is the
specialization of the modular field given by evaluation at the
point in the upper half plane parameterizing the ground state
$\varphi$.

\smallskip

\n Another important nuance is that the algebra $A$ is no longer
unital while ${\mathcal A}_\Q$ is a subalgebra of the algebra of
unbounded multipliers of $A$. Just as an ordinary function need
not be bounded to be integrable, so states can be evaluated on
unbounded multipliers. In our case, the rational subalgebra
${\mathcal A}_\Q$ is not self-adjoint.

\smallskip

\n Finally the action of the symmetry group on the ground states
is obtained via the action on states at positive temperature.
Given a ground state, one warms it up below the critical
temperature and acts on it by endomorphisms. When taking the limit
to zero temperature of the resulting state, one obtains the
corresponding transformed ground state. In our framework, the
correct notion of ground states is given by a stronger form of the
KMS$_\infty$ condition, where we also require that these are weak
limits of KMS$_\beta$ states for $\beta \to \infty$.

\smallskip

\n We then consider the ``dual'' system of the $\GL_2$-system,
which describes the space of commensurability classes of
2-dimensional $\Q$-lattices (not up to scaling). The corresponding
algebra is closely related to the modular Hecke algebra of
\cite{CM10}. As in the 1-dimensional case, where the corresponding
space is compactified by removing the non-zero condition for the
archimedean component of the ad\`ele, the compactification of the
two-dimensional system amounts to replacing the archimedean
component $\GL_2(\R)$ with matrices $M_2(\R)$. This corresponds to
the noncommutative compactification of modular curves considered
in \cite{ManMar}. In terms of $\Q$-lattices this corresponds to
degenerations to pseudo-lattices, as in \cite{Man1}.

\smallskip

\n It is desirable to have a concrete physical (experimental)
system realizing the BC symmetry breaking phenomenon (as suggested
in \cite{planat}). In fact, we shall show that the explicit
presentation of the BC algebra not only exhibits a strong analogy
with phase states, as in the theory of optical coherence, but it
also involves an action on them of a discrete scaling group,
acting by integral multiplication of frequencies.

\bigskip

\n {\bf Acknowledgements.} We are very grateful to Niranjan
Ramachandran for many extremely useful conversations on class
field theory and KMS states, that motivated the $\GL_2$ system
described here, whose relation to the theory of complex
multiplication is being investigated in \cite{CMR}. We thank
Marcelo Laca for giving us an extensive update on the further
developments on \cite{BC}. We benefited from visits of the first
author to MPI and of the second author to IHES and we thank both
institutions for their hospitality. The second author is partially
supported by a Sofja Kovalevskaya Award of the Humboldt Foundation
and the German Government.

\bigskip
\section{Quantum Statistical Mechanics}\label{stat}

\noindent In classical statistical mechanics a state is a probability
measure $\mu$ on the phase space that assigns to each observable
$f$ an expectation value, in the form of an average
\begin{equation}\label{average-stat} \int f \, d\mu. \end{equation}

\smallskip

\noindent In particular for a Hamiltonian system, the Gibbs canonical
ensemble is a measure
defined in terms of the Hamiltonian and the
symplectic structure on the phase space. It depends on a
parameter $\beta$, which is an inverse temperature, $\beta=1/kT$
with $k$ the Boltzmann constant. The Gibbs measure is given by
\begin{equation}\label{Gibbs} d\mu_{G} = \frac{1}{Z} e^{-\beta H}
d\mu_{Liouville}, \end{equation} normalized by $Z= \int e^{-\beta
H} d\mu_{Liouville}$.

\smallskip

\n When passing to infinitely many degrees of freedom, where the
interesting phenomena of phase transitions and symmetry breaking
happen, the definition of the Gibbs states becomes more involved (\cf
\cite{ruelle}). In the quantum mechanical framework, the analog of the
Gibbs condition is given by the KMS condition at inverse temperature
$\beta$ (\cite{hhw}). This is simpler in formulation than its classical
counterpart, as it relies only on the involutive algebra $A$ of
observables and its time evolution $\sigma_t\in \Aut(A)$, and does
not involve any additional structure like the symplectic structure or
the approximation by regions of finite volume.

\smallskip

\noindent In fact, quantum mechanically, the observables form a
$C^\ast$-algebra $A$, the Hamiltonian is the
infinitesimal generator of the (pointwise norm continuous)
one parameter group of automorphisms
$\sigma_t \in \Aut(A)$, and the analog of a probability
measure, assigning to every observable a certain average, is given
by a {\em state}.

\smallskip

\begin{defn}\label{defstate}
 A {\it state} on a $C^\ast$-algebra
 $A$ is a linear form on ${ A}$
such that
\begin{equation}
\label{eq2}
\varphi (1) = 1 \, , \quad \varphi (a^* a) \geqq 0 \qquad \forall \, a \in
{A} \, .
\end{equation}
\end{defn}

\smallskip

\n When the $C^\ast$-algebra
 $A$ is non unital the condition $\varphi (1) = 1 $
is replaced by $||\varphi||=1$ where
\begin{equation}
\label{eq1}
||\varphi||:= {\rm sup}_{x\in A,||x||\leq 1}|\varphi(x)| \,.
\end{equation}
Such states are restrictions of states on the unital
$C^\ast$-algebra $\tilde{A}$ obtained by adjoining a unit.

\smallskip

\n The evaluation $\varphi (a)$ gives the expectation value
of the observable $a$ in the statistical state
$\varphi$. The Gibbs relation between a thermal state at
inverse temperature $\beta = \frac{1}{kT}$ and the time evolution
\begin{equation}
\label{eq3}
\sigma_t \in {\rm Aut} ({ A})
\end{equation}
is encoded by the KMS condition (\cite{hhw}) which reads
\begin{equation}
\label{eq4}
\forall \, a,b \in { A} \, , \ \exists \, F \
\hbox{bounded holomorphic in the strip}\; \{z\, |\,{\rm Im}\,z\in [0,\beta]\}
\end{equation}
$$ F(t) = \varphi (a \, \sigma_t (b)) \ \ \ \ \  F(t+i\beta) = \varphi
(\sigma_t (b) a) \ \ \  \forall t\in {\mathbb R}. $$

\smallskip

\n In the case of a system with finitely many quantum degrees of
freedom, the algebra of observables is the algebra of operators in a
Hilbert space $\cH$ and the time evolution is given by
$\sigma_t(a)=e^{itH} \, a \, e^{-itH}$, where $H$ is a positive
self-adjoint operator such that $\exp(-\beta H)$ is trace class for
any $\beta >0$. For such a system, the analog of \eqref{Gibbs} is
\begin{equation}\label{GibbsQ}
\varphi(a) = \frac{1}{Z} \, \Tr \left(a\, e^{-\beta H} \right) \ \ \ \
\forall a\in A,
\end{equation}
with the normalization factor $Z=\Tr(\exp(-\beta H))$.
It is easy to see that \eqref{GibbsQ} satisfies the KMS condition
\eqref{eq4} at inverse temperature $\beta$.

\smallskip

\begin{figure}
\begin{center}
\epsfig{file=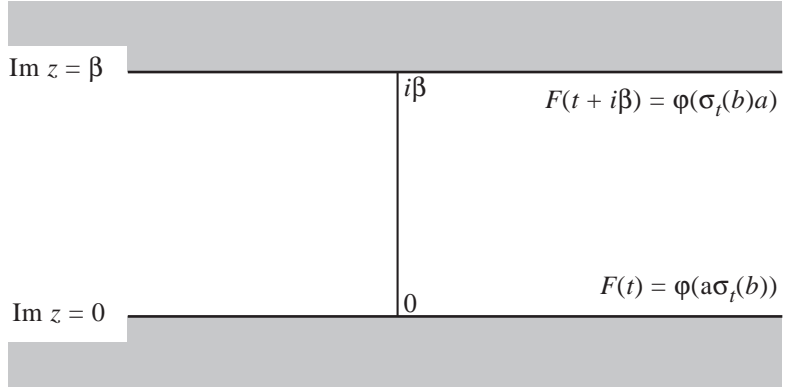}
\end{center}
\caption{The KMS condition.
\label{KMSfig}}
\end{figure}

\medskip

\n In the nonunital case, the KMS condition is defined in the same
way by (\ref{eq4}). Let $M(A)$ be the multiplier algebra of $A$ and let $B\subset
M(A)$ be the $C^*$-subalgebra of elements $x\in M(A)$ such that
$t\mapsto \sigma_t(x)$ is norm continuous.

\begin{prop}\label{KMSnonunital}
Any state $\varphi$ on $A$ admits a canonical extension to a state
still noted $\varphi$ on the multiplier algebra $M(A)$ of $A$. The
canonical extension of a KMS state still satisfies the KMS
condition on $B$.
\end{prop}

\n {\em Proof.} For the first statement we refer to
\cite{pedersen}. The proof of the second statement illustrates a
general density argument, where the continuity of $t\mapsto
\sigma_t(x)$ is used to control the uniform continuity in the
closed strip, in order to apply the Montel theorem of normal
families. Indeed, by weak density of $A$ in $M(A)$, one obtains a
sequence of holomorphic functions, but one only controls their
uniform continuity on smooth elements of $B$. $\Box$

\smallskip

\n As we shall see, it will also be useful to extend whenever
possible the integration provided by a state to unbounded
multipliers of $A$.

\smallskip

\noindent In the unital case, for any given value of $\beta$, the
set $\Sigma_{\beta}$ of KMS$_\beta$ states on ${A}$ forms a convex
compact Choquet simplex (possibly empty and in general infinite
dimensional). In the nonunital case, given a $\sigma_t$-invariant
subalgebra $C$ of $B$, the set $\Sigma_\beta(C)$ of KMS$_\beta$
states on $C$ should be viewed as a compactification of the set of
KMS$_\beta$ states on $A$. The restriction from $C$ to $A$ maps
$\Sigma_\beta(C)$ to KMS$_\beta$ positive linear forms on $A$ of
norm less than or equal to one (quasi-states).

\smallskip

\n The typical pattern for a system with a single phase transition
is that this simplex consists of a single point for $\beta \leqq
\beta_c$ \ie when the temperature is larger than the critical
temperature $T_c$, and is non-trivial (of some higher dimension in
general) when the temperature lowers. Systems can exhibit a more
complex pattern of multiple phase transitions, where no KMS state
exists above a certain temperature. The $\GL_2$ system, which is
the main object of study in this paper, will actually exhibit this
more elaborate behavior.

\smallskip

\n We refer to the books (\cite{brarob}, \cite{haag}) for the
general discussion of KMS states and phase transitions. The main
technical point is that for finite $\beta$ a $\beta$-KMS state is
extremal iff the corresponding GNS representation is factorial.
The decomposition into extremal $\beta$-KMS states is then the
primary decomposition for a given $\beta$-KMS state.

\smallskip

\n At $0$ temperature ($\beta = \infty$) the interesting notion is
that of weak limit of $\beta$-KMS states for $\beta \rightarrow
\infty$. It is true that such states satisfy a weak form of the
KMS condition. This can be formulated by saying that, for all
$a,b\in A$, the function
$$ F(t)= \varphi(a\, \sigma_t(b)) $$
extends to a bounded holomorphic function in the upper half plane
$\Hb$. This implies that, in the Hilbert space of the GNS
representation of $\varphi$ (\ie the completion of $A$ in the inner
product $\varphi(a^*b)$), the generator $H$ of the one-parameter group
$\sigma_t$ is a positive operator (positive energy condition).
However, this condition is too weak in general to be
interesting, as one sees by taking the trivial evolution
($\sigma_t = {\rm id} \qqq t\in \R$). In this case
any state fulfills it, while weak limits of
$\beta$-KMS states are automatically tracial states.
Thus, we shall define  $\Sigma_{\beta=\infty}$
as the set of weak limit points of the sets $\Sigma_{\beta}$ of
$\beta$-KMS states for $\beta \rightarrow \infty$.

\smallskip

\n The framework for spontaneous symmetry breaking (\cite{haag})
involves a
(compact) group of automorphisms $G \subset {\rm Aut}
({ A})$ of ${ A}$
commuting with the time evolution,
\begin{equation}
\label{eq5}
\sigma_t \, \alpha_g = \alpha_g \, \sigma_t \qquad \forall \, g \in G
\, , \ t  \in {\mathbb R} \, .
\end{equation}
The group $G$ is the symmetry group of the system, and the choice
of an equilibrium state $\varphi$ may break it to a smaller subgroup
given by the isotropy group of $\varphi$
\begin{equation}
\label{eq6}
G_\varphi=
\{ g \in G \, , \ g \varphi = \varphi \}.
\end{equation}
The group $G$ acts on $\Sigma_{\beta}$ for any $\beta$, hence on its extreme
points ${\mathcal E} (\Sigma_{\beta}) = {\mathcal E}_{\beta}$.

\n The  unitary group $\Uc$ of the
fixed point algebra of $\sigma_t$ acts by inner
automorphisms of the dynamical system $({A},\sigma_t)$: for $u\in \Uc$,
$$
({\rm Ad}u)\,(x):=\,u\,x\,u^\ast \qqq x\in {A}\,.
$$
These {\em inner} automorphisms of $({A},\sigma_t)$
act trivially on KMS$_\beta$ states,
as one checks using the KMS
condition. This gives us the freedom
to wipe out the group
${\rm Int}({A},\sigma_t)$ of inner symmetries
and to define an action {\em modulo inner}
of a group $G$
on the system $({A},\sigma_t)$
as a map
$$
\a:G\to {\rm Aut}({A},\sigma_t)
$$
fulfilling the condition
$$
\a(g_1g_2)\,\a(g_2)^{-1}\,\a(g_1)^{-1}\in {\rm Int}({A},\sigma_t)
\qqq g_j \in G\,.
$$
Such an action gives an action of the group $G$
on the set $\Sigma_{\beta}$ of KMS$_\beta$ states since the
ambiguity coming from ${\rm Int}({A},\sigma_t)$
disappears in the action on $\Sigma_{\beta}$.
In fact there is one more generalization of the above obvious
notion of symmetries that we shall crucially need later --
it involves actions by endomorphisms.
This type of symmetry plays a key role in the theory
of superselection sectors developed by
Doplicher-Haag-Roberts (cf.\cite{haag}, Chapter IV).

\begin{defn} \label{end} An {\em endomorphism} $\rho$ of the dynamical system
$({A},\sigma_t)$
is a $\ast$-homomorphism $\rho:{A} \to {A}$
commuting with $\sigma_t$.
\end{defn}
It follows then that $\rho(1)=e$ is an idempotent fixed by $\sigma_t$.
Given a KMS$_\beta$ state $\varphi$ the equality
$$
\rho^\ast(\varphi):= Z^{-1}\,\varphi \circ \rho \,,\quad Z=\phi(e)
$$
gives a KMS$_\beta$ state, provided that $\varphi(e)\neq 0$.
Exactly as above for unitaries, consider an isometry
$$
u\in {A}\,,\quad u^\ast\,u=1
$$
which is an eigenvector for $\sigma_t$, \ie that
fulfills, for some $\lambda \in\R^\ast_+$ ($\lambda\geq 1$), the condition
$$
\sigma_t(u)=\, \lambda^{it} \,u\qqq t\in \R\, .
$$
Then $u$ defines an {\em inner} endomorphism ${\rm Ad}u$
of the dynamical system
$({A},\sigma_t)$
by the equality
$$
({\rm Ad}u)\,(x):=\,u\,x\,u^\ast \qqq x\in {A}\,,
$$
and one obtains the following.

\smallskip

\begin{prop} \label{endo1} The  inner endomorphisms
of the dynamical system
$({A},\sigma_t)$ act trivially on the set of KMS$_\beta$ states,
$$
({\rm Ad}u)^\ast(\varphi)=\varphi \qqq \varphi \in \Sigma_{\beta}\,.
$$

\end{prop}

\noindent  {\it Proof.} The KMS$_\beta$ condition shows that
$\varphi(u\,u^\ast)=\lambda^{-\beta}>0$ so that
$({\rm Ad}u)^\ast(\varphi)$ is well defined. The same
KMS$_\beta$ condition applied to the pair $(x\,u^\ast,\,u)$
shows that  $({\rm Ad}u)^\ast(\varphi)=\varphi$. $\Box$

\smallskip

\n At $0$ temperature ($\beta = \infty$) it is no longer
true that the endomorphisms act directly on  the set
$\Sigma_\infty$ of KMS$_\infty$ states, but one can
use their action on KMS$_\beta$-states together
with the ``warming up'' operation. This is defined as the
map
\begin{equation}\label{warming}
W_\beta(\varphi)(x)=\,Z^{-1}\,
{\rm Trace }(\pi(x)\,e^{-\beta \,H})\qqq
x\in A\,,
\end{equation}
where $H$ is the positive energy Hamiltonian, implementing the time
evolution in the representation $\pi$ associated to the KMS$_\infty$
state $\varphi$ and $Z={\rm Trace }(\,e^{-\beta \,H})$.
Typically, $W_\beta$ gives a bijection
$$ W_\beta: \Sigma_\infty \to \Sigma_\beta, $$
for $\beta$ larger than critical. Using the bijection $W_\beta$, one
can transfer the action back to zero temperature states.

\smallskip

\n Another property of KMS states that we shall need later is the
following functoriality. Namely, besides the obvious functoriality
under pullback, discussed above, KMS states push forward under
equivariant surjections, modulo normalization.

\begin{prop} \label{image} Let $({A},\sigma_t)$
be a $C^\ast$-dynamical system ($A$ separable)
and $J$ a norm closed two sided ideal of $A$
globally invariant under $\sigma_t$.
Let $u_n$ be a quasi central approximate
unit for $J$. For any KMS$_\beta$-state $\varphi$
on $({A},\sigma_t)$
the  following sequence converges and defines a KMS$_\beta$
positive linear form
on $(A/J,\sigma_t)$,
$$
\psi(x)=\lim_{n \rightarrow \infty}\varphi((1-u_n)\,x) \qqq x\in A\,.
$$
\end{prop}

\noindent  {\it Proof.} Let $A''$ be the double dual of $A$
and $p\in A''$ the central open projection corresponding
to the ideal $J$ (cf. \cite{pedersen}). By construction the
$u_n$ converge weakly to $p$ (cf. \cite{pedersen} 3.12.14)
so the convergence follows as well as the positivity of $\psi$.
By construction $\psi$ vanishes on $J$. To get the KMS$_\beta$
condition one applies (\ref{eq4}) with $a=(1-u_n)\,x$, $b=y$
where $y$ is a smooth element in $A$. Then one gets a
bounded uniformly continuous sequence $F_n(z)$ of holomorphic
functions in the strip $\{z\, |\,{\rm Im}\,z\in [0,\beta]\}$
with
$$ F_n(t) = \varphi ((1-u_n)\,x \, \sigma_t (y)) \ \ \ \ \
F_n(t+i\beta) = \varphi
(\sigma_t (y) (1-u_n)\,x) \ \ \  \forall t\in {\mathbb R}. $$
Using the Montel theorem on normal families and the quasi-central
property of $u_n$ one gets the KMS$_\beta$
condition for $\psi$. $\Box$

\medskip
\section{$\Q^{ab}$ and KMS states}

\n We shall now describe an explicit system (\cf \cite{bos-con-CR},
\cite{BC}) that will make contact between the general framework above
and arithmetic. The algebra ${\mathcal A}$ of this system is defined
over the rationals,
\begin{equation}
\label{eq8}
{\mathcal A} = {\mathcal A}_{\mathbb Q} \otimes_{\mathbb Q} {\mathbb
C}\, ,
\end{equation}
where ${\mathcal A}_{\mathbb Q}$ is a ${\mathbb Q}$-algebra and is
of countable (infinite) dimension as a vector space over ${\mathbb
Q}$. The algebra ${\mathcal A}$ has a $C^*$-completion $A$ and a
natural time evolution $\sigma_t$.

\smallskip

\n To any vacuum state $\varphi \in {\mathcal E}_{\infty}$ for
$(A,\sigma_t)$ we attach the
${\mathbb Q}$-vector space of complex numbers,
\begin{equation}
\label{eq9}
V_\varphi:=\{ \varphi (a) \, ; \ a \in {\mathcal A}_{\mathbb Q} \}
\end{equation}
that is of countable dimension over ${\mathbb Q}$. It turns out that
$V_\varphi$  is included in algebraic numbers, so that
one can act on these numbers by the Galois group
\begin{equation}
\label{eq10}
\Gal \ (\overline{\mathbb Q} / {\mathbb Q}) \, .
\end{equation}

\smallskip

\n The symmetry group $G$ is the inverse limit with
the profinite topology
\begin{equation}\label{invlimGL1}
G = \hat\Z^* =\varprojlim_n \GL_1(\Z/n\Z) .
\end{equation}
This can also be described as  the quotient
of the id\`ele class group of ${\mathbb Q}$ by the connected component
of the identity,
\begin{equation}
\label{ideleCl}
G = {\rm GL}_1 ({\mathbb Q}) \backslash {\rm GL}_1 ({\mathbb
A}) / {\mathbb R}^*_+ = C_\Q/D_\Q\, .\end{equation}
Here $\A=\A_\Q$ denotes the ad\`eles of $\Q$, namely
$\A=\A_f\times \R$, with $\A_f=\hat\Z\otimes \Q$.

\smallskip

\n The following amazing fact holds:
\begin{equation}
\label{Galpos}
\hbox{For any $\varphi \in {\mathcal E}_{\infty}$ and any $\gamma \in
\Gal  (\overline{\mathbb Q} / {\mathbb Q})$, the composition}
\end{equation}
$$
\gamma \circ \varphi \ \hbox{defined on ${\mathcal A}_{\mathbb Q}$
does extend to a {\it state} on ${\mathcal A}$.}
$$

\medskip

\n What is ``unreasonable'' in this property defining
``fabulous'' states is that, though elements
\begin{equation}
\label{eq13}
\gamma \in \Gal (\overline{\mathbb Q} / {\mathbb Q})
\end{equation}
extend to automorphisms of $\C$, these are extremely discontinuous and
not even Lebesgue measurable (except for $z\mapsto \bar z$), and
certainly do not preserve positivity.

\smallskip

\n It follows from \eqref{Galpos} that the composition $\varphi \mapsto
\gamma \circ \varphi$ defines uniquely an action
of ${\rm Gal}
(\overline{\mathbb Q} / {\mathbb Q})$ on ${\mathcal E}_{\infty}$ and
the equation
\begin{equation}
\label{eq14}
\gamma \circ \varphi = \varphi \circ g
\end{equation}
gives a relation between Galois automorphisms and elements of $G$,
\ie id\`ele classes \eqref{ideleCl}, which is in fact the class field
theory isomorphism $C_\Q/D_\Q\cong \Gal(\Q^{ab}/\Q)$.

\medskip

\n Let us now concretely describe our system, consisting of the
algebra ${\mathcal A}$ (defined over ${\mathbb Q}$) and of the time
evolution $\sigma_t$.

\smallskip

\n The main conceptual steps involved in the construction of this
algebra are:
\begin{itemize}
\item The construction, due to Hecke, of convolution algebras
associated to double cosets on algebraic groups over the rational
numbers;
\item The existence of a canonical time evolution on a von Neumann
algebra.
\end{itemize}

\smallskip

\n More concretely, while Hecke was considering the case of $\GL_2$,
where Hecke operators appear in the convolution algebra associated to
the almost normal subgroup $\GL_2(\Z)\subset \GL_2(\Q)$, the
BC system arises from the Hecke algebra associated to the
corresponding pair of parabolic subgroups.

\smallskip

\n Indeed, let $P$ be the
algebraic group ``$ax+b$'', \ie the functor which to any abelian ring $R$
assigns the group $P_R$ of 2 by 2 matrices over $R$ of the form
\begin{equation}
\label{eq15}
P_R = \left\{ \left[ \begin{matrix} 1 &b \\ 0 &a \end{matrix} \right] \, ; \
a,b \in R \, , \ a \ \hbox{invertible} \right\} \, .
\end{equation}
By construction $P_{\mathbb Z}^+ \subset P_{\mathbb Q}^+$ is an
inclusion $\Gamma_0 \subset \Gamma$ of countable groups, where
$P_R^+$ denotes the restriction to $a>0$. This inclusion fulfills
the following commensurability condition:
\begin{equation}
\label{orbleft}
\hbox{The orbits of the left action of $\Gamma_0$ on $\Gamma / \Gamma_0$ are
all {\it finite}.}
\end{equation}
For obvious reasons the same holds for orbits of $\Gamma_0$ acting on
the right  on $\Gamma_0 \backslash \Gamma$.

\smallskip

\n The Hecke algebra ${\mathcal A}_{\mathbb Q} = {\mathcal H}_{\mathbb
Q} (\Gamma, \Gamma_0)$ is by definition the convolution algebra of
functions of finite support
\begin{equation}
\label{eq17}
f : \Gamma_0 \backslash \Gamma \rightarrow {\mathbb Q}\, ,
\end{equation}
which fulfill the $\Gamma_0$-invariance condition
\begin{equation}
\label{eq18}
f (\gamma \gamma_0) = f(\gamma) \qquad \forall \, \gamma \in \Gamma \, , \
\gamma_0 \in \Gamma_0
\end{equation}
so that $f$ is really defined on $\Gamma_0 \backslash \Gamma /
\Gamma_0$. The convolution product is then given by
\begin{equation}
\label{convhecke}
(f_1 * f_2)(\gamma) = \sum_{\Gamma_0 \backslash \Gamma} f_1
(\gamma
\gamma_1^{-1}) f_2 (\gamma_1) \, .
\end{equation}

\smallskip

\n The time evolution appears from the analysis of the {\em regular
representation} of the pair  $(\Gamma, \Gamma_0)$. It is trivial
when $\Gamma_0$ is normal, or in the original case of Hecke, but it
becomes interesting in the parabolic case, due to the lack of
unimodularity of the parabolic group, as will become clear in the
following.

\smallskip

\n The regular representation
\begin{equation}
\label{eq22}
(\pi (f) \xi)(\gamma) = \sum_{\Gamma_0 \backslash \Gamma} f (\gamma
\gamma_1^{-1}) \xi (\gamma_1)
\end{equation}
in the Hilbert space
\begin{equation}
\label{eq21}
{\mathcal H} = \ell^2 (\Gamma_0 \backslash \Gamma)
\end{equation}
extends to the complexification
\begin{equation}
\label{eq20}
{\mathcal A}_{\mathbb C} = {\mathcal A}_{\mathbb Q} \otimes_{\mathbb Q}
{\mathbb C}
\end{equation}
of the above algebra, which
inherits from this representation the involution $a \mapsto
a^*$, uniquely defined so that $\pi (a^*) = \pi (a)^*$ (the Hilbert space
adjoint), namely
\begin{equation}
\label{eq24}
f^* (\gamma) := \overline{f(\gamma^{-1})} \qquad \forall \, \gamma \in
\Gamma_0 \backslash \Gamma / \Gamma_0.
\end{equation}

\smallskip

\n It happens that the time evolution (\cf \cite{TT}) of the von
Neumann algebra generated by $\sA$ in the regular representation
restricts to the dense subalgebra $\sA$. This implies that
there is a uniquely determined time evolution $\sigma_t\in \Aut(\sA)$,
such that the state $\varphi_1$ given by
\begin{equation}
\label{eq23}
\varphi_1 (f) = \langle \pi (f) \varepsilon_e , \varepsilon_e \rangle
\end{equation}
is a KMS$_1$ state
\ie a KMS state at inverse temperature
$\beta=1$. Here $\varepsilon_e$ is the
cyclic and separating vector for the
regular representation given by the
left coset $\{ \Gamma_0 \} \in \Gamma_0 \backslash
\Gamma$.

\smallskip

\n Explicitly, one gets the following formula for the time evolution:
\begin{equation}
\label{eq25}
\sigma_t (f)(\gamma) = \left( \frac{L(\gamma)}{R(\gamma)} \right)^{-it}
f(\gamma) \qquad \forall \gamma \in \Gamma_0 \backslash \Gamma /
\Gamma_0 \, ,
\end{equation}
where the integer valued functions $L$ and $R$ on the double coset space are
given respectively by
\begin{equation}
\label{eq26}
L(\gamma) = \ \hbox{Cardinality of left $\Gamma_0$ orbit of $\gamma$
in $\Gamma
/ \Gamma_0$}\,,\quad R(\gamma) = L(\gamma^{-1})\, .
\end{equation}

\smallskip

\n Besides the conceptual description given above, the algebra
$\cA_\Q$ also has a useful explicit presentation in terms of
generators and relations (\cf \cite{BC} \S 4, Prop.18). We recall it
here, in the slightly simplified version of \cite{Laca}, Prop.24.

\smallskip

\begin{prop}\label{presentation}
The algebra $\cA_\Q$ is generated by elements $\mu_n$, $n\in
\N^\times$ and $e(r)$, for $r\in \Q/\Z$, satisfying the relations
\begin{itemize}
\item $\mu_n^*\mu_n =1$, for all $n\in \N^\times$,
\item $ \mu_k \mu_n=\mu_{kn} $, for all $k,n\in \N^\times$,
\item $e(0)=1$, $e(r)^*=e(-r)$, and $e(r)e(s)=e(r+s)$, for all $r,s\in
\Q/\Z$,
\item For all $n\in \N^\times$ and all $r\in \Q/\Z$,
\begin{equation}\label{scaling}
 \mu_n\, e(r)\, \mu_n^* = \frac{1}{n} \sum_{ns=r} e(s).
\end{equation}
\end{itemize}
\end{prop}

\smallskip

\n In this form the time evolution preserves pointwise the subalgebra
$R_\Q=\Q[\Q/\Z]$ generated by the $e(r)$ and acts on the $\mu_n$'s as
$$ \sigma_t(\mu_n)=n^{it} \, \mu_n. $$

\medskip

\n The Hecke algebra considered above admits an automorphism
$\a,\,\a^2=1$ whose fixed point algebra is the Hecke
algebra of the pair $P_{\mathbb Z} \subset P_{\mathbb Q}$.
The latter admits an
equivalent
description\footnote{This procedure holds more generally (\cf
\cite{schli1} \cite{schli2}) for arbitrary
almost normal inclusions $(\Gamma_0, \Gamma)$.}, from the pair
$$ (P_{\RR}, P_{{\mathbb A}_f}), $$
where $\RR$ is the maximal compact subring of the ring of finite
ad\`eles
\begin{equation}
\label{eq28}
{\mathbb A}_f = \prod_{\text{res}} {\mathbb Q}_p \, .
\end{equation}
This ad\`elic description displays, as a natural symmetry group, the
quotient $G$ of the
id\`ele class group of ${\mathbb Q}$ by the connected component of
identity \eqref{ideleCl}.

\medskip

\n Let $\overline{\Q}$ be an algebraic closure of $\Q$ and
$\Q^{ab}\subset \overline{\Q}$ be the maximal abelian extension of
$\Q$. Let $r\mapsto \zeta_r$ be a (non-canonical) isomorphism of $\Q/\Z$ with
the multiplicative group of roots of unity inside $\Q^{ab}$.

\smallskip

\n We can now state the basic result that gives content to the
relation between phase transition and arithmetic (BC \cite{BC}):

\begin{thm}\label{BCthm}
\begin{enumerate}
\item  For $0 < \beta \leqq 1$ there exists a unique ${\rm KMS}_{\beta}$
state $\varphi_{\beta}$ for the above system. Its restriction to
$R_\Q=\Q[{\mathbb Q} / {\mathbb Z}] \subset {\mathcal A}$ is given by
\begin{equation}\label{BC-KMS01}
\varphi_{\beta} \left( e(a/b) \right) = b^{-\beta}
\prod_{p \, {\rm prime,} \, p \mid b} \left( \frac{1 - p^{\beta -
1}}{1 - p^{-1}}
\right) \, .
\end{equation}
\item  For $\beta > 1$ the extreme ${\rm KMS}_{\beta}$ states are
parameterized by embeddings $\rho \,:{\mathbb Q}^{ab} \to {\mathbb C}$ and
\begin{equation}\label{BC-KMS1}
\varphi_{\beta , \rho} \left( e(a/b) \right) = Z (\beta)^{-1}
\sum_{n=1}^{\infty} n^{-\beta} \rho \left( \zeta^n_{a/b} \right)
\, ,
\end{equation}
where the partition function $Z(\beta)=\zeta(\beta)$ is the
Riemann zeta function.
\item  For $\beta = \infty$, the Galois group ${\rm Gal}
(\overline{\mathbb Q} / {\mathbb Q})$ acts by  composition on
${\mathcal E}_{\infty}$. The action factors through the abelianization
${\rm Gal}({\mathbb Q}^{ab} / {\mathbb Q})$, and the class field theory
isomorphism $\theta\,: G\to  {\rm Gal}({\mathbb Q}^{ab} / {\mathbb Q})$
intertwines the actions,
$$
\alpha \circ \varphi = \varphi \circ \theta^{-1}(\alpha)\,, \qquad \alpha
\in {\rm Gal} ({\mathbb Q}^{ab} / {\mathbb Q}) \, .
$$
\end{enumerate}
\end{thm}

\bigskip
\section{Further Developments}\label{paperslist}

\noindent The main theorem of class field theory provides a
classification of finite abelian extensions of a local or global
field $K$ in terms of subgroups of a locally compact abelian group
canonically associated to the field. This is the multiplicative
group $K^*=\GL_1(K)$ in the local nonarchimedean case, while in the
global case
it is the quotient $C_K/D_K$ of the id\`ele class group $C_K$ by the
connected component of the identity. The construction of the group
$C_K$ is at the origin of the theory of id\`eles and ad\`eles.

\smallskip

\noindent Hilbert's 12th problem can be formulated as the question
of providing an explicit set of generators of the maximal abelian
extension $K^{ab}$ of a number field $K$, inside an algebraic
closure $\bar K$, and of the action of the Galois group
$\Gal(K^{ab}/K)$. The typical example where this is achieved,
which motivated Hilbert's formulation of the explicit class field
theory problem, is the Kronecker--Weber case: the construction of
the maximal abelian extension of $\Q$. In this case the torsion
points of $\C^*$ (roots of unity) generate $\Q^{ab}\subset \C$.

\smallskip

\noindent Remarkably, the only other case for number fields where
this program has been carried out completely is that of imaginary
quadratic fields, where the construction relies on the theory of
elliptic curves with complex multiplication (\cf \eg
\cite{Steven}). Generalizations to other number fields
involve other remarkable problems in number theory like the Stark conjectures.
Recent work of Manin \cite{Man1} \cite{Man2} suggests a close relation
between the real quadratic case and noncommutative geometry.

\smallskip

\n To better appreciate the technical difficulties underlying any attempt
to address the Hilbert 12th problem of explicit class field theory via
the BC approach, in view of the problem of fabulous states that we
shall formulate in \S \ref{fabulous}, we first summarize briefly the state
of the art
(to this moment and to our knowledge) in the study of $C^*$-dynamical
systems with phase transitions associated to number fields.

\smallskip

\n Some progress from the original BC paper followed in
various directions, and some extensions of the
BC construction  to other global fields (number fields and
function fields) were obtained.
Harari and Leichtnam \cite{HaLe} produced a $C^*$-dynamical system
with phase transition for function fields and algebraic number
fields. In the number fields case a localization is used in order to
deal with lack of unique factorization into primes, and this changes
finitely many Euler factors in the zeta function.
The construction is based on the inclusion
${\mathcal O} \rtimes 1 \subset K \rtimes K^*$, where ${\mathcal O}$ is
the ring of integers of $K$. The symmetry group $G$ of \eqref{ideleCl} is
replaced by the group
$$  G = \hat\Oc^*=\GL_1(\hat\Oc), $$
with $\hat\Oc$ the profinite completion of the ring of integers $\Oc$.
Their result on fixed point
algebras of compact group actions shows that the contribution of the
group of units $\mathcal
O^*$ can be factored out in the cases when this group is finite.
There is a group homomorphism $s:G\to
C_K/D_K$, but it is in general neither injective nor surjective,
hence, even in the case of imaginary quadratic fields,
the construction does not capture the action of the Galois group
$\Gal(K^{ab}/K)$, except in the very special class number one case.

\smallskip

\n P.~Cohen gave in \cite{Cohen} a construction of a $C^*$-dynamical system
associated to a number field $K$, which has spontaneous symmetry
breaking and recovers the full Dedekind  zeta function as partition
function. The main point of her approach
is to involve the semigroup of all ideals
rather than just the principal ideals used in other approaches
as the replacement
of the semi-group of positive integers involved in BC. Still, the group of
symmetries is $G = \hat\Oc^*$ and not the desired $C_K/D_K$.

\smallskip

\n Typically, the extensions of the number field $K$ obtained via these
constructions are given by roots of unity, hence they do not recover
the maximal abelian extension.

\smallskip

\n The Hecke algebra of the inclusion $\mathcal O \rtimes 1 \subset K
\rtimes K^*$ for an arbitrary algebraic number field $K$ was also
considered by  Arledge,  Laca, and  Raeburn in \cite{alr}, where they
discuss its structure and representations, but not the problem of
KMS states.

\smallskip

\n Further results on this Hecke algebra
have been announced by Laca and van
Frankenhuysen \cite{lvf}: they obtain
some general results on the structure
and representations for all number fields,
while they analyze the structure of KMS
states only for the class number one
case. In this case, their
announced result is that there are
enough ground states to support a transitive
free action of $\Gal (K^{ab}/K)$ (up to a copy of
$\{\pm 1\}$ for each real embedding).
However, it appears that the construction does
not give embeddings of $K^{ab}$ as actual
values of the ground states on the Hecke
algebra over $K$, hence it does not seem
suitable to treat the class field theory
problem of providing explicit generators
of $K^{ab}$.

\smallskip

\n The structure of the Hecke algebra of the inclusion ${\mathcal O}
\rtimes {\mathcal O}^* \subset K \rtimes K^*$ was further clarified by
Laca and Larsen in \cite{lala}, using a decomposition of the Hecke
algebra of a semidirect  product as the cross product of the Hecke
algebra of an intermediate (smaller) inclusion by an action of a
semigroup.

\smallskip

\n The original BC algebra was also studied in much greater
details in several following papers. It was proved by Brenken in
\cite{bre} and by Laca and Raeburn in \cite{LacaRaeb} that the
BC algebra can be written as a semigroup cross product.
Brenken also discusses the case of
Hecke algebras from number fields of the type considered in
\cite{LacaRaeb,alr}.

\smallskip

\n Laca then re-derived the original BC result from the point
of view of semigroup cross products in \cite{Laca}. This allows for
significant simplifications of the argument in the case of $\beta >1$,
by looking at the conditional expectations and the KMS condition at
the level of the ``predual'' (semigroup) dynamical system.
A further simplification of the original phase transition theorem of
BC was given by Neshveyev in \cite{nes}, via a direct
argument for ergodicity, which implies uniqueness of the KMS states for
$0\leq \beta \leq 1$.

\smallskip

\n The BC algebra can also be realized as a full corner in the
cross product of the finite ad\`{e}les by the multiplicative
rationals, as was shown by Laca in \cite{minau}, by dilating the
semigroup action to a minimal full group action. Laca and Raeburn used
the dilation results of \cite{minau} to calculate explicitly the
primitive (and maximal) ideal spaces of the BC algebra as
well as of the cross product of the full ad\`{e}les by the action of
the multiplicative rationals.

\smallskip

\n Using the cross product description of the BC algebra,
Leichtnam and Nistor computed Hochschild, cyclic, and periodic cyclic
homology groups of the BC algebra, by computing the
corresponding groups for the $C^*$-dynamical system algebras arising
from the action of $\Q^*$ on the ad\`eles of $\Q$. The calculation for
the BC algebra then follows by taking an increasing sequence
of smooth subalgebras and an inductive limit over certain Morita
equivalent subalgebras.

\smallskip

\n Further results related to aspects of the BC construction
and generalizations can be found in \cite{boc-zah}, \cite{glo-wil},
\cite{quasilat}, \cite{twoprime}, \cite{semi}, \cite{tza}.

\medskip
\section{Fabulous States}\label{fabulous}

Given a number field $K$,
we let $\A_K$ denote the ad\`eles of $K$ and
$J_K=\GL_1(\A_K)$ be the group of id\`eles of $K$. We write $C_K$ for
the group of id\`eles classes $C_K=J_K/K^*$ and $D_K$ for the
connected component of the identity in $C_K$.

\medskip

\n If we remain close to the spirit of the Hilbert 12th problem,
we can formulate a general question, aimed at extending the results
of \cite{BC} to other number fields $K$.
 Given a number field $K$, with a choice of an embedding $K\subset
\C$, the ``problem of fabulous states'' consists in constructing a
$C^*$-dynamical system $(A,\sigma_t)$ and an ``arithmetic''
subalgebra ${\mathcal A}$, which satisfy the following properties:
\begin{enumerate}
\item The id\`eles class group $G=C_K/D_K$
acts by symmetries on $(A,\sigma_t)$ preserving the subalgebra
${\mathcal A}$.
\item The states $\varphi \in \sE_\infty$, evaluated on elements of
${\mathcal A}$, satisfy:
\begin{itemize}
\item $\varphi(a)\in \bar K$, the algebraic
closure of $K$ in $\C$;
\item the
elements of $\{ \varphi(a): a\in {\mathcal A} \}$, for $\varphi\in
\sE_\infty$ generate $K^{ab}$.
\end{itemize}
\item The class field theory isomorphism
$$
\theta:C_K/D_K \stackrel{\simeq}{\longrightarrow} \Gal (K^{ab}/K)
$$
intertwines the actions,
\begin{equation}\label{CFTiso}
 \alpha \circ \varphi = \varphi \circ \theta^{-1}(\alpha),
\end{equation}
for all $\alpha \in \Gal (K^{ab}/K)$ and for all $\varphi \in
\sE_\infty$.
\end{enumerate}

\smallskip

\n Notice that, with this formulation, the problem of the construction
of fabulous states is intimately related to Hilbert's 12th
problem. This question will be pursued in \cite{CMR}.

\smallskip

\n
 We shall construct here a system which is
the analog of the BC system for $\GL_2(\Q)$ instead of
$\GL_1(\Q)$. This will extend
 the results of
\cite{BC} to this non-abelian $\GL_2$ case and will exhibit many new
features which have no counterpart in the abelian case. Our construction
 involves the explicit description of the automorphism group of the
modular field,
\cite{Shimura}.
The construction of fabulous states for imaginary quadratic fields,
which will be investigated with N. Ramachandran in \cite{CMR}, involves
specializing the $\GL_2$ system to
a subsystem compatible with complex multiplication in a given imaginary
quadratic field. \\

\n The   construction of the $\GL_2$ system gives
a $C^*$-dynamical system $(A,\sigma_t)$ and
an involutive subalgebra ${\mathcal A}_\Q$ defined over $\Q$,
satisfying the following properties:

\begin{itemize}
\item The quotient group $
S:=\,\Q^\ast\backslash {\rm GL}_2 ({\A_f})
$ of the finite ad\`elic group of $\GL_2$ acts
as symmetries of the dynamical system $(A,\sigma_t)$
preserving the subalgebra ${\mathcal A}_\Q$.
\item For generic $\varphi \in \sE_\infty$,
 the
values $\{ \varphi(a)\in \C: a\in {\mathcal A}_\Q \}$
generate a subfield $F_\varphi \subset \C$
which is an extension of
$\Q$ of transcendence degree $1$.
\item For generic $\varphi \in
\sE_\infty$, there exists an  isomorphism
$$
\theta:S\stackrel{\simeq}{\longrightarrow} \Gal (F_\varphi/\Q)
$$
which intertwines the actions
\begin{equation}\label{CFTiso1}
 \alpha \circ \varphi = \varphi \circ \theta^{-1}(\alpha),
\end{equation}
for all $\alpha \in \Gal (F_\varphi/\Q)$.
\end{itemize}

\n There are a number of important nuances between the abelian
case above and the non-abelian one. For instance, in the abelian case
the field generated by $\varphi({\mathcal A})$ does not depend on the
choice of $\varphi \in
\sE_\infty$ and the isomorphism $\theta$ is also independent of $\varphi$.
This no longer holds in the non-abelian case, as is clear from the
presence of inner automorphisms of the symmetry group $S$. Also, in the latter
case, the action of $S$ on $\sE_\infty$ is not transitive and the
corresponding invariant of the orbit of
$\varphi$ under $S$ is the subfield $F_\varphi\subset \C$.
Another important nuance is that the algebra $A$
is no longer  unital while  ${\mathcal A}_\Q$ is an
 algebra of unbounded  multipliers of $A$. Finally, the symmetries
require the full framework of endomorphisms as explained above in
\S \ref{stat}.

\bigskip
\section{The subalgebra $\cA_\Q$ and Eisenstein Series}\label{eisenstein}

\smallskip

\n In this section we shall recast
the BC algebra in terms of the trigonometric analog of the
Eisenstein series, following the analogy developed by Eisenstein and
Kronecker between
trigonometric and elliptic functions, as outlined by A.Weil in
\cite{WeilEll}.

\n This will be done by first giving a geometric
interpretation in terms of $\Qb$-lattices
of the noncommutative space $X$ whose algebra of continuous
functions $C(X)$ is the BC $C^*$-algebra.
The space $X$ is by construction the quotient of
 the Pontrjagin dual of the abelian group
$\Q/\Z$ by the equivalence relation generated by the action
by multiplication of the
semi-group $\N^\times$.

\n Let
$$
\RR=\,\prod_p \,\Z_p
$$
be the compact ring product of the rings $\Z_p
$ of $p$-adic integers. It is the maximal compact subring of
the locally compact ring of finite ad\`eles
$$
{\mathbb A}_f =
 \prod_{{\rm res}} {\mathbb Q}_p
$$
We recall the following standard fact

\begin{prop}
\begin{itemize}
\item The inclusion $\Q \subset {\mathbb A}_f$ gives
an isomorphism of abelian groups
$$
\Q/\Z = {\mathbb A}_f/\RR  \,.
$$

\item The following map is an isomorphism of compact rings
$$
j:\;\RR \rightarrow
{\rm Hom}(\mathbb{Q}\slash\mathbb{Z}, \,\mathbb{Q}\slash\mathbb{Z})\,,
\quad j(a)(x)= a \,x \qqq x\in {\mathbb A}_f/\RR \qqq a \in \RR\,.
$$

\end{itemize}
\end{prop}

\n We shall use $j$ from now on
to identify $\RR$ with
${\rm Hom}(\mathbb{Q}\slash\mathbb{Z}, \,\mathbb{Q}\slash\mathbb{Z})$.
Note that by construction
${\rm Hom}(\mathbb{Q}\slash\mathbb{Z}, \,\mathbb{Q}\slash\mathbb{Z})$
is endowed with the topology of pointwise convergence.
It is identified with $\varprojlim \Z/N\Z
$ using the restriction to
$N$-torsion elements.

\n For every $r\in \Q/\Z$ one gets a function $e(r)\in C(\RR)$ by,
$$
e(r)( \rho):={\rm exp}2\pi i \rho(r)\,\quad \forall \rho \in
{\rm Hom}(\mathbb{Q}\slash\mathbb{Z}, \,\mathbb{Q}\slash\mathbb{Z})
$$
and this gives the identification of $\RR$ with the Pontrjagin dual of
$\Q/\Z$ and of $C(\RR)$ with the group $C^\ast$-algebra
$C^\ast(\mathbb{Q}\slash\mathbb{Z})$.

\n One can then describe the
BC $C^*$-algebra as the cross product of $C(\RR)$ by the
semigroup action of $\N^\times$ as follows.
 For each integer $n\in \N^\times$ we let $n \RR \subset \RR$
be the range of multiplication by $n$. It is an open and closed subset
of $\RR$ whose characteristic function $\pi_n$ is a projection
$\pi_n \in C(\RR)$. One has by construction
$$
 \pi_n \, \pi_m = \pi_{n \vee m} \,,\quad \forall n\,,\;m \in \N^\times
$$
where $n \vee m$ denotes the lowest common multiple of $n$ and $m$.

\n The semigroup action of $\N^\times$ on $C(\RR)$ corresponds to the
isomorphism
\begin{equation}\label{end1}
\a_n(f)(\rho):= f(n^{-1}\, \rho) \qqq \rho \in n \RR \,.
\end{equation}
of $C(\RR)$ with the reduced algebra
$C(\RR)_{\pi_n}$ of
$C(\RR)$ by the projection $\pi_n$.
In the BC algebra one has
\begin{equation}\label{end2}
\mu_n \, f\, \mu_n^*= \, \a_n (f) \qqq f\in C(\RR).
\end{equation}
There is an equivalent description of the  BC algebra
in terms of the \'etale groupoid $G$
of pairs $(r, \rho)$, where $r\in \Q^*_+$, $\rho \in \RR$ and
$r\, \rho \in \RR$.
The composition in $G$ is given by
\begin{equation}
\label{comp01}
(r_1, \rho_1)\circ (r_2, \rho_2)=\,(r_1\,r_2, \rho_2)\,,
\quad {\rm if}\;r_2\, \rho_2=\,\rho_1\,,
\end{equation}
and the convolution of functions by
\begin{equation}
\label{conv01}
f_1 \ast f_2(r,\rho):=\,\sum \,f_1(r s^{-1},s\,\rho)\,f_2(s,\rho)\,,
\end{equation}
while the adjoint of $f$ is
\begin{equation}
\label{adjoint0}
f^\ast(r,\rho):=\,\overline{f(r^{-1},r\,\rho)}\,.
\end{equation}

\n All of this is implicit in (\cite{BC})
and has been amply described in the subsequent papers mentioned in
\S \ref{paperslist}.
In the description above, $\mu_n$ is given by the function
$\mu_n(r,\rho)$ which vanishes unless $r=n$ and is equal
to $1$ for $r=n$. The time evolution is given by
\begin{equation}
\label{sigmat}
\sigma_t (f)(r, \rho):= r^{it}\,f(r, \rho)\qqq f\in C^\ast(G)\,.
\end{equation}

\n We shall now describe a geometric interpretation of this
groupoid $G$ in terms of commensurability of $\Q$-lattices.
In particular, it will pave the way to the generalization
of the BC system to higher dimensions.
The basic simple geometric objects are
$\Qb$-lattices in $\Rb^n$,
defined as follows.

\begin{defn}\label{Defn-Qlat}
A $\Qb$-lattice in $\Rb^n$ is a
 pair $ \, ( \L , \phi) \,$, with $
\L $ a lattice in $\mathbb{R}^n$, and
 $\displaystyle \phi :  \mathbb{Q}^n\slash\mathbb{Z}^n
\longrightarrow\mathbb{Q}\L \slash \L $
an homomorphism of abelian groups.
\end{defn}

\n Two lattices $\L_j$ in $\Rb^n$ are
commensurable iff their intersection
$\L_1 \cap \L_2$ is of finite index in
$\L_j$. Their sum $\L=\L_1+\L_2$
is then a lattice and, given two
homomorphisms of abelian groups
$\displaystyle \phi_j :  \Q^n\slash\Z^n
\longrightarrow\mathbb{Q}\,\L_j \slash \L_j $,
the difference $\phi_1-\phi_2$ is well
defined modulo $\L=\L_1+\L_2$.

\smallskip

\n Notice that in Definition \ref{Defn-Qlat} the homomorphism
$\phi$, in general, is not an isomorphism.

\begin{defn}\label{Defn-invQlat}
A $\Q$-lattice $(\L,\phi)$ is {\em invertible} if the map $\phi$
is an isomorphism of abelian groups.
\end{defn}

\smallskip

\n We consider a natural equivalence relation on the set of
$\Q$-lattices defined as follows.

\begin{prop} \label{comm}The following
defines an equivalence relation
called {\em commensurability} between
 $\Qb$-lattices:  $ \, ( \L_1 , \phi_1) \,,\;( \L_2 , \phi_2)$
are commensurable iff $ \L_j$ are commensurable and
$\phi_1-\phi_2=0$ modulo $\L=\L_1+\L_2$.
\end{prop}

\noindent  {\it Proof.} Indeed, let $ \, ( \L_j , \phi_j) \,$
be three $\Qb$-lattices and assume commensurability
between the pairs $(1,2)$ and $(2,3)$. Then the lattices
$ \L_j$ are commensurable and are of finite index in
$\L=\L_1+\L_2+\L_3$. One has $\phi_1-\phi_2=0$ modulo $\L$,
$\phi_2-\phi_3=0$ modulo $\L$ and thus
$\phi_1-\phi_3=0$ modulo $\L$. But $\L'=\L_1+\L_3$
is of finite index in $\L$ and thus $\phi_1-\phi_3$
gives a group homomorphism
$$ \Q^n\slash\Z^n
\longrightarrow \,\L \slash \L' $$
which is zero since $\Q^n\slash\Z^n$
is infinitely divisible and
$\L \slash \L' $ is finite. This shows that
$\phi_1-\phi_3=0$ modulo $\L'=\L_1+\L_3$
and hence that the pair $(1,3)$
is commensurable. $\Box$

\smallskip

\no Notice that every $\Qb$-lattice in $\R$ is uniquely of the form
\begin{equation}\label{1dQlat}
 ( \L , \phi) \,
= (\lambda\, \Z,\lambda\,\rho)\,,\quad\lambda>0,
\end{equation}
with $\rho \in {\rm Hom}(\mathbb{Q}\slash\mathbb{Z},
\,\mathbb{Q}\slash\mathbb{Z}) =\RR$.

\smallskip

\begin{prop} \label{groupoid1}
The  map
$$
\g(r,\rho)=((r^{-1}\, \Z,\,\rho)\,,( \Z,\rho))\qqq (r,\rho) \in G\,,
$$
defines an isomorphism of locally compact
\'etale groupoids between $G$ and the
quotient $\Rc\,/\,\R^\ast_+$
of the equivalence relation $\Rc$ of commensurability
on the space of $\Qb$-lattices in $\R$
by the natural scaling action of $\R^\ast_+$.
\end{prop}

\noindent  {\it Proof.} First since $r\,\rho \in \RR$ the
pair $(r^{-1}\, \Z,\,\rho)=r^{-1}( \Z,\,r\,\rho)$
is a $\Q$-lattice and is commensurable to $( \Z,\rho)$.
Thus, the map $\g$ is well defined. Using the identification
\eqref{1dQlat}, we see that
the restriction of $\g$ to the objects $G^{(0)}$
of $G$ is an isomorphism of $\RR$ with the quotient of
the space of $\Qb$-lattices in $\R$
by the natural scaling action of $\R^\ast_+$. The freeness
of this action shows that the quotient $\Rc\,/\,\R^\ast_+$
is still a groupoid, and one has
$$
\g(r_1,\rho_1)\circ \g(r_2,\rho_2)= \,\g(r_1\,r_2,\rho_2)
\quad {\rm if}\;r_2\, \rho_2=\,\rho_1\,.
$$
Finally, up to scaling, every element of $\Rc$ is of
the form
$$
((r^{-1}\, \Z,\,r^{-1}\,\rho')\,,( \Z,\rho))
$$
where both $\rho'$ and $\rho$ are in $\RR$ and $r=\frac{a}{b}\in \Q^\ast_+$.
Moreover since $r^{-1}\,\rho'=\rho$ modulo $\frac{1}{a}\Z$
one gets $a \rho- b \rho'=0$ and $r^{-1}\,\rho'=\rho$.
Thus $\g$ is surjective and is an isomorphism. $\Box$\\

\n This geometric description of the BC algebra
allows us to generate in a natural manner a rational subalgebra
 which will generalize to the two dimensional case.
In particular the algebra $C(\RR)$ can be viewed as the algebra
of homogeneous functions of ``weight $0$" on the space of $\Qb$-lattices
for the natural scaling action of the multiplicative group
$\R^\ast_+$ where weight $k$ means
$$
f(\lambda\, \L , \lambda\,\phi)=\,\lambda^{-k} \,f( \L , \phi)\qqq
\lambda \in\R^\ast_+\,.
$$
We let the function $c(\Lambda)$ be the multiple of
the covolume $|\L|$ of the lattice, specified by
\begin{equation}\label{cnorm}
2\,\pi\,i\,c(\Zb)= 1
\end{equation}
The function $c$ is homogeneous  of
weight $-1$ on the space of $\Qb$-lattices.
For $a \in  \mathbb{Q}\slash\mathbb{Z}\,$, we then define
a function $e_{1,a}$ of weight $0$ by
\begin{equation}\label{epsilon1}
e_{1,a} ( \L , \phi)= c(\L) \sum_{y\in  \Lambda +\phi(a)} y^{-1}
\, ,
\end{equation}
where one uses Eisenstein summation \ie
$\lim_{N \rightarrow \infty}\sum_{-N}^N$ when $\phi(a)\neq 0$
and one lets $e_{1,a} ( \L , \phi)=0$ when $\phi(a)= 0$.

\smallskip

\n The main result of this section is the following

\medskip

\begin{thm}\label{ratiosubalg} \begin{itemize}

\item The $e_{1,a}\,,  a \in  \mathbb{Q}\slash\mathbb{Z}\,$
generate $\Q[\mathbb{Q}\slash\mathbb{Z}]$.

\item  The rational algebra $\cA_\Qb$ is the
 subalgebra of $A=C^\ast(G)$ generated by the
$e_{1,a}\,, \, a \in  \mathbb{Q}\slash\mathbb{Z}\,$
and the
 $\mu_n$, $\mu_n^\ast$.
\end{itemize}
\end{thm}

\medskip

\n We  define more generally for each  weight $k\in \Nb$ and each
$a \in \mathbb{Q}\slash\mathbb{Z}\,$ a function
$\epsilon_{k,a}$ on the space of $\Qb$-lattices in $\Rb$ by
\begin{equation}\label{epsilon0k}
\epsilon_{k,a} ( \L , \phi)= \sum_{y\in  \Lambda +\phi(a)} y^{-k}\, .
\end{equation}
This is well defined  provided $\phi(a)\neq 0$. For $\phi(a)= 0$
we let
\begin{equation}\label{lambda1}
\epsilon_{k,a} ( \L , \phi)= \lambda_k \;c(\Lambda)^{-k},
\end{equation}
where we shall fix the constants $\lambda_k$ below in (\ref{lambdak}).
The function $\epsilon_{k,a}$ has weight $k$ \ie it satisfies the
homogeneity condition
$$
\epsilon_{k,a} (\lambda \,\L, \, \lambda\, \phi )=\lambda^{-k}
\epsilon_{k,a} ( \L , \phi) \qqq \lambda \in\R^\ast_+\,.
$$
When $a= \frac{b}{N}$ the function $\epsilon_{k,a}$ has level $N$ in
that it only
uses the restriction $\phi_N$ of $\phi$ to $N$-torsion points
of $\mathbb{Q}\slash\mathbb{Z}$,
$$\displaystyle \phi_N : \frac{1}{N} \mathbb{Z}\slash \mathbb{Z}
\longrightarrow\frac{1}{N}\L \slash \L .
$$

\smallskip

\n The products
\begin{equation}\label{epsilonk}
e_{k,a}:= c^k \, \epsilon_{k,a}
\end{equation}
are of weight $0$ and satisfy two types of relations.

\smallskip

\n The first relations are multiplicative
and express $e_{k,a}$ as a polynomial in $e_{1,a}$,
\begin{equation}\label{mul-rel}
e_{k,a}=\,P_k( e_{1,a})
\end{equation}
where the $P_k$ are the polynomials with rational coefficients
uniquely determined by the equalities
$$
P_1(u)=\,u\,,\quad P_{k+1}(u)=\frac{1}{k} (u^2-\frac{1}{4})
 \,\partial_uP_k(u)\,.
$$
This follows for $ \phi(a)\notin \L$ from the elementary
formulas for the trigonometric analog of the Eisenstein
series (\cite{WeilEll} Chapter II).
 Since $e_{1,a}( \L , \phi)=0$ is the natural choice for
$ \phi(a)\in \L$, the validity of (\ref{mul-rel})
uniquely dictates the choice
of the normalization constants $\lambda_k$ of (\ref{lambda1}). One gets
\begin{equation}\label{lambdak}
\lambda_k=P_k(0)=\,(2^k-1)\,\g_k\, ,
\end{equation}
where $\g_k = 0$ for odd $k$ and $\displaystyle \g_{2j} =
(-1)^j\,\frac{B_j}{(2j)!}$
with $B_j \in \Qb$ the Bernoulli numbers. Equivalently,
$$
\frac{x}{e^x -1}=1-\frac{x}{2}-\sum_1^\infty\,\g_{2j}\,x^{2j}.
$$

\n One can express the $e_{k,a}$ as $\Q$-linear
combinations of the generators
$e(r)$. We view $e(r)$ as the function
on $\Qb$-lattices which assigns
to $ \, ( \L , \phi) \,
= (\lambda\, \Z,\lambda\,\rho)$, $\lambda>0$,
the value
$$
e(r)( \L , \phi):={\rm exp}2\pi i \rho(r).
$$
One then has

\begin{lem} Let $a\in \Q/\Z$, and $n>0$
with $n\,a=\,0$. Then
\begin{equation}\label{etoe}
e_{1,a}=\,\sum_{k=1}^{n-1}(\frac{k}{n}-\frac{1}{2})\,e(k\,a).
\end{equation}
\end{lem}

\noindent  {\it Proof.} We evaluate both sides on
$ \, ( \L , \phi) \,
= (\lambda\, \Z,\lambda\,\rho)$, $\lambda>0$.
 Both sides only depend on
the restriction $\displaystyle x \mapsto c\,x$
 of $\rho$ to $n$-torsion elements
of $\Q/\Z$
which we  write as multiplication by $c\in \Z/n\Z$.
Let $a=\frac{b}{n}$.
If $b c =0 (n)$ then $\phi(a)=0$ and both sides
vanish since $e(k\,a)( \L , \phi)=
{\rm exp}2\pi i (\frac{k b c}{n})=1$
for all $k$.
If $b c \neq 0 (n)$ then $\phi(a)\neq 0$ and the left side
is $\frac{1}{2}(U+1)/(U-1)$ where
$U={\rm exp}2\pi i \frac{ b c}{n} $, $U^n=1$, $U\neq 1$.
The right hand side is
$$
\,\sum_{k=1}^{n-1}(\frac{k}{n}-\frac{1}{2})\,U^k\, ,
$$
which gives $\frac{1}{2}(U+1)$ after multiplication by $U-1$.
$\Box$\\

\smallskip

\n This last equality shows that $e_{1,a}$ is (one half of)
the Cayley transform of $e(a)$ with care taken where $e(a)-1$
fails to be invertible. In particular while $e(a)$
is unitary, $e_{1,a}$ is skew-adjoint,
$$
e_{1,a}^\ast=\;-\,e_{1,a}\,.
$$
 We say that a $\Qb$-lattice $ \, ( \L , \phi) \,$
 is divisible by an integer $n\in \Nb$
when $\phi_n=0$.
We let $\pi_n$ be the characteristic function of the set
of $\Qb$-lattices divisible by $n$. It corresponds to the
characteristic function of $n\,\RR\subset \RR$.

\n Let $N>0$ and $ \, ( \L , \phi) \,$ a $\Qb$-lattice
with $\phi_N(a)= c\, a$
for $c\in \Z/N\Z$. The order of the kernel of
 $\phi_N$ is $m={\rm gcd}(N,c)$. Also a divisor  $b|N$
divides $ \, ( \L , \phi) \,$ iff it divides $c$.
Thus for any function $f$ on $\N^\ast$ one has
$$
\sum_{b|N}f(b)\,\pi_b ( \L , \phi) = \sum_{b|{\rm gcd}(N,c)}f(b)\, ,
$$
which allows one to express any function of the order
$m={\rm gcd}(N,c)$ of the kernel
of $\phi_N$ in terms of the projections $\pi_b\,,\, b|N$. In order to
obtain the function
$m\mapsto m^j$ we  let
$$
f_j(n):=\sum_{d|n} \mu(d) (n/d)^j
 \, ,
$$
where $\mu$ is the M\"obius function so that
$$
f_j(n)=\,n^j\,
\prod_{p \, {\rm prime,} \, p \mid n}\,(1 - p^{-j})\,.
$$
Notice that $f_1$
is the Euler totient function and that the ratio $f_{-\beta +1}/f_1$
gives the r.h.s. of \eqref{BC-KMS01} in Theorem \ref{BCthm}.

\smallskip

\no The M\"obius inversion formula  gives
\begin{equation}\label{div}
\sum_{b|N}f_j(b)\,\pi_b ( \L , \phi) = m^j \,,\quad m={\rm gcd}(N,c)\,.
\end{equation}
We can now write division relations fulfilled by
the functions \eqref{epsilonk}.

\begin{lem} Let  $N>0$
 then
\begin{equation}\label{div-rel}
\sum_{N\,a = 0} e_{k,a}= \g_k\,\sum_{d|N}((2^k-2)\,f_{1}(d)+N^k
f_{-k+1}(d))\,\pi_d .
\end{equation}
\end{lem}

\noindent  {\it Proof.}
For a given $\Qb$-lattice $ \, ( \L , \phi) \,$
with ${\rm Ker}\,\phi_N$ of order $m |N$, $N= m \, d$, the result
follows from
$$
\sum_{N\,a = 0} \epsilon_{k,a} ( \L , \phi)=m\,
\sum_{y\in \frac{1}{d}\L\backslash \L}y^{-k}\,+m\,
(2^k-1)\,\g_k \,c^{-k}( \L)
=m\,(d^{\,k}\,+2^k-2)\g_k \,c^{-k}
$$
together with (\ref{div}) applied for $j=1$ and $j=1-k$. $\Box$

\smallskip

\n The semigroup action of $\N^\times$ is given on functions
of $\Q$-lattices by the endomorphisms
\begin{equation}\label{end11}
\a_n(f)(\L , \,\phi):= f( n\,\L , \phi) \qqq (\L , \phi) \in \pi_n \,,
\end{equation}
while $\a_n(f)(\L , \,\phi)=0$ outside $\pi_n$. This
semigroup action preserves the rational subalgebra $\Bc_\Qb$
generated by the $e_{1,a}\,,  a \in  \mathbb{Q}\slash\mathbb{Z}\,$,
since one has
\begin{equation}\label{muepsilon}
\a_n( e_{k,a}) =\,\pi_n\, \, e_{k,a/n}\, ,
\end{equation}
(independently of the choice of the solution $b=a/n$
of $n\,b=a$) and we shall now show that the projections
$\pi_n$ belong to $\Bc_\Qb$.\\

\noindent  {\it Proof of Theorem \ref{ratiosubalg}}

\n Using (\ref{div-rel}) one can
express $\pi_n$ as a rational linear combination of the $e_{k,a}$,
with $k$ even, but special care is needed when $n$ is a  power of two.
The coefficient of $\g_k\,\pi_N$ in (\ref{div-rel}), when $N=p^b$
is a prime power, is given by $(2^k-2)(p-1)p^{b-1} - p^b(p^{k-1}-1)$,
which does not vanish unless $p=2$, and is
$-p^{b-1}(2-3p+p^2)$ for $k=2$. Thus, one can express $\pi_N$
as a linear combination of the $e_{2,a}$
by induction on $b$.  For $p=2$, $N=2^b$, $b>1$
the coefficient of $\g_k\,\pi_N$ in (\ref{div-rel}) is zero but the
coefficient
of $\g_k\,\pi_{N/2}$ is $-2^{b-2}(2^k-1)(2^k-2)\neq 0$ for $k$ even.
This allows one to express $\pi_N$ as a linear combination of the $e_{2,a}$
by induction on $b$.
Thus, for instance, $\pi_2$ is given by
$$
\pi_2=\, 3 + 2\, \sum_{4\,a = 0} e_{2,a}.
$$
In general, $\pi_{2^n}$ involves $\displaystyle \sum_{2^{n+1}\,a = 0}
e_{2,a}$.
\smallskip

\n Since for relatively prime integers $n,\,m$ one has
$\pi_{nm}=\pi_n\,\pi_m$,
we see that the algebra $\Bc_\Qb$ generated over $\Qb$
by the $e_{1,a}$ contains all the projections $\pi_n$.
In order to show that $\Bc_\Qb$ contains the $e(r)$ it is enough
to show that for any prime power $N=p^b$ it contains $e(\frac{1}{N})$.
This is proved by induction on $b$. Multiplying
(\ref{div-rel}) by $1-\pi_p$ and using $(1-\pi_p)\,\pi_{p^l}=0$
for $l>0$ we get the equalities
$$
(1-\pi_p)\,\sum_{N\,a = 0} e_{k,a}= (N^k\,+2^k-2)\,\g_k \,(1-\pi_p)\, .
$$
Let then $z(j)=\,(1-\pi_p)\,e_{1,\frac{j}{N}}$. The above relations together
with(\ref{mul-rel})   show
that in the reduced algebra $(\Bc_\Qb)_{1-\pi_p}$ one has, for all $k$,
$$
\sum_{j=1}^{N-1}\,P_k(z(j))=\,(N^k-1)\,\g_k \,.
$$
Thus, for $j\in\{1,..., N-1\}$, the symmetric functions
of the $z(j)$ are fixed rational numbers $\sigma_h$.
In particular $z=z(1)$ fulfills
$$
Q(z)=\,z^{N-1}+\,\sum_1^{N-1} (-1)^h\,\sigma_h\, z^{N-1-h}=\,0
$$
and $\pm \,\frac{1}{2}$ is not a root of this equation,
whose roots are the $\displaystyle \frac{1}{2i}\,{\rm cot}(\frac{\pi\,j}{N})$.
This allows us, using the companion matrix of $Q$, to express the
 Cayley
transform of $2 z$ as a polynomial with  rational coefficients,
$$
\frac{2 z+1}{2 z-1}=\,\sum_0^{N-2} \a_n \,z^n.
$$
One then has
$$
\sum_0^{N-2} \a_n \,z^n=\,(1-\pi_p)\,e(\frac{1}{N}),
$$
where the left-hand side belongs to $\Bc_\Qb$ by construction.
Now $\pi_p\,e(\frac{1}{N})$ is equal to $\a_p( e(\frac{p}{N}))
$. It follows from the induction hypothesis
on $b$, ($N=p^b$), that $e(\frac{p}{N})\in \Bc_\Qb$
and therefore using (\ref{muepsilon}) that $\a_p( e(\frac{p}{N}))
\in  \Bc_\Qb$. Thus, we get $e(\frac{1}{N})\in \Bc_\Qb$
as required. This proves the first part. To get the second
notice that the cross product by $\N^\times$
is obtained by adjoining to the rational group ring of $\Q/\Z$
the isometries $\mu_n$ and their adjoints $\mu_n^\ast$
with the relation
$$
\mu_n\, f\,\mu_n^\ast=\, \a_n(f)\qqq f\in
\Q[\mathbb{Q}\slash\mathbb{Z}]\, ,
$$
which gives the rational algebra $\cA_\Qb$.
$\Box$

\medskip

\n It is not true, however, that the division relations (\ref{div-rel})
combined with the multiplicative relations (\ref{mul-rel})
suffice to present the algebra. In particular there are more elaborate
division relations which we did not need in the above proof.
In order to formulate them, we let for $d|N$, 
 $\pi(N,d)$ be the projection belonging to the algebra
generated by the $\pi_b,\, b|N$, and corresponding to the subset
$$
 {\rm gcd }(N, \,  ( \L , \phi))=\,N/d
$$
so that
$$
\pi(N,d)=\,\pi_{N/d}\,\prod_{k|d} (1-\pi_{k\,N/d}),
$$
where the product is over non trivial divisors $k\neq 1$ of $d$.

\smallskip

\begin{prop}\label{div1} The $e_{k,a}\,, \, a \in
\mathbb{Q}\slash\mathbb{Z}\,$
, $k$ odd,
fulfill for any $x\in  \mathbb{Q}\slash\mathbb{Z}\,$ and any
integer $N$ the relation
$$
\frac{1}{N}\,\sum_{N\,a = 0} e_{k,x+a}=\,
\sum_{d|N}\,\pi(N,d)\,d^{\,k-1}\,e_{k,d\, x}.
$$
\end{prop}

\noindent  {\it Proof.}  To prove this, let $( \L , \phi)$
be such that ${\rm gcd }(N, \,  ( \L , \phi))=\,N/d=m$
and assume by homogeneity that $\Lambda=\Zb$. Then when
$a$ ranges through the $\frac{j}{N},\, j\in\{0,...,N-1\}$,
the $\phi(a)$ range $m$-times through the
$\frac{j}{d},\, j\in\{0,...,d-1\}$. Thus the left-hand side of
(\ref{div1}) gives $m$-times
$$
c(\Zb)^k\,\sum_{j=0}^{d-1}\,\sum_{y\in  \Zb +\phi(x)+\frac{j}{d}} y^{-k}
=\,c(\Zb)^k\,d^{\,k} \,\sum_{y\in  \Zb +\phi(d\,x)} y^{-k}.
$$
This is clear when $y=0$ does not appear in the sums involved. When it does
one has, for $\epsilon \notin \frac{\Zb}{d}$,
$$
\sum_{j=0}^{d-1}\,\sum_{y\in  \Zb +\phi(x)+\frac{j}{d}} (y+\epsilon)^{-k}
=\,d^{\,k} \,\sum_{y\in  \Zb +\phi(d\,x)} (y+\,d\,\epsilon)^{-k}.
$$
Subtracting the pole part on both sides and equating the
finite values gives the desired equality, since for odd $k$
the value of $\epsilon_{k,a} ( \L , \phi)$ for $\phi(a)=0$
can be written as the finite value of
$$
 \sum_{y\in  \Lambda +\phi(a)}  (y+\epsilon)^{-k}.
$$

\n For even
$k$ this no longer holds and
the finite value $\g_k c(\Lambda)^k$ is replaced by
$(2^k-1)\,\g_k c(\Lambda)^k$.
 Thus when  $\phi(d\,x)\in \Zb$
one gets an additional term which is best taken care of
by multiplying the right hand side in Proposition \ref{div1}
by $(1-\pi_{\delta(d\,x)})$, with $\delta(y)$
the order of $y$ in $\mathbb{Q}\slash\mathbb{Z}$,
and adding
corresponding terms to the formula, which becomes
\begin{eqnarray}
\frac{1}{N}\,\sum_{N\,a = 0} e_{k,x+a}&= &
\, \sum_{d|N}\,\pi(N,d)(1-\pi_{\delta(d\,x)})\,d^{\,k-1}\,e_{k,d\, x}
\\ \,&+ &\, \g_k
\sum_{d|N}\,(d^{\,k-1}+\,d^{-1}(2^k-2))\,\pi(N,d)\,\pi_{\delta(d\,x)}
 \nonumber
\end{eqnarray}
These relations are more elaborate than the division relations for
trigonometric functions. They restrict to the latter on the subset
of invertible $\Qb$-lattices, for which all $\pi_n$, $n\neq 1$ are
zero and the only non-zero term in the r.h.s. is the term in
$d=N$. The standard discussion of Eisenstein series in higher
dimension is restricted to invertible $\Qb$-lattices, but in our
case the construction of the endomorphisms implemented by the
$\mu_n$ requires the above extension to non-invertible
$\Qb$-lattices. We shall now proceed to do it in dimension $2$.

\section{The Determinant part of the ${\rm GL}_2$-System}

 As we recalled in the previous sections,
the algebra of the 1-dimensional system can be described
as the  semigroup cross product
$$ C(\RR)\rtimes \N^\times. $$
Thus, one may wish to follow a similar
approach for the 2-dimensional case, by replacing
$C(\RR)$ by $C(M_2(\RR))$ and the semigroup
action of $\N^\times$ by the semigroup action
of $M_2(\Z)^+$. Such construction can be
carried out, as we discuss in this section, and it
corresponds to the ``determinant part''
of the $\GL_2$ system. It is useful to
analyze what
happens in this case first, before we
discuss the full $\GL_2$-system in
the next section. In fact, this will
show quite clearly where some important
technical issues arise.

\n For instance, just as in the case of the BC algebra,
where the time evolution acts on the isometries
$\mu_n$ by $n^{it}$ and leaves the elements of
$C(\RR)$ fixed, the time evolution here is
given by ${\rm Det}\,(m)^{it}$ on the isometries
implementing the semigroup action of
$m\in M_2(\Z)^+$, while leaving $C(M_2(\RR))$ pointwise fixed.
In this case, however, the vacuum state of
the corresponding Hamiltonian is highly degenerate,
because of the presence of the $\SL_2(\Z)$ symmetry.
This implies that the partition function and the
KMS states below critical temperature can only
be defined via the type II$_1$ trace $\Trace_\Gamma$.

\n This issue will be taken care more naturally
in the full $\GL_2$-system, by first taking the classical
quotient by $\Gamma=\SL_2(\Z)$ on the space
$M_2(\RR)\times \H$. This will resolve the
degeneracy of the vacuum state and the counting
of modes of the Hamiltonian will be on the
coset classes $\Gamma\backslash M_2(\Z)^+$.

\n The whole discussion of this section extends to
${\rm GL} (n)$ for arbitrary $n$ and we shall
briefly indicate how this is
done, but we stick to $n=2$ for definiteness.

\n We start with the
 action of the semigroup
\begin{equation}
\label{M2Z+}
 M_2 ({\mathbb Z})^+ = \{ m \in M_2 ({\mathbb Z}) , {\rm Det}\, (m) > 0 \}
={\rm GL}_2^+ ({\Q})\cap M_2(\RR)
\end{equation}
on the  compact space $ M_2 (R) $,
given by left multiplication
\begin{equation}
\label{actmr}
\rho \mapsto \; m\,\rho,
\end{equation}
where the product $m\rho$ takes place in $M_2 (R)$ using the natural
homomorphism
\begin{equation}
\label{inc}
M_2 ({\mathbb Z})^+ \to M_2 (R),
\end{equation}
which is the extension to two by two matrices of the inclusion of the ring
${\mathbb Z}$ in $\hat{\mathbb Z} = R$.
 The relevant $C^\ast$-algebra is the semi-group cross product
\begin{equation}
\label{cstar}
A = C(M_2(\RR))\rtimes M_2 ({\mathbb Z})^+ .
\end{equation}
It can be viewed as the $C^\ast$-algebra
$C^\ast(G_2)$ of the \'etale groupoid $G_2$
of pairs $(r, \rho)$, with $r\in {\rm GL}_2^+ ({\Q})$,
$\rho \in M_2(\RR)$ and
$r\, \rho \in M_2(\RR)$, where the product takes place
in $M_2(\A_f)$. The composition in $G_2$ is given by
$$
(r_1, \rho_1)\circ (r_2, \rho_2)=\,(r_1\,r_2, \rho_2)\,,
\quad {\rm if}\;r_2\, \rho_2=\,\rho_1\,
$$
and the convolution of functions by
\begin{equation}
\label{conv1}
f_1 \ast f_2(r,\rho):=\,\sum \,f_1(r s^{-1},s\,\rho)\,f_2(s,\rho)\, ,
\end{equation}
while the adjoint of $f$ is
\begin{equation}
\label{adjoint1}
f^\ast(r,\rho):=\,\overline{f(r^{-1},r\,\rho)}\,
\end{equation}
(\cf the analogous expressions \eqref{comp01}, \eqref{conv01},
\eqref{adjoint0} in the 1-dimensional case).

\smallskip

\n A homomorphism $G_2 \to H$ of the groupoid $G_2$ to an abelian
group $H$ determines a {\em dual action} of the Pontrjagin dual of
$H$ on the algebra of $G_2$, as in the case of the time evolution
$\sigma_t$, with $H=\R^*_+$ and its dual identified with $\R$. We
shall use the same term ``dual action'' for nonabelian $H$.

\smallskip

\n The main structure is given by the  {\em dual action} of
${\rm GL}_2^+ ({\R})$ corresponding to the groupoid
homomorphism $j$
\begin{equation}\label{jjj}
j\,:\;G_2\to {\rm GL}_2^+ ({\R})\,,\quad j(r,\rho)=r
\end{equation}
obtained from the inclusion ${\rm GL}_2^+ ({\Q})\subset{\rm GL}_2^+ ({\R})$.
As a derived piece of structure one gets the one parameter
group of automorphisms $\sigma_t \in {\rm Aut}(A)$ which is dual
to the determinant of the homomorphism $j$,
\begin{equation}
\label{sigt}
\sigma_t (f)(r, \rho):= {\rm Det}(r)^{it}\,f(r, \rho)\qqq f\in A\,.
\end{equation}
The obtained $C^\ast$-dynamical system $(A,\sigma_t)$ only
involves ${\rm Det}\, \circ j$ and it does not fully correspond to the BC
system. We shall make use of the full dual action of $\GL_2^+(\R)$
later in the construction of the full $\GL_2$ system.

\smallskip

\n The algebra $C(M_2(\RR))$ embeds as a $\ast$-subalgebra of
$A$. The analogs of the isometries $\mu_n,\,n\in \N^\times$
are the isometries $\mu_m,\,m\in M_2 ({\mathbb Z})^+$
given by
$$
\mu_m(m,\rho)=1\,,\quad \mu_m(r,\rho)=0 \qqq r\neq m\,.
$$
The range $\mu_m\,\mu_m^\ast$ of $\mu_m$ is the projection
given by the characteristic function of the
subset $P_m=\,m\,M_2(\RR)\subset M_2(\RR)$. It depends only
on the lattice $L=m(\Z^2)\subset \Z^2$.
Indeed, if $m,m' \in M_2 ({\mathbb Z})^+$ fulfill $m({\mathbb Z}^2) =
m'({\mathbb Z}^2)$, then  $m' =
m\gamma$ for some $\gamma \in \Gamma$, hence
$m \, M_2 (R) = m' \, M_2 (R)$. Thus, we shall label
this analog of the $\pi_n$ by lattices
\begin{equation}
\label{latticelabel0}
L \subset {\mathbb Z}^2 \mapsto \pi_L \in C(M_2(\RR)\, ,
\end{equation}
where $\pi_L$ is the characteristic function of $P_m$,
for any $m$ such that  $ m({\mathbb Z}^2)=L$.
The algebra generated by the
$\pi_L$ is then governed by
\begin{equation}
\label{eq11}
\pi_L \, \pi_{L'} = \pi_{L \cap L'} \, , \qquad \pi_{{\mathbb Z}^2} = 1 \,.
\end{equation}
In fact, the complete rules are  better expressed in terms of partial
isometries
$\mu(g,L)$, with $g \in {\rm GL}_2^+ ({\mathbb Q})$, $L \subset
{\mathbb Z}^2$ a lattice, and $g(L) \subset {\mathbb Z}^2$, satisfying
$$
\mu(g,L)(g,\rho)=\pi_L(\rho)\,,\quad
\mu(g,L)(r,\rho)=0 \qqq r\neq g\,.
$$
One has
\begin{equation}
\label{preslat}
 \mu (g_1 , L_1) \, \mu (g_2 ,
L_2) = \mu (g_1 g_2 , g_2^{-1} (L_1) \cap L_2) \, ,
\end{equation}
and
\begin{equation}
\label{preslat1}
\mu (g,L)^* = \mu (g^{-1} , g(L)) \,.
\end{equation}
The $\mu(g,L)$ generate the semi-group $C^\ast$-subalgebra
$C^* (M_2 ({\mathbb Z})^+)\subset A$ and together with
$C(M_2(\RR))$ they generate $A$.
The additional relations are
\begin{equation} \label{crossrel}
f\,  \mu(g,L)=\, \mu(g,L)\,f^g \qqq f\in C(M_2(\RR))\,,\quad
g \in {\rm GL}_2^+ ({\mathbb Q}),
\end{equation}
where $f^g(y):=f(gy)$ whenever $gy$ makes sense.

\smallskip

\n The  action of ${\rm GL}_2 ({\RR})$ on $M_2(\RR)$
by right multiplication commutes with the semi-group action
\eqref{actmr}
of $M_2 ({\mathbb Z})^+$ and with  the time evolution $\sigma_t$.
They define symmetries
$$
\a_\theta \in {\rm Aut}(A,\sigma)\,.
$$

\n Thus, we have a $C^\ast$-dynamical system with a compact group of
symmetries. The following results show how to construct
KMS$_\beta$-states for $\beta >2$.
We first describe a specific positive energy representation of the
$C^*$-dynamical sub-system ($C^* (M_2 ({\mathbb Z})^+)$, $\sigma_t$).
We let ${\mathcal H} = \ell^2 (M_2 ({\mathbb Z})^+)$ with cano\-nical basis
$\varepsilon_m$, $m \in M_2 ({\mathbb Z})^+$. We define $\pi (\mu (g,L))$
as the partial isometry in ${\mathcal H}$ with initial domain given by the
span of
\begin{equation}
\label{eps-m}
\varepsilon_m \, , \quad m = \begin{bmatrix} m_{11} &m_{12} \\ m_{21}
&m_{22} \end{bmatrix} \, , \quad (m_{11} , m_{21}) \in L \, , \ (m_{12} ,
m_{22}) \in L \, ,
\end{equation}
\ie matrices  $m$ whose columns
 belong to the lattice $L \subset {\mathbb
Z}^2$. On this domain we define the action of $\pi (\mu (g,L))$ by
\begin{equation}
\label{defrep}
\pi (\mu (g,L)) \, \varepsilon_m = \varepsilon_{gm} \, .
\end{equation}
Notice that the columns of $gm$ belong to $gL$.

\medskip

\begin{prop} \label{reppi} 1)  $\pi$ is an involutive representation
of $C^* (M_2 ({\mathbb Z})^+)$ in ${\mathcal H}$.

\smallskip

\noindent 2) The Hamiltonian $H$ given by $H \varepsilon_m = \log {\rm
Det} (m) \, \varepsilon_m$ is positive and implements the time
evolution $\sigma_t$:
$$
\pi (\sigma_t (x)) = e^{itH} \, \pi (x) \, e^{-itH} \qquad \forall \, x \in
C^* (M_2 ({\mathbb Z})^+) \, .
$$

\noindent 3)
$\Gamma = {\rm SL}_2 ({\mathbb Z})$ acts
on the right in ${\mathcal H}$ by
\begin{equation}
\nonumber
\rho(\g)\,\epsilon_m := \epsilon_{m\,\g^{-1}}\qqq \g \in \Gamma\,,\quad m\in
M_2 ({\mathbb Z})^+\,.
\end{equation}
and this action commutes with
$\pi(C^* (M_2 ({\mathbb Z})^+)$.
\end{prop}

\noindent {\it Proof.} The map $m \mapsto gm$ is injective so that $\pi
(\mu (g,L))$ is a partial isometry. Its range is the set of $h \in M_2
({\mathbb Z})^+$ of the form $g \, m$ where ${\rm Det} (m) > 0$ and the
columns of $m$ are in $L$. This means that ${\rm Det} (h) > 0$ and the
columns of $h$ are in $gL \subset {\mathbb Z}^2$. This shows that
\begin{equation}
\label{pimugL}
\pi (\mu (g,L))^* = \pi (\mu (g^{-1} , gL)),
\end{equation}
so that $\pi$ is involutive on these elements.

\n Then the support of $\pi (\mu (g_1 , L_1)) \, \pi (\mu (g_2 , L_2))$ is
formed by the $\varepsilon_m$ with columns of $m$ in $L_2$, such that the
columns of $g_2 m$ are in $L_1$. This is the same as the support of $\pi
(\mu (g_1 g_2 , g_2^{-1} L_1 \cap L_2))$ and the two partial isometries
agree there. Thus, we get
\begin{equation}
\label{replat}
\pi (\mu (g_1 , L_1) \, \mu (g_2 , L_2)) = \pi (\mu (g_1 , L_1)) \, \pi
(\mu (g_2 , L_2)) \, .
\end{equation}
Next, using (\ref{defrep}) we see that
\begin{equation}
\label{hamil}
H \, \pi (\mu (g,L)) - \pi (\mu (g,L)) \, H = \log ({\rm Det} \, g) \, \pi
(\mu (g,L)),
\end{equation}
since both sides vanish on the kernel while on the support one can use the
multiplicativity of Det.

\n Now $\Gamma = {\rm SL}_2 ({\mathbb Z})$ acts
on the right in ${\mathcal H}$ by
\begin{equation}
\label{sl2}
\rho(\g)\,\epsilon_m := \epsilon_{m\,\g^{-1}}\qqq \g \in \Gamma\,,\quad m\in
M_2 ({\mathbb Z})^+\,
\end{equation}
and this action commutes by construction with the algebra
$\pi(C^* (M_2 ({\mathbb Z})^+)$. $\Box$

\smallskip

\n The image $\rho(C^* (\Gamma))$ generates a type II$_1$
factor in ${\mathcal H}$, hence one can evaluate
the corresponding trace ${\rm Trace}_{\Gamma}$
on any element of its commutant. We let
\begin{equation}
\label{betakms}
\varphi_\beta(x):= {\rm Trace}_{\Gamma} (\pi(x)\;e^{-\beta H})\qqq x\in
C^* (M_2 ({\mathbb Z})^+)
\end{equation}
and we define the normalization factor by
\begin{equation}
\label{ZetaTrGamma}
Z(\beta) = {\rm Trace}_{\Gamma} (e^{-\beta H}) \, .
\end{equation}
We then have the following:

\medskip

\begin{lem} 1) The normalization factor $Z(\beta)$ is given by
$$
Z(\beta) = {\z} (\beta) \, {\z} (\beta-1),
$$
where $\z$ is the Riemann ${\z}$-function.

\smallskip

\n 2)  For all $\beta >2$, $Z^{-1}\;\varphi_\beta$ is a KMS$_\beta$
state on  $C^* (M_2 ({\mathbb Z})^+)$.
\end{lem}

\medskip

\noindent {\it Proof.} Any sublattice $L \subset {\mathbb Z}^2$ is uniquely
of the form $L = \begin{bmatrix} a &b \\ 0 &d \end{bmatrix} \, {\mathbb
Z}^2$, where $a,d \geq 1$, $0 \leq b < d$ (\cf \cite{Serre} p.~161). Thus, the
type II$_1$ dimension of the action of $\Gamma$ in the subspace of
${\mathcal H}$ spanned by the $\varepsilon_m$ with ${\rm Det} \, m  = N$
is the same as the cardinality of the quotient of $\{ m \in M_2 ({\mathbb
Z})^+ , {\rm Det} \, m  = N \}$ by $\Gamma$ acting on the right. This
is equal to the cardinality of the set of matrices $\begin{bmatrix} a
&b \\ 0 &d \end{bmatrix}$ as above with determinant $=N$. This gives
$\sigma_1 (N) = \underset{d\mid N}{\sum} \ d$. Thus, $Z(\beta)$ is
given by
\begin{equation}
\label{Zeta12}
\sum_{N=1}^{\infty} \ \frac{\sigma_1 (N)}{N^{\beta}} = {\z} (\beta)
\, {\z} (\beta-1) \, .
\end{equation}
One checks the KMS$_\beta$-property of $\varphi_\beta$
using the trace property of ${\rm Trace}_{\Gamma}$
together with the second equality in Proposition \ref{reppi}. $\Box$

\begin{prop} \label{gamman}
1) For any $\theta\in {\rm GL}_2 ({\RR})$
 the formula
$$
\pi_{\theta}(f)\,\epsilon_m:=\,f(m\,\theta)\,\epsilon_m
\qqq m \in M_2 ({\mathbb Z})^+\,
$$
extends
the representation $\pi$ to an involutive representation
$\pi_{\theta}$ of the cross product
$A=\,C(M_2(\RR))\rtimes M_2 ({\mathbb Z})^+$
 in ${\mathcal H}$.

\smallskip

\n 2) Let $f\in C(M_2(\Z/N\Z))\subset C(M_2(\RR))$. Then
$\displaystyle
\pi_{\theta}(f)\,\in \rho(\G_N)'\,
$
where $\G_N$ is the congruence subgroup of level $N$.

\smallskip

\n 3) For each  $\beta >2$
the formula
$$
\psi_\beta(x):= {\rm Lim_{N\to \infty}}\quad
Z_N^{-1}\,  {\rm Trace}_{\Gamma_N} (\pi_\theta(x)\;e^{-\beta H})\qqq x\in A
$$
defines a KMS$_\beta$
state on $A$,
where $Z_N:={\rm Trace}_{\Gamma_N} (e^{-\beta H})$.
\end{prop}

\noindent  {\it Proof.} 1) The invertibility of
$\theta$ shows that
letting $f_L$ be the characteristic function of $P_L $ one has
$\pi_{\theta}(f_L)= \pi_L$ independently of $\theta$.
Indeed $f_L(m\,\theta )=1$ iff $m\,\theta \in P_L$
and this holds iff $m(\Z^2)\subset L$.

\n To check (\ref{crossrel}) one uses
$$
f(g m\,\theta)=f^g(m\,\theta))\qqq g \in {\rm GL}_2^+ ({\mathbb Q})\,.
$$

\smallskip
\n 2) Let $p_N:M_2(\RR)\to M_2(\Z/N\Z)$ be the canonical
projection. It is a ring homomorphism. Let then
$f = h \circ p_N$ where $h$ is a function on  $M_2(\Z/N\Z)$.
One has, for any $\g \in \G_N$,
\begin{eqnarray}
\pi_{\theta}(f)\,\rho(\g)\,\epsilon_m  &= &\pi_{\theta}(f)\,
\epsilon_{m\,\g^{-1}} \,=\nonumber \\
f(m\,\g^{-1}\theta) \,
\epsilon_{m\,\g^{-1}}&= &h(p_N(m\,\g^{-1}\theta))\,
\epsilon_{m\,\g^{-1}}\, .\nonumber
\end{eqnarray}
The equality $ p_N(\g)=1$ shows that
$$
p_N(m\,\g^{-1}\theta)=p_N(m)\,p_N(\g^{-1})\,p_N(\theta)
=p_N(m)\,p_N(\theta)=p_N(m\,\theta),
$$
hence
$$
\pi_{\theta}(f)\,\rho(\g)\,\epsilon_m
=\,\rho(\g)\,\pi_{\theta}(f)\,\epsilon_m \,.
$$

\smallskip
\n 3) One uses 2)
to show that, for all $N$ and  $f \in
 C(M_2(\Z/N\Z))\subset C(M_2(\RR))$, the products
$f\,  \mu(g,L)$ belong to the commutant of
$\rho(\G_N)$. Since $\G_N$ has finite index in $\G$
it follows that for $\beta >2$ one has
$Z_N:={\rm Trace}_{\Gamma_N} (e^{-\beta H})<\infty$.
Thus, the limit defining $\psi_\beta(x)$
makes sense on a norm dense subalgebra of $A$
and extends to a state on $A$ by uniform continuity.
One checks the KMS$_\beta$ condition on the dense
subalgebra in the same way as above. $\Box$\\

\n It is not difficult to extend the above discussion
to arbitrary $n$ using (\cite{Shimura}).
The normalization factor is then given by
$$
Z(\beta) = \prod_0^{n-1} \, {\z} (\beta-k)\,.
$$

\n What happens, however, is that the states $\psi_\beta$
only depend on the determinant of $\theta$. This
shows that the above construction should be extended
to involve not only the one-parameter group $\sigma_t$
but in fact the whole dual action given by the
groupoid homomorphism (\ref{jjj}).

\smallskip

\begin{defn}\label{Defn-groupoid}
Given a groupoid $G$ and a homomorphism $j:G \to H$ to a group
$H$, the ``cross product" groupoid $G\times_j H$ is defined as the
product $G\times H$ with units $G^{(0)}\times H$, range and source
maps
$$
r(\g,\a):=(r(\g),j(\g)\,\a)\,,\quad s(\g,\a):=(s(\g),\a)
$$
and composition
$$
(\g_1,\a_1)\circ (\g_2,\a_2):= (\g_1\circ \g_2,\a_2)\,.
$$
\end{defn}

\smallskip

\n In our case this cross product
$\tilde{G_2}=G_2\times_j {\rm GL}_2^+ ({\R})$ corresponds
to the groupoid of the partially defined
action of ${\rm GL}_2^+ ({\Q})$
on the locally compact space $Z_0$ of pairs
$(\rho, \a) \in M_2(\RR)\times {\rm GL}_2^+ ({\R})$
given  by
$$
g (\rho, \a):= (g\, \rho, g\,\a) \qqq g \in {\rm GL}_2^+ ({\Q})\,,
\quad g\, \rho\in M_2(\RR)\,.
$$
Since the subgroup $\G \subset {\rm GL}_2^+ ({\Q})$
acts freely and properly by translation on ${\rm GL}_2^+ ({\R})$,
one obtains a Morita equivalent groupoid $S_2$ by dividing
$\tilde{G_2}$ by the following action of $\G \times \G$:
\begin{equation}\label{S2def1}
(\g_1,\g_2).(g, \rho,\a):=\,(\g_1 \,g\g_2^{-1},\,\g_2 \rho,\g_2\a).
\end{equation}
The space of units $S_2^{(0)}$
is the quotient
$\G\backslash Z_0$. We let $p: Z_0\to \G\backslash Z_0$
be the quotient map. The range and source maps are given by
$$
r(g, \rho,\a):=p(g\rho,g\a)\,,\quad s(g, \rho,\a):=p(\rho,\a)
$$
and the composition  is given by
$$
(g_1, \rho_1,\a_1)\circ (g_2, \rho_2,\a_2)=(g_1g_2, \rho_2,\a_2),
$$
which passes to $\G \times \G$-orbits.

\section{Commensurability of $\Qb$-Lattices in $\Cb$
and the full ${\rm GL}_2 $-System}\label{CQlattices}

\n We shall now describe the full ${\rm GL}_2$ $C^\ast$-dynamical
system $(A,\sigma_t)$. It is obtained from the system of the previous
section by taking a cross product with the dual action of
${\rm GL}_2^+ ({\R})$ \ie from the groupoid
$S_2$ that we just described.  It admits an equivalent
and more geometric description in terms of
the notion of
{\em commensurability} between $\Q$-lattices
developed in section 1.6 above and we shall
follow both points of view. The $C^\ast$-algebra $A$
is a Hecke algebra, which is a variant of
the {\em modular Hecke algebra} defined in
(\cite{CM10}).
 Recall from section 1.6:

\begin{defn} 1) A $\Qb$-lattice in $\Cb$ is a
 pair $ \, ( \L , \phi) \,$, with $
\L $ a lattice in $\mathbb{C}$, and
 $\displaystyle \phi :  \Q^2\slash\Z^2
\longrightarrow\mathbb{Q}\,\L \slash \L $
an homomorphism of abelian groups.

\n 2) Two $\Qb$-lattices  $ \, ( \L_j , \phi_j) \,$
are commensurable iff $ \L_j$ are commensurable and
$\phi_1-\phi_2=0$ modulo $\L=\L_1+\L_2$.
\end{defn}

\n This is an equivalence relation $\mathcal R$ between $\Qb$-lattices
(Proposition \ref{comm}).
We use the basis $\{e_1=1,e_2=-i\}$
of the $\R$-vector space $\C$ to let ${\rm GL}_2^+ ({\R})$
act on $\C$ as $\R$-linear transformations. We let
$$
\L_0:=\,\Z\,e_1 + \Z\,e_2= \Z + i \Z
$$
Also we view $\rho \in M_2(\RR)$ as the homomorphism
$$
\rho :\Q^2\slash\Z^2
\longrightarrow\Q.\L_0
\slash \L_0\,,\quad \rho(a)= \rho_1(a)e_1 +\,\rho_2(a)e_2\,.
$$

\begin{prop} \label{groupoid2}
The  map
$$
\g(g,\rho,\a)=((\a^{-1}\, g^{-1}\L_0,\,\a^{-1}\,\rho)
\,,(\a^{-1}\,\L_0,\,\a^{-1}\,\rho))\qqq (g,\rho,\a) \in S_2\,
$$
defines an isomorphism of locally compact
\'etale groupoids between $S_2$ and the
 equivalence relation $\Rc$ of commensurability
on the space of $\Qb$-lattices in $\C$.
\end{prop}

\noindent  {\it Proof.} The proof is the same
as for Proposition \ref{groupoid1}. $\Box$

\smallskip

\begin{figure}
\begin{center}
\epsfig{file=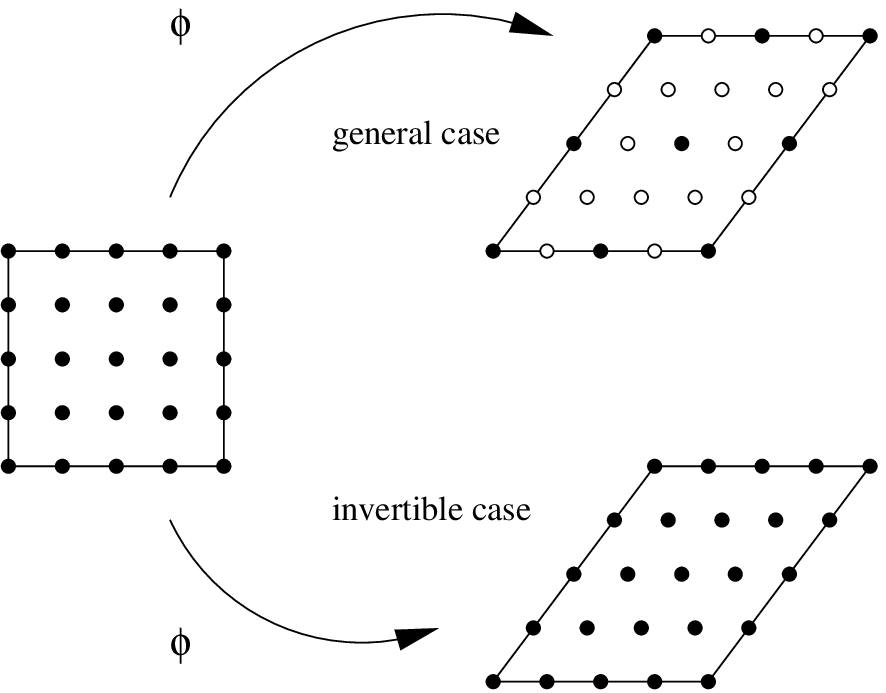}
\end{center}
\caption{$\Q$-Lattices in $\Cb$.
\label{Qlattices}}
\end{figure}

\n We shall now describe the quotient of $S_2\sim \Rc$
by the natural scaling action of $\Cb^\ast$.
We view $\Cb^\ast$ as a subgroup of ${\rm GL}_2^+ ({\R})$ by the map
\begin{equation}\label{inclGL2R}
\lambda= a+i\,b \in \Cb^\ast \mapsto \begin{bmatrix} a &b \\ -b &a
\end{bmatrix}
\in {\rm GL}_2^+ ({\R})
\end{equation}
and identify the quotient ${\rm GL}_2^+ ({\R})/\Cb^\ast$
with $\H$ by the map
\begin{equation}\label{upper}
\a\in {\rm GL}_2^+ ({\R})\mapsto \tau=\a(i)\in \H.
\end{equation}

\n Given a pair $ \, ( \L_j , \phi_j) \,$
of commensurable $\Qb$-lattices and a
non zero complex number $\lambda \in \Cb^\ast $
the pair $ \, ( \lambda\L_j , \lambda\phi_j) \,$
is still a pair
of commensurable $\Qb$-lattices. Moreover, one has
\begin{equation}\label{homogeneous}
\g(g,\rho,\a\,\lambda^{-1})=\,\lambda\,\g(g,\rho,\a) \qqq \lambda
\in\C^\ast\,.
\end{equation}
The scaling action of $\C^\ast$ on ${\mathbb Q}$-lattices in
${\mathbb C}$ is not free, since the lattice $\L_0$ for instance
is invariant under multiplication by $i$. It follows that
the quotient $S_2/\C^\ast\sim \Rc/\C^\ast$ is not a groupoid.
One can nevertheless define its convolution algebra
in a straightforward manner by restricting the convolution
product on  $S_2\sim \Rc$ to functions which are homogeneous
of {\em weight} $0$, where weight $k$ means
\begin{equation}\label{weightk}
f(g,\rho,\a\,\lambda)=\,\lambda^k\, f(g,\rho,\a)\qqq \lambda \in \Cb^\ast\,.
\end{equation}
Let
\begin{equation}\label{Yspace}
Y=M_2 (R) \times \H,
\end{equation}
endowed with the natural action of ${\rm GL}_2^+ ({\Q})$ by
\begin{equation}
\label{actm2}
\g \cdot (\rho ,
\tau) = \left( \g\,\rho , \;\frac{a\tau + b}{c\tau + d} \right),
\end{equation}
for $\g = \begin{bmatrix} a &b \\ c &d \end{bmatrix} \in
{\rm GL}_2^+ ({\Q})$ and $(\rho ,
\tau) \in Y$.
Let then
\begin{equation}\label{Zspace}
Z\subset \G \backslash
{\rm GL}_2^+ ({\Q})\times_\G Y
\end{equation}
 be the locally compact
space quotient of
$
\{(g,y)\in {\rm GL}_2^+ ({\Q})\times Y\, ,\; g\,y \in Y\}
$
by the following action of $\G\times \G$ :
$$
(g,y)\mapsto (\g_1 \,g\,\g_2^{-1},
\g_2\,y) \qqq \g_j\in \G\,.
$$
The natural lift of the quotient map (\ref{upper}),
together with proposition \ref{groupoid2}, first gives
the identification of the quotient of  $Y$ by $\Gamma$
with the space of ${\mathbb Q}$-lattices in ${\mathbb C}
$ up to scaling, realized by the map
\begin{equation}
\label{maptheta}
\t:\Gamma \backslash Y \to
 (\hbox{Space of ${\mathbb Q}$-lattices in
${\mathbb C}$}) / {\mathbb C}^* \,=\,X\, ,
\end{equation}
$$
\t(\rho , \tau)= (\Lambda, \phi)\,,\quad
\Lambda = {\mathbb Z} + \, {\mathbb Z} \tau \,,\quad
\phi (x) =
\rho_1 (x) - \tau \rho_2 (x)\,.
$$
It also gives
 the
isomorphism $\theta: S_2/\C^\ast
= Z \to \Rc/\C^\ast$,
\begin{equation}\label{maptheta1}
\theta(g,y)=(\lambda \t(gy),\t(y)),
\end{equation}
where $\lambda={\rm Det}(g)^{-1}(c\tau + d)$
for $\g = \begin{bmatrix} a &b \\ c &d \end{bmatrix} \in
{\rm GL}_2^+ ({\Q})$ and $y=(\rho ,
\tau) \in Y$.

\n We let $\Ac=C_c(Z)$ be the space of continuous
functions with compact support on $Z$. We view elements
$f\in \Ac$ as functions on ${\rm GL}_2^+ ({\Q})\times Y$
such that
$$
f(\g \,g,y)=f(g,y)\quad f(g \,\g,y)=f(g,\g\,y)\qqq \g\in\G\,,\;g\in
{\rm GL}_2^+ ({\Q})\,,\;y\in Y\,.
$$
 This does not imply that
$f(g,y)$ only depends on the orbit $\G.y$ but that
it only depends on the orbit of $y$ under the congruence
subgroup $\G \, \cap g^{-1} \G \,
g \,$.

\n We define the
convolution product of two such functions by
\begin{equation} \label{convpro}
(f_1 * f_2)(g,y) \,:= \, \sum_{ h \in \G \backslash {\rm GL}_2^+
({\Q}),\,hy\in Y}
\, f_1(g h^{-1},h y) \; f_2(h,y)
\end{equation}
and the adjoint by
\begin{equation} \label{adjoint}
f^\ast(g,y)\,:= \,\overline{f(g^{-1},g\,y)}.
\end{equation}
Notice that these rules combine (\ref{convhecke}) and (\ref{conv1}).

\smallskip

\n For any $x\in X$ we let $c(x)$ be the commensurability
class of $x$. It is a countable subset of $X$ and we
want to define a natural representation in $l^2(c(x))$.
We let $p$ be the quotient map
from $Y$ to $X$. Let $y\in Y$ with $p(y)=x$ be an element in
the preimage of $x$. Let
$$G_y=\{g\in {\rm GL}_2^+ ({\Q})\,|\,gy\in Y\}\,.$$
The natural map $g\in G_y\mapsto p(gy)\in X$ is a surjection from
$\G \backslash G_y$ to $c(x)$ but it fails to be injective in
degenerate cases such as $y=(0,\tau)$ with $\tau\in \H$ a complex
multiplication point (\cf Lemma \ref{fixedpt}). This corresponds
to the phenomenon of holonomy in the context of foliations
(\cite{Cosurvey}). To handle it one defines the representation
$\pi_y$ directly in the Hilbert space $\Hc_y=l^2(\G \backslash
G_y)$ of left $\G$-invariant functions on $G_y$ by
\begin{equation} \label{repl}
(\pi_y(f)\, \xi)(g ):=
\, \sum_{ h \in \G \backslash G_{y}}
\, f(g h^{-1},hy) \, \xi(h)\qqq g\in G_y\,,
\end{equation}
for $f\in \Ac$ and $\xi\in \Hc_y$.

\begin{prop}\label{groupoidG}
1) The vector space $\Ac$ endowed with the product $\ast$
and the adjoint $f\mapsto f^\ast$ is an involutive algebra.

\smallskip
\n 2) For any $y\in Y$, $\pi_y$ defines a unitary
representation of  $\Ac$ in $\Hc_y$ whose unitary
equivalence class only depends on $x=p(y)$.

\smallskip
\n 3) The completion of $\Ac$ for the norm given by
$$
||f||:= {\rm Sup}_{y\in Y}\;||\pi_y(f)||
$$
is a $C^\ast$-algebra.
\end{prop}

\noindent  The proof of (1) and (2) is similar to (\cite{CM10},
Proposition 2). Using the compactness of the support of $f$, one
shows that the supremum is finite for any $f\in \Ac$ (cf.
\cite{Cosurvey}). $\Box$\\

\begin{rem}\label{M2Rremark}{\em
The locally compact space $Z$ of \eqref{Zspace} is not a groupoid,
due to the torsion elements in $\Gamma$, which give nontrivial
isotropy under scaling, for the square and equilateral lattices.
Nonetheless, Proposition \ref{groupoidG} yields a well defined
$C^*$-algebra. This can be viewed as a subalgebra of the
$C^*$-algebra of the groupoid obtained by replacing $\Gamma$ by
its commutator subgroup in the definition of $S_2$ as in
\eqref{S2def1}.}
\end{rem}

\n  We let
$\sigma_t$ be the one parameter
group of automorphisms of $A$ given by
\begin{equation}
\label{modaut2}
\sigma_t (f)(g,y) = ({\rm Det} \, g)^{it} \, f(g,y) \, .
\end{equation}
Notice that since $X$ is not compact (but still locally compact)
the  $C^\ast$-algebra $A$ does not have a unit, hence the
discussion of Proposition \ref{KMSnonunital} applies.

\smallskip

\n The one parameter group $\sigma_{2t}$ (\ref{modaut2}) is the
modular automorphism group associated to the regular
representation of $\Ac$. To obtain the latter we endow $X=\G
\backslash Y$ with the measure
$$
dy=\,d\rho \times d\mu(\tau),
$$
where $d\rho=\,\prod d\rho_{ij}$ is the normalized
Haar measure of the additive compact group $M_2(\RR)$
and $d\mu(\tau)$ is the Riemannian volume form in
$\H$ for the Poincar\'e metric, normalized so that $\mu(\G \backslash \H)=1$.
We then get the following result.

\begin{prop} \label{regrep}
The expression
\begin{equation}\label{KMS2state}
\varphi(f)=\int_X\, f(1,y)\,dy\,.
\end{equation}
defines a state on $A$, which
is a KMS$_2$ state for the one parameter group $\sigma_t$.
\end{prop}

\noindent {\it Proof.} At the measure theory level, the
quotient $X=\G \backslash Y$ is the total space over
$\G \backslash \H$ of a bundle with fiber the probability
space $M_2(\RR)/\{\pm 1\}$, thus the total mass
$\int_X\,dy\,=1$. One gets
$$
\varphi(f^\ast \ast f)=\int_X\,
\sum_{ h \in \G \backslash {\rm GL}_2^+ ({\Q}),\,hy\in Y}
 \overline{f(h,y)} f(h,y) \,dy \qqq f\in \Ac \,,
$$
which suffices to get the Hilbert space $\Hc$ of the
regular representation and the cyclic vector $\xi$
implementing the state $\varphi$, which corresponds to
$$
\xi(g,y)=0 \qqq g\notin \G\,,\quad \xi(1,y)=1 \qqq y\in X \,.
$$
The measure $d\rho$ is the product of the additive Haar
measures on column vectors, hence one gets
$$
d(g\,\rho)=\, ({\rm Det}\,g)^{-2}\, d\rho\qqq g\in {\rm GL}_2^+ ({\Q})\,.
$$
Let us prove that $\varphi$ is a KMS$_2$ state. The above equality shows
that, for any compactly supported continuous function $ \a$ on
$ \G \backslash {\rm GL}_2^+ ({\Q}) \times_\G Y$, one has
\begin{equation}\label{parts}
\int_X\,
\sum_{ h \in \G \backslash {\rm GL}_2^+ ({\Q})}
 \a(h,y) \,dy\,=
\int_X\,
\sum_{ k \in \G \backslash {\rm GL}_2^+ ({\Q})}
 \a(k^{-1},k\,y) \,({\rm Det}\,k)^{-2} \,dy\, .
\end{equation}
Let then $f_j\in \Ac$ and define $\a(h,y)=0$ unless $hy\in Y$ while otherwise
$$
\a(h,y) = \,  f_1(h^{-1},h\,y)\, f_2(h,y)
({\rm Det}\,h)^{it-2} .
$$
The l.h.s. of (\ref{parts}) is then equal to
$\varphi(f_1\,\sigma_z(f_2))$ for $z=t +2i$.
The r.h.s. of (\ref{parts}) gives
$\varphi(\sigma_t(f_2)\,f_1)$ and (\ref{parts})
gives the desired equality
$\varphi(f_1\,\sigma_{t +2i}(f_2))=\varphi(\sigma_t(f_2)\,f_1)$. $\Box$\\

\smallskip

\n We can now state the main result on the analysis of KMS states
on the $C^\ast$-dynamical system $(A,\sigma_t)$. Recall that a
$\Q$-lattice $l=( \L , \phi)$ is invertible if $\phi$ is an
isomorphism (Definition \ref{Defn-invQlat}). We have the following
result.

\smallskip

\begin{thm} \label{thmKMSst}
1) For each invertible $\Q$-lattice $l=( \L , \phi)$,
the representation $\pi_l$ is a positive energy
representation of the $C^\ast$-dynamical system $(A,\sigma_t)$.

\smallskip
\n 2) For $\beta >2$ and $l=( \L , \phi)$ an invertible
$\Q$-lattice, the formula
$$
\varphi_{\beta,l}(f)= Z^{-1}\, \sum_{\Gamma \backslash M_2(\Z)^+ }
\,f(1, m \,\rho, m(\tau))\,{\rm Det}(m)^{-\beta}\,,
$$
defines an extremal KMS$_\beta$ state $\varphi_{\beta,l}$ on
$(A,\sigma_t)$, where $Z=\z(\beta)\,\z(\beta-1)$ is the partition
function.

\smallskip
\n 3) For $\beta >2$ the map $l\mapsto \varphi_{\beta,l}$
is a bijection from the space of invertible  $\Q$-lattices
(up to scaling)
to the space $\Ec_\beta$ of extremal KMS$_\beta$
states  on $(A,\sigma_t)$.
\end{thm}

\n The proof of 1) reflects the following fact, which in essence
shows that the invertible $\Q$-lattices are {\em ground states}
for our system.
\begin{lem} 1) Let $s: \Q^2/\Z^2 \to \Q^2$
be a section of the projection
$\pi: \Q^2 \to \Q^2/\Z^2$. Then the set of $s(a+b)-
s(a)-s(b)\,, \; a, b \in \Q^2/\Z^2$, generates  $\Z^2$.

\smallskip
\n 2) Let $l=( \L , \phi)$ be an invertible $\Q$-lattice
and $l'=( \L' , \phi')$ be commensurable with $l$. Then
$\L \subset \L'$.
\end{lem}

\noindent  {\it Proof.} 1) Let $L \subset \Z^2$ be the subgroup
generated by the $s(a+b)- s(a)-s(b)\,, \; a, b \in \Q^2/\Z^2$. If
$L \neq \Z^2$ we can assume, after a change of basis, that for
some prime number $p$ one has $L \subset p \, \Z \oplus \Z$.
Restricting $s$ to the $p$-torsion elements of $\Q^2/\Z^2$ and
multiplying it by $p$, we get a morphism of groups
$$\Z/p\,\Z\oplus\Z/p\,\Z\to \Z/p^2\Z\oplus\Z/p\,\Z,$$
which is a section of the projection
$$\Z/p^2\Z\oplus\Z/p\,\Z\to \Z/p\,\Z\oplus\Z/p\,\Z.$$
This gives a contradiction, since the group
$\Z/p^2\Z\oplus\Z/p\,\Z$ contains elements of order $p^2$.

\n 2) Let $s$ (resp. $s'$) be a lift of $\phi$ modulo $\L$ (resp.
of $\phi'$ modulo $\L'$). Since $\phi-\phi'=0$ modulo
$\L''=\L+\L'$ one has $s(a)-s'(a) \in \L+\L'$ for all $a \in
\Q^2/\Z^2$. This allows one to correct $s$ modulo $\L$ and $s'$
modulo $\L'$ so that $s=s'$. Then for any $a, b \in \Q^2/\Z^2$ one
has $s(a+b)-s(a)-s(b) \in \L \cap \L'$ and the first part of the
lemma together with the invertibility of $\phi$ show that $\L \cap
\L'=\L$. $\Box$\\

\n Given $y\in Y$ we let $H_y$
be the diagonal operator in $\Hc_y$
given by
\begin{equation}
(H_y \, \xi)(h):= \log ({\rm Det}(h))\,\; \xi(h)\qqq h\in G_y
\end{equation}
It implements the one parameter group $\sigma_t$ \ie
\begin{equation}\label{implement}
\pi_y(\sigma_t(x))=\, e^{itH_y}\,\pi_y(x)\,e^{-itH_y}\qqq x\in A\,.
\end{equation}
In general the operator $H_y$ is not positive but when the lattice
$l=( \L , \phi)=\theta(p(y))$ is invertible one has $${\rm Det}(h)
\in \N^\ast \qqq h\in G_y,$$ hence $H_y\geq 0$. This proves the
first part of the theorem. The basis of the Hilbert space $\Hc_y$
is then labeled by the lattices $\L'$ containing $\L$ and the
operator $H_y$ is diagonal with eigenvalues the logarithms of the
orders $\L':\L$. Equivalently, one can label the orthonormal basis
$\epsilon_m$ by the coset space $\Gamma \backslash M_2(\Z)^+$.
Thus, the same counting as in the previous section
(cf.\cite{Shimura}) shows that $$Z=\,{\rm Trace} (e^{-\beta
H_y})=\z(\beta)\,\z(\beta-1)$$ and in particular that it is finite
for $\beta>2$. The KMS$_\beta$ property of the functional
$$
\varphi_{\beta,l}(f)=
Z^{-1}\,  {\rm Trace} (\pi_y(f)\;e^{-\beta H_y})
 $$
then follows from (\ref{implement}).
One has, using (\ref{repl}) for $y=(\rho,\tau)\in Y$,
$$
\langle \pi_y(f)(\epsilon_m),\;\epsilon_m \rangle\, = \,f(1, m
\,\rho, m(\tau))\,,
$$
hence we get the following formula for $\varphi_{\beta,l}$:
\begin{equation}\label{kmsexplicit}
\varphi_{\beta,l}(f)= Z^{-1}\,
\sum_{\Gamma \backslash M_2(\Z)^+
} \,f(1, m \,\rho, m(\tau))\,{\rm Det}(m)^{-\beta}\,.
\end{equation}
Finally, the irreducibility of the representation $\pi_y$ follows
as in \cite{Cosurvey} p.562 using the absence of holonomy for
invertible $\Q$-lattices. This completes the proof of 2) of Theorem
\ref{thmKMSst}.

\smallskip

\n In order to prove 3) of Theorem \ref{thmKMSst} we shall proceed in
two steps. The first
shows (Proposition \ref{measure} below) that KMS$_\beta$ states
are given by measures on the space $X$ of $\Q$-lattices (up to
scaling). The second shows that when $\beta > 2$ this measure is
carried by the commensurability classes of invertible
$\Q$-lattices.

\smallskip

\n We first describe the stabilizers of the action of $\GL^+_2(\Q)$
on the space of $\Q$-lattices in $\C$.

\smallskip

\begin{lem} \label{fixedpt} Let $g  \in
{\rm GL}_2^+ ({\Q}),\, g\neq 1$ and $y\in Y,\, y=(\rho,\tau)$ such that
$gy=\,y$. Then $\rho=\,0$. Moreover $g\in \Q^\ast\subset
{\rm GL}_2^+ ({\Q})$ unless $\tau$ is an imaginary quadratic number
in which case $g\in K^\ast\subset
{\rm GL}_2^+ ({\Q})$ where $K=\Q(\tau)$ is the corresponding quadratic field.
\end{lem}

\noindent {\it Proof.} Let $g =\begin{bmatrix} a &b \\ c &d \end{bmatrix} $,
then $g(\tau)=\,\tau$ means $a\tau + b=\,\tau\,(c\tau + d)$. If
$c\neq 0$ this shows that $\tau$ is an imaginary quadratic number.
Let $K=\Q(\tau)$ be the corresponding field and let
$\{\tau,1\}$ be the natural basis of $K$ over $\Q$.
Then the multiplication by $(c\tau + d)$ is given by
the transpose of the matrix $g$. Since $g\neq 1$ and $K$ is a field
we get
$$
g -1=\begin{bmatrix} a-1 &b \\ c &d-1 \end{bmatrix} \in {\rm GL}_2 ({\Q})
$$
and thus $(g-1)\rho =0$ implies $\rho =0$. If $c=0$ then
$a\tau + b=\,\tau\,(c\tau + d)$ implies
$a=d,\, b=0$ so that
$g\in \Q^\ast\subset
{\rm GL}_2^+ ({\Q})$.  Since $g\neq 1$ one gets $\rho =0$.  $\Box$\\

\n We let $X=\Gamma \backslash Y $ be the quotient and
$p: \,Y \to X$ the quotient map. We let $F$ be the closed
subset of $X$, $F=p(\{(0,\tau);\tau\in \H\})$.

\begin{lem} \label{corresp} Let $g  \in
{\rm GL}_2^+ ({\Q}),\, g\notin \G$ and $x\in X,\, x\notin F $.
There exists  a neighborhood $V$ of $x$, such that
$$
p(g\,p^{-1}(V))\cap V=\, \emptyset
$$
\end{lem}

\noindent {\it Proof.} Let $\G_0=\G\cap g^{-1}\G g$
and $X_0=\,\Gamma_0 \backslash Y$. For $x_0\in X_0$
the projections $p_1(x_0)=p(y)$ and $p_2(x_0)=p(g y)$
are independent of the representative $y\in Y$.
Moreover if $p_1(x_0)\notin F$ then $p_1(x_0)\neq p_2(x_0)$
by Lemma \ref{fixedpt}. By construction $\G_0$ is of finite index
$n$
in $\G$ and the fiber $p_1^{-1}(x)$
has at most $n$ elements. Let then $W\subset X_0$ be a compact
neighborhood of $p_1^{-1}(x)$ in $X_0$ such that $x\notin p_2(W)$.
For $z \in V\subset X$ sufficiently close to $x$ one has
$p_1^{-1}(z)\subset W$ and thus $p_2(p_1^{-1}(z))\subset V^c$
which gives the result. $\Box$

\smallskip

\n We can now prove the following.
\begin{prop} \label{measure}
 Let $\beta >0$ and $\varphi$ a  KMS$_\beta$
state  on $(A,\sigma_t)$. Then there exists a probability measure
$\mu$ on $X=\G\backslash Y$ such that
$$
\varphi(f)=\,\int_X\,f(1,x)\,d\mu(x)\qqq f\in A\,.
$$
\end{prop}

\noindent {\it Proof.} Let $g_0  \in
{\rm GL}_2^+ ({\Q})$ and
$f\in C_c(Z)$ such that
$$
f(g,y)=0 \qqq g\notin \G\,g_0\,\G\,.
$$
Since any element of $C_c(Z)$ is a finite linear combination of
such functions, it is enough to show that $\varphi(f)=0$ provided
$ g_0\notin \G$. Let $h_n\in C_c(X),\,0\leq h_n\leq 1$ with
support disjoint from $F$ and converging pointwise to $1$ in the
complement of $F$. Let $u_n \in A$ be given by
$$
u_n(1,y):=h_n(y)\,,\quad u_n(g,y)=0 \qqq g\notin \G\,.
$$

\smallskip

\n The formula
\begin{equation}\label{Phif}
\Phi(f)(g,\tau):= f(g,0,\tau) \ \ \ \forall f\in A
\end{equation}
defines a homomorphism of $(A,\sigma_t)$ to the $C^\ast$ dynamical
system $(B,\sigma_t)$ obtained by specialization to $\rho=0$, with
convolution product
$$ f_1 * f_2 (\rho,\tau)=\sum_h f_1(gh^{-1},h(\tau)) f_2(h,\tau), $$
where now we have no restriction on the summation, as in
\cite{CM10}.

\smallskip

\n For each $n\in\N^\ast$ we let
\begin{equation}\label{mubig}
\mu_{[n]}(g,y)=1 \;{\rm if} \;g\in \G.[n]\,, \quad
\mu_{[n]}(g,y)=0 \;{\rm if} \; g\notin \G.[n]\,.
\end{equation}
One has $\mu_{[n]}^\ast\,\mu_{[n]}=1$ and $\sigma_t(\mu_{[n]})=
n^{2it}\,\mu_{[n]} \qqq t\in \R$. Moreover, the range
$\pi(n)=\mu_{[n]}\,\mu_{[n]}^\ast$ of $\mu_{[n]}$ is the
characteristic function of the set of $\Q$-lattices that are
divisible by $n$, \ie those of the form $( \L ,n\, \phi)$.

\smallskip

\n Let $\nu_{[n]}=\Phi(\mu_{[n]})$. These are unitary multipliers
of $B$. Since they are eigenvectors for $\sigma_t$, the system
$(B,\sigma_t)$ has no non-zero KMS$_\beta$ positive functional.
This shows that the pushforward
of $\varphi$ by $\Phi$ vanishes and
by Proposition \ref{image}  that, with the
notation introduced above,
$$
\varphi(f)= \lim_n \,\varphi(f * u_n)\,.
$$
Thus, since $(f * u_n)(g,y)=\,f(g,y)\,h_n(y)$, we can assume that
$f(g,y)=0$ unless $p(y)\in K$, where $K\subset X$ is a compact
subset disjoint from $F$. Let $x\in K$ and $V$ as in lemma
\ref{corresp} and let $h\in C_c(V)$. Then, upon applying the
KMS$_\beta$ condition (\ref{eq4}) to the pair  $a,b$ with $a=f$
and
$$
b(1,y):=h(y)\,,\quad b(g,y)=0 \qqq g\notin \G\,,
$$
one gets $\varphi(b*f)=\varphi(f*b)$. One has
$(b*f)(g,y)=\,h(gy)\,f(g,y)$. Applying this to $f*b$ instead of
$f$ and using $h(gy)\,h(y)=0 \qqq y\in X$ we get
$\varphi(f*b^2)=0$ and $\varphi(f)=0$, using a partition of unity
on $K$. $\Box$

\smallskip

\n We let ${\rm Det}$ be the continuous map from
$M_2(\RR)$ to $\RR$ given by the determinant
$$
\begin{bmatrix} a &b \\ c &d \end{bmatrix} \in M_2(\RR)
\mapsto a\,d -\,b\,c \in \RR\,.
$$
For each $n\in\N^\ast$, the composition $\pi_n\circ {\rm Det}$
defines a projection $\pi'(n)$, which is the characteristic
function of the set of $\Q$-lattices whose determinant is
divisible by $n$. If a $\Q$-lattice is divisible by $n$ its
determinant is divisible by $n^2$ and one controls divisibility
using the following family of projections $\pi_p(k,l)$. Given a
prime $p$ and a pair $(k,l)$ of integers $k\leq l$, we let
\begin{equation}\label{projp}
\pi_p(k,l):=\,(\pi(p^k)-\pi(p^{k+1}))\,(\pi'(p^{k+l})-\pi'({p^{k+l+1}})).
\end{equation}
This corresponds, when working modulo $N=p^b, b>l$, to matrices in the
double class of
$$
\begin{bmatrix} a &0 \\ 0 &d \end{bmatrix}\,,\quad v_p(a)=k
\,,\quad v_p(d)=l\,,
$$
where $v_p$ is the $p$-adic valuation.

\smallskip

\begin{lem} \label{proba}
\begin{itemize}
\item Let $\varphi$ be a  KMS$_\beta$
state  on $(A,\sigma_t)$. Then, for any prime $p$ and pair $(k,l)$
of integers $k< l$, one has
$$
\varphi(\pi_p(k,l))=\,p^{-(k+l)\beta}\,p^{l-k}\,(1+p^{-1})
\,(1-p^{-\beta})\,(1-p^{1-\beta})\,
$$
while for $k=l$ one has
$$
\varphi(\pi_p(l,l))=\,p^{-2l\beta}
\,(1-p^{-\beta})\,(1-p^{1-\beta})\,.
$$

\item For distinct primes $p_j$ one has
$$
 \varphi(\prod
 \pi_{p_j}(k_j,l_j))=\,\prod\,\varphi(\pi_{p_j}(k_j,l_j)).
$$
\end{itemize}
\end{lem}

\noindent  {\it Proof.} For each $n\in\N^\ast$ we let $\nu_n\in
M(A)$ be given by
$$
 \nu_n(g,y)=1 \qqq g\in \G \begin{bmatrix} n &0 \\ 0
&1 \end{bmatrix}\,\G\,,\quad \nu_n(g,y)=0 \quad {\rm otherwise}\,.
$$
One has $\sigma_t(\nu_n)= n^{it}\,\nu_n \qqq t\in \R$.
The double class $\G \begin{bmatrix} n &0 \\ 0
&1 \end{bmatrix}\,\G$ is the union of the left $\G$-cosets
of the matrices
$
\begin{bmatrix} a &b \\ 0
&d \end{bmatrix}
$
where $a\,d=n$ and ${\rm gcd}(a,b,d)=1$. The number of these
left cosets
is $$\omega(n):=\, n \prod_{p \, {\rm prime,} \, p \mid n}\,(1 + p^{-1})$$
and
\begin{equation}\label{nunstar}
\nu_n^\ast * \nu_n(1,y)=\,\omega(n)\qqq y\in Y.
\end{equation}
 One has
$$
\nu_n * \nu_n^\ast(1,y)=\, \sum_{ h \in \G \backslash {\rm GL}_2^+ ({\Q}),\,hy\in Y}
\, \nu_n( h^{-1},h y)^2\,.
$$
With $y=(\rho,\tau)$, the r.h.s. is independent of $\tau$ and only depends
upon the $\SL_2(\Z/n\Z)-\GL_2(\Z/n\Z)$ double class of $\rho_n=\,p_n(\rho)
\in M_2(\Z/n\Z)$.

\n Let us  assume that $n=p^l$ is a prime
power. We can assume that $\rho_n=\,p_n(\rho)$ is of the form
$$
\rho_n=\,\begin{bmatrix} p^a &0 \\ 0
&p^b \end{bmatrix}\,,\quad 0\leq a\leq b\leq l\,.
$$
We need to count the number $\omega(a,b)$ of left $\G$-cosets $\G\,h_j$
in the double class $\G \begin{bmatrix} n^{-1} &0 \\ 0
&1 \end{bmatrix}\,\G$ such that $h_j y\in Y$ \ie
$h_j \rho\in M_2(\RR)$. A full set of representatives
of the double class is given by $h_j=(\a_j^t)^{-1}$ where
the $\a_j$ are $$\a_0=\begin{bmatrix} n &0 \\ 0
&1 \end{bmatrix}\,\quad\a(s)=\begin{bmatrix} 1 &s \\ 0
&n \end{bmatrix}\,,s\in\{0,1,...,n-1\}$$
 and for $x\in \{1,2,...,l-1\}$, $s\in \Z/ p^{l-x}\Z$ prime to $p$
$$
\a(x,s)=\begin{bmatrix} p^x &s \\ 0
&p^{l-x} \end{bmatrix}\,.
$$
The counting gives
\begin{itemize}
\item $\omega(a,b)=0$ if $b<l$.
\item $\omega(a,b)=p^a$ if $a< l,b\geq l$.
\item $\omega(a,b)=p^l(1 + p^{-1})$ if $a\geq l$.
\end{itemize}

\n Let $e_p(i,j),\, (i\leq j)$ be the projection corresponding to $a\geq i,\,b\geq j$.
Then for $i<j$ one has
\begin{equation}\label{projpie}
\pi_p(i,j)= e_p(i,j)-e_p(i+1,j)-e_p(i,j+1)+e_p(i+1,j+1)
\end{equation}
while
\begin{equation}\label{projpie1}
\pi_p(j,j)= e_p(j,j)-e_p(j,j+1).
\end{equation}
The computation above gives
\begin{equation}\label{nunustar}
\nu_n * \nu_n^\ast(1,y)=\,p^l(1 + p^{-1})e_p(l,l)+\,
 \sum_0^{l-1}\,p^k\,(e_p(k,l)-e_p(k+1,l))\,,
\end{equation}
where we omit the variable $y$ in the r.h.s.

\n Let $\varphi$ be a KMS$_\beta$ state, and $
\sigma(k,l):=\,\varphi(e_p(k,l)) $. Then, applying the KMS$_\beta$
condition to the pair $(\mu_{[p]}f,\mu_{[p]}^\ast)$ for
 $f\in C(X)$, one gets
$$
\sigma(k,l)=p^{-2k\beta}\,\sigma(0,l-k)\,.
$$
Let $\sigma(k)=\sigma(0,k)$. Upon applying the KMS$_\beta$
condition to $(\nu_n,\nu_n^\ast)$, one gets
$$
p^l(1 + p^{-1})p^{-l\beta}=\,p^l(1 + p^{-1})p^{-2l\beta}
+\sum_0^{l-1}\,p^k\,(p^{-2k\beta}\sigma(l-k)-p^{-2(k+1)\beta}
\sigma(l-k-1)\,.
$$
Since $\sigma(0)=1$, this determines the $\sigma(n)$ by induction
on $n$ and gives
$$ \sigma(n)= a \, p^{n(1-\beta)} + (1-a)\, p^{-2n\beta}, $$
with
$$ a= (1+p) \frac{p^\beta -1}{p^{1+\beta}-1} \,. $$
Combined with
(\ref{projpie}) and (\ref{projpie1}), this gives the required
formulas for $\varphi(\pi_p(k,l))$ and the first part of the lemma
follows.

\smallskip

\n To get the second part, one proceeds by induction on the number
$m$ of primes $p_j$. The function $f=\prod_1^{m-1}
\pi_{p_j}(k_j,l_j)$ fulfills
$$f(h\,y)=f(y)\qqq y\in Y\qqq h\in
\G \begin{bmatrix} n^{-1} &0 \\ 0 &1 \end{bmatrix}\,\G,$$ where
$n=p_m^l$. Thus, when applying the KMS$_\beta$ condition to
$(\nu_n\,f,\nu_n^\ast)$, the above computation applies with no
change to give the result. $\Box$\\

\n Let us now complete the proof of 3) of Theorem \ref{thmKMSst}.
Let $\varphi$ be a KMS$_\beta$ state.
Proposition \ref{measure} shows that there is a
probability measure $\mu$ on $X$
such that
$$
\varphi(f)=\,\int_X\,f(1,x)\,d\mu(x)\qqq f\in A\,.
$$
With $y=(\rho, \tau) \in X$, Lemma \ref{proba} shows that the
probability $\varphi(e_p(1,1))=\sigma(1,1)$
that a prime $p$ divides $\rho$ is $p^{-2\beta}$. Since
 the series $\sum
p^{-2\beta}$ converges ($\beta>\frac{1}{2}$ would suffice here),
it follows (\cf \cite{Rudin} Thm. 1.41)
that, for almost all $y\in X$, $\rho$ is only divisible
by a finite number of primes. Next, again by Lemma \ref{proba},
the probability that the determinant of $\rho$ is divisible by $p$
is
$$\varphi(e_p(0,1))=\sigma(1) =(1+p)p^{-\beta}-\,p^{1-2\beta}.
$$
For $\beta >2$ the corresponding series
$\sum\,((1+p)p^{-\beta}-\,p^{1-2\beta} )$ is convergent. Thus, we
conclude that with probability one
$$\rho_p \in \GL_2(\Z_p)\,,\quad {\rm for \;almost\; all}\; p.
$$
Moreover, since  $\sum\varphi(\pi_p(k,l))=1$, one gets with
probability one
$$\rho_p \in \GL_2(\Q_p)\qqq p \,.
$$
In other words, the measure $\mu$ gives measure one to finite
id\`eles. (Notice that finite id\`eles form a Borel subset which
is not closed.) However, when $\rho$ is a finite id\`ele the
corresponding $\Q$-lattice is commensurable to a unique invertible
$\Q$-lattice. Then the KMS$_\beta$ condition shows that the
measure $\mu$ is entirely determined by  its restriction to
invertible $\Q$-lattices, so that, for some probability measure
$\nu$,
$$
\varphi=\, \int \,\varphi_{\beta,l}\,\,d\nu(l).
$$
It follows that the Choquet simplex of extremal KMS$_\beta$ states
is the space of probability measures on the locally
compact space
$$
{\rm GL}_2 ({\Q})\backslash {\rm GL}_2 ({\A})/\C^\ast
$$
of invertible $\Q$-lattices \footnote{\cf \eg \cite{Mu} for the
standard identification of the set of invertible $\Q$-lattices
with the above double quotient}
and its extreme points are the $\varphi_{\beta,l}$. $\Box$\\

\smallskip
\n In fact Lemma \ref{proba}
 admits the following corollary:
\begin{cor}
For $\beta \leq 1$ there is no KMS$_\beta$ state
on $(A,\sigma_t)$.
\end{cor}

\noindent  {\it Proof.} Indeed the value of
$\varphi(\pi_p(k,l))$ provided by the lemma
is strictly negative for $\beta < 1$ and
vanishes for $\beta=1$.
In the latter case this shows that the measure
$\mu$ is supported by $\{0\}\times \H\subset Y$
and one checks that no such measure
fulfills the KMS condition for $\beta=1$. $\Box$

\smallskip

\n In fact the measure provided by Lemma \ref{proba} allows us to
construct a specific KMS$_\beta$ state on $(A,\sigma_t)$ for
$1<\beta\leq 2$. We shall analyze this range of values in Chapter
III in connection with the renormalization group.

\smallskip

\n To get some feeling about what happens when $\beta \rightarrow
2$ from above, we shall show that, on functions $f$ which are
independent of $\tau$, the states $\varphi_{\beta,l}$ converge
weakly to the KMS$_2$ state $\varphi$ of \eqref{KMS2state},
independently of the choice
of the invertible $\Q$-lattice $l$. Namely, we have
$$
\varphi_{\beta,l}(f) \rightarrow \int_{M_2(\RR)}\,f(a)\,da\, .
$$
Using the density of functions of the form $f\circ p_N$ among left
$\G$-invariant continuous functions on $M_2(\RR)$, this follows
from:
\begin{lem} \label{equi}
For $N\in\N$, let $\G(N)$ be the congruence subgroup of level $N$
and
$$
Z_\beta= \sum_{\Gamma(N) \backslash M_2(\Z)^+
} \,{\rm Det}(m)^{-\beta}\,.
 $$
When $\beta \rightarrow 2$ one has, for any function $f$ on
$M_2(\Z/N\Z)$,
$$
Z_\beta^{-1}\, \sum_{\Gamma(N) \backslash M_2(\Z)^+
} \,f(p_N(m))\,{\rm Det}(m)^{-\beta}\,\rightarrow
N^{-4}\,\sum_{M_2(\Z/N\Z)} \,f(a)  \,.
$$
\end{lem}

\noindent  {\it Proof.} For  $x\in M_2(\Z/N\Z)$ we let
$$h(x)=\,\lim_{\beta \rightarrow 2}
Z_\beta^{-1}\, \sum_{m\in\Gamma(N) \backslash M_2(\Z)^+,\,p_N(m)=x
} \,{\rm Det}(m)^{-\beta}\,
$$ be the limit of the above expression,
with $f$  the characteristic function of the subset
$\{x\}\subset M_2(\Z/N\Z)$.
We want to show that
\begin{equation} \label{limit}
h(x)= N^{-4}\qqq x\in M_2(\Z/N\Z)\,.
\end{equation}
Since $p_N$ is a surjection $SL_2(\Z)\to SL_2(\Z/N\Z)$ and $\G(N)$
a normal subgroup of $\G$, one gets
\begin{equation} \label{invar}
h(\g_1\,x\g_2)=\,h(x)
\qqq \g_j \in SL_2(\Z/N\Z)\,.
\end{equation}
Thus, to prove (\ref{limit}) we can assume that $x$ is a diagonal
matrix
$$x= \begin{bmatrix} n &0 \\ 0
&n\,\ell \end{bmatrix}\in M_2(\Z/N\Z)\,.$$ Dividing both $n$ and
$N$ by their g.c.d. $k$ does not affect the validity of
(\ref{limit}), since all $m\in\Gamma(N) \backslash M_2(\Z)^+$ with
$p_N(m)=x$ are of the form $k\,m'$, while ${\rm
Det}(m)^{-\beta}=k^{-2\beta}\,{\rm Det}(m')^{-\beta}$. This shows
that (\ref{limit}) holds for $n=0$ and allows us to assume that
$n$ is coprime to $N$. Let then $r$ be the g.c.d. of $\ell$ and
$N$. One can then assume that $\displaystyle x= \begin{bmatrix} n &0 \\
0 &n'\,r \end{bmatrix}$, with $r|N$ and with $n$ and $n'$ coprime
to $N$. Let $\Delta\subset SL_2(\Z/N\Z)$ be the diagonal subgroup.
The left coset $\Delta\, x\,\subset M_2(\Z/N\Z)$ only depends on
$r$ and the residue $\delta\in (\Z/N'\Z)^\ast$ of $n\,n'$ modulo
$N'=N/r$. It is the set of all diagonal matrices of the form
$$
y= \begin{bmatrix} n_1 &0 \\ 0 &n_2\,r \end{bmatrix}\, ,\quad
n_1\in (\Z/N\Z)^\ast \, ,\quad n_1\,n_2=\delta \;(N')\,.$$ Let
$\Gamma_\Delta(N)\subset\G$ be the inverse image of $\Delta$ by
$p_N$. By (\ref{invar}) $h$ is constant on $\Delta\,x$, hence
\begin{equation} \label{ressum}
h(x)=\,\lim_{\beta \rightarrow 2}Z_\beta^{-1}\,
\sum_{m\in\Gamma_\Delta(N) \backslash M_2(\Z)^+,\,p_N(m)\in \Delta \,x
} \,{\rm Det}(m)^{-\beta}\,.
\end{equation}
In each left coset $m\in\Gamma_\Delta(N) \backslash M_2(\Z)^+$
with $p_N(m)\in \Delta \,x$ one can find a unique triangular
matrix $ \displaystyle
\begin{bmatrix} a &b \\ 0
&d \end{bmatrix} $ with $a>0$ coprime to $N$, $d>0$ divisible by
$r$, $a\,d/r=\delta \,(N')$ and $b=N\,b'$ with $0\leq b'<d$. Thus,
we can rewrite (\ref{ressum}) as
\begin{equation} \label{ressum1}
h(x)=\,\,\lim_{\beta \rightarrow 2}\;Z_\beta^{-1}\, \sum_{Y }
\,d\,(ad)^{-\beta},
\end{equation}
where $Y$ is the set of pairs of
positive integers $(a,d)$ such that
$$
p_N(a)\in (\Z/N\Z)^\ast \,,\quad r|d
\,,\quad a\,d=\,r\,\delta \,(N).
$$
To prove (\ref{limit}) we assume first that $r<N$ and write
$N=N_1\,N_2$, where $N_1$ is coprime to $N'$ and $N_2$ has the
same prime factors as $N'$. One has $r=N_1\,r_2$, with $r_2|N_2$.
An element of $\Z/N_2\Z$ is invertible iff its image in $\Z/N'\Z$
is invertible. To prove (\ref{limit}) it is enough to show that,
for any of the $r_2$ lifts $\d_2\in \Z/N_2\Z$ of $\d$, one has
\begin{equation} \label{ressum2}
\lim_{\beta \rightarrow 2}\;Z_\beta^{-1}\, \sum_{Y' }
\,d\,(ad)^{-\beta}=N_1\,N^{-4},
\end{equation}
where $Y'$ is the set of pairs of
positive integers $(a,d)$ such that
$$p_N(a)\in (\Z/N\Z)^\ast \,,\quad
p_{N_2}(a\,d)=\d_2.$$ We let $1_N$ be the trivial Dirichlet
character modulo $N$. Then when $\Xc_{N_2}$ varies among Dirichlet
characters modulo $N_2$ one has
$$
\sum_{Y' } \,d\,(ad)^{-\beta}=\,\varphi(N_2)^{-1}\, \sum
\Xc_{N_2}(\d_2)^{-1}\,L(1_{N_1}\times \Xc_{N_2}, \beta) \,L(
\Xc_{N_2}, \beta-1),
$$
where $\varphi$ is the Euler totient function. Only the trivial
character $\Xc_{N_2}=1_{N_2}$ contributes to the limit
(\ref{ressum2}), since the other $L$-functions are regular at $1$.
Moreover, the residue of $L( 1_{N_2}, \beta-1)$ at $\beta=2$ is
equal to $\frac{\varphi(N_2)}{N_2}$ so that, when ${\beta
\rightarrow 2}$, we have
$$
\sum_{Y'
} \,d\,(ad)^{-\beta}\sim \,N_2^{-1}\, L(1_{N}, 2)\,(\beta -2)^{-1}\,.
$$
By construction one has
$$
Z_\beta \sim  |\G:\Gamma(N)|\,\zeta(2)\,(\beta -2)^{-1}\,,
$$
where the order of the quotient group $\G:\Gamma(N)$
is $N^3\,\prod_{p|N}\,(1-p^{-2})$ (\cite{Shimura}).
Since
$$
L(1_{N}, s)=\,\prod_{p|N}\,(1-p^{-s})\;\zeta(s)
$$
one gets (\ref{ressum2}). A similar argument handles the
case $r=N$. $\Box$ \\

\n The states $\varphi_{\beta,l}$ converge when $\beta \rightarrow
\infty$ and their limits restrict to $C_c(X)\subset A$ as
characters given by evaluation at $l$:
$$
\varphi_{\infty,l}(f)= f(l)\qqq f \in C_c(X).
$$
These characters are all distinct and we thus get
a bijection of the space
$$
{\rm GL}_2 ({\Q})\backslash {\rm GL}_2 ({\A})/\C^\ast
$$
of invertible $\Q$-lattices
with the space ${\mathcal E}_{\infty}$ of extremal
KMS$_\infty$ states.

\smallskip

\n We shall now describe the natural symmetry group of the above
system, from an action of the quotient group $S$
$$
S:=\,\Q^\ast\backslash {\rm GL}_2 ({\A_f})
$$
as symmetries of our dynamical system.
Here the finite ad\`elic group of $\GL_2$ is given by
$$ \GL_2(\A_f) = \prod_{{\rm res}} \GL_2(\Q_p), $$
where in the restricted product the $p$-component lies in
$\GL_2(\Z_p)$ for all but finitely many $p$'s. It satisfies
$$ \GL_2(\A_f) = \GL_2^+(\Q) \,\GL_2(\RR). $$
The action of the subgroup
$$
{\rm GL}_2 ({\RR})\subset S
$$
is defined in a straightforward manner using the following right
action of ${\rm GL}_2 ({\RR})$ on $\Q$-lattices:
$$
(\L,\phi). \g\,= (\L,\phi \circ \g)\qqq \g \in {\rm GL}_2 ({\RR})\,.
$$
By construction this action preserves the commensurability
relation for pairs of $\Q$-lattices and preserves the value of the
ratio of covolumes for such pairs. We can view it as the action
$$
(\rho,\tau).\g\,=\,(\rho \circ \g, \tau)
$$
of ${\rm GL}_2 ({\RR})$ on $Y=M_2(\RR)\times\H$, which commutes
with the left action of ${\rm GL}^+_2 ({\Q})$. Thus, this action
defines automorphisms of the dynamical system $(A,\sigma_t)$ by
$$
\t_\g(f)(g,y):=\, f(g,y.\g)\qqq f\in \Ac\,,\g\in {\rm GL}_2 ({\RR})\,,
$$
and one has
$$
\t_{\g_1}\;\t_{\g_2}\;=\;\t_{\g_1\g_2}\qqq \g_j \in {\rm GL}_2 ({\RR})\,.
$$

\smallskip

\n The complementary action of ${\rm GL}^+_2 ({\Q})$ is more
subtle and is given by endomorphisms of the dynamical system
$(A,\sigma_t)$, following Definition \ref{end}.

\smallskip

\n For $m \in M_2 ({\mathbb Z})^+$, let $\tilde{m}={\rm Det}(m)\;
m^{-1}\in M_2 ({\mathbb Z})^+$. The range $R_m$ of the map $\rho
\rightarrow \rho\,\tilde{m}$  only depends on $L = m({\mathbb
Z}^2)$. Indeed if $m_j \in M_2 ({\mathbb Z})^+$ fulfill
$m_1({\mathbb Z}^2) = m_2({\mathbb Z}^2)$ then  $m_2 = m_1\,\gamma
\, $ for some $\gamma \in \Gamma$, hence $ M_2 (R)\,\tilde{m_1} =
M_2 (R)\,\tilde{m_2}$. Let then
\begin{equation}
\label{latticelabel}
 e_L \in C(X)
\end{equation}
be the characteristic function of $\Gamma \backslash (R_m \times
\H) \subset \Gamma \backslash (M_2 (R) \times \H) $, for any $m$
such that  $ m({\mathbb Z}^2)=L$. Equivalently, it is the
characteristic function of the open and closed subset $E_L \subset
X$ of $\Q$-lattices of the form $(\L, \phi \circ \tilde{m})$. One
has
\begin{equation}
\label{rho}
e_L \, e_{L'} = e_{L \cap L'} \, , \qquad e_{{\mathbb Z}^2} = 1 \,
.\end{equation}
For $l=( \L , \phi)\in E_L\subset X$ and $m \in M_2 ({\mathbb Z})^+$,
$ m({\mathbb Z}^2)=L$
we let
$$
l\circ \tilde{m}^{-1}:= ( \L , \phi\circ \tilde{m}^{-1})\in X.
$$
This map preserves commensurability of $\Q$-lattices.
On $Y$ it is given by
$$
(\rho,\tau)\circ \tilde{m}^{-1}:=\,(\rho \circ \tilde{m}^{-1}, \tau)\qqq
(\rho,\tau)\in R_m \times \H
$$
and it commutes with the left action of ${\rm GL}^+_2 ({\Q})$. The
formula
\begin{equation}
\label{theta} \t_m(f)(g,y):= f(g,y\circ \tilde{m}^{-1})\qqq y \in
R_m \times \H,
\end{equation}
extended by $\t_m(f)(g,y)=0$ for $y \notin R_m \times \H$, defines
an endomorphism $\t_m$ of $A$ that commutes with the time
evolution $\sigma_t$. Notice that $\t_m(1)=e_L \in M(A)$ is a
multiplier of $A$ and that $\t_m$ lands in the reduced algebra
$A_{e_L}$, so that (\ref{theta}) is unambiguous. Thus one obtains
an action of the semigroup $M_2 ({\mathbb Z})^+$ by endomorphisms
of the dynamical system $(A,\sigma_t)$, fulfilling Definition
\ref{end}.

\begin{prop} The above actions of the group ${\rm GL}_2 ({\RR})\subset S$
and of the semigroup $M_2 ({\mathbb Z})^+\subset S$ assemble to an
action of the group $S=\,\Q^\ast\backslash {\rm GL}_2 ({\A_f})$ as
symmetries of the dynamical system $(A,\sigma_t)$.
\end{prop}

\noindent {\it Proof.}  The construction above applies to give an
action by endomorphisms of the semigroup ${\rm GL}_2 ({\A_f})\cap
M_2(\RR)$, which contains both ${\rm GL}_2 ({\RR})$ and $M_2
({\mathbb Z})^+$. It remains to show that the sub-semigroup
$\N^\times \subset M_2 ({\mathbb Z})^+$ acts by inner
endomorphisms of $(A,\sigma_t)$. Indeed for any $n\in
\N^\ast$, the endomorphism $\t_{[n]}$ (where $[n]= \begin{bmatrix}
n &0 \\ 0 &n \end{bmatrix} \in M_2 ({\mathbb Z})^+$) is inner and
implemented by the
 multiplier $\mu_{[n]}\in M(A)$ which was defined in
(\ref{mubig}) above \ie one has
$$
\t_{n}(f)= \mu_{[n]}\,f\,\mu_{[n]}^\ast \qqq f\in A\,.
$$
$\Box$

\section{The subalgebra $\Ac_\Q$ and the Modular Field}

\n The strategy outlined in \S \ref{eisenstein} allows us to find,
using Eisenstein series, a suitable {\em arithmetic} subalgebra
$\Ac_\Q$ of the algebra of unbounded multipliers of the basic
Hecke $C^\ast$-algebra $A$ of the previous section. The extremal
KMS$_\infty$ states $\varphi \in {\mathcal E}_\infty $  extend to
$\Ac_\Q$ and the image $\varphi(\Ac_\Q)$ generates, in the generic
case, a specialization $F_\varphi\subset\C$ of the modular field
$F$. The state $\varphi$ will then intertwine the symmetry group
$S$ of the system $(A,\sigma_t)$ with the Galois group of the
modular field \ie we shall show that there exists an isomorphism
$\theta$ of $S$ with $\Gal (F_\varphi/\Q)$ such that
\begin{equation}\label{CFTiso2}
 \alpha \circ \varphi = \varphi \circ \theta^{-1}(\alpha)
\qqq \alpha \in \Gal (F_\varphi/\Q)\,.
\end{equation}

\smallskip

\n Let us first define $\Ac_\Q$ directly without any reference to
Eisenstein series and check directly its algebraic
properties. We let $Z\subset \G \backslash {\rm GL}_2^+
({\Q})\times_\Gamma Y$ be as above and $f\in C(Z)$ be a function
with {\em finite support} in the variable  $g\in\G \backslash {\rm
GL}_2^+ ({\Q})$. Such an $f$ defines an unbounded multiplier of
the $C^\ast$-algebra $A$ with the product given as above by
$$
(f_1 * f_2)(g,y) \,:= \, \sum_{ h \in \G \backslash {\rm GL}_2^+
({\Q}),hy\in Y} \,f_1(g h^{-1},h y)  \; f_2(h,y).
$$
One has $Y = M_2(\RR)\times \H$ and we write $f(g,y)=f(g,\rho,z)$,
with $(g,\rho,z) \in {\rm GL}_2^+ ({\Q}) \times M_2(\RR)\times \H
$. In order to define the {\em arithmetic} elements $f \in \Ac_\Q$
we first look at the way $f$ depends on $\rho \in M_2(\RR)$. Let
as above $p_N:M_2(\RR)\to M_2(\Z/N\Z)$ be the canonical
projection. It is a ring homomorphism. We say that $f$ has level
$N$ iff $f(g,\rho,z)$ only depends upon $(g,p_N(\rho),z)\in {\rm
GL}_2^+ ({\Q}) \times M_2(\Z/N\Z)\times \H $. Then specifying $f$
amounts to assigning the finitely many continuous functions
$f_{g,m} \in C(\H)$ with $m \in M_2(\Z/N\Z)$ and
$$
f(g,\rho,z)=\, f_{g,p_N(\rho)}(z).
$$
The invariance condition
\begin{equation}\label{modinv}
\quad f(g \,\g,y)=f(g,\g\,y)\qqq \g\in\G\,,\;g\in {\rm GL}_2^+
({\Q})\,,\; y\in Y
\end{equation}
then shows that
$$
f_{g,m}\vert \g=\, f_{g,m} \qqq \g \in \G(N)\cap g^{-1}\G g,
$$
with standard notations for congruence subgroups and for the slash
operation in weight $0$ (\cf \eqref{act1}).

\smallskip

\n We denote by $F$ the field of modular functions which are
rational over $\Q^{ab}$, \ie the union of the fields $F_N$ of
modular functions of level $N$ rational over $\Q(e^{2\pi i/N})$.
Its elements are modular functions $h(\tau)$ whose
$q^{\frac{1}{N}}$-expansion has all its coefficients in
$\Q(e^{2\pi i/N})$ (\cf \cite{Shimura}).

\smallskip

\n The first requirement for arithmetic elements is that
\begin{equation}\label{arithm1}
f_{g,m}\in F \ \ \ \forall (g,m).
\end{equation}

\smallskip

\n This condition alone, however, is not sufficient. In fact, the
modular field $F_N$  of level $N$  contains (\cf \cite{Shimura}) a
primitive $N$-th root of $1$. Thus, the condition \eqref{arithm1}
alone allows the algebra $\Ac_\Q$ to contain the cyclotomic field
$\Q^{ab}\subset \C$, but this would prevent the existence of
``fabulous states", because the ``fabulous'' property would not be
compatible with $\C$-linearity. We shall then impose an additional
condition, which forces the spectrum of the corresponding elements
of $\Ac_\Q$ to contain all Galois conjugates of such a root, so
that no such element can be a scalar. This is, in effect, a
consistency condition on the roots of unity that appear in the
coefficients of the $q$-series, when $\rho$ is multiplied on the
left by a diagonal matrix.

\smallskip

\n Consider elements $g\in \GL_2^+(\Q)$ and $\alpha\in
\GL_2(\Z/N\Z)$, respectively of the form
\begin{equation}\label{galpha}
 g = r \left[ \begin{array}{cc} n & 0 \\ 0 & 1
\end{array} \right] \ \ \ \text{ and } \ \ \
\alpha=\left[ \begin{array}{cc} k & 0 \\ 0 & 1
\end{array} \right],
\end{equation}
with $k$ prime to $N$ and $n| N$.

\begin{defn} \label{arith}
We shall say that $f$ of level $N_0$ is {\em arithmetic} ($f \in
\Ac_\Q$) iff for any multiple $N$ of $N_0$ and any pair
$(g,\alpha)$ as in \eqref{galpha} we have $f_{g,m}\in F_N$ for all
$m\in M_2(\Z/N\Z)$ and the $q$-series of $f_{g,\a\,m}$ is obtained
from the $q$-series for $f_{g,m}$ by raising to the power $k$ the
roots of unity that appear as coefficients.
\end{defn}

\smallskip

\n The arithmetic subalgebra $\cA_\Q$ enriches the structure of
the noncommutative space to that of a ``noncommutative arithmetic
variety''. As we shall prove in Theorem \ref{main}, a generic
ground state $\varphi$ of the system, when evaluated on $\cA_\Q$
generates an embedded copy $F_\varphi$ of the modular field in
$\C$. Moreover, there exists a unique isomorphism
$\theta=\theta_\varphi$ of the symmetry group $S$ of the system
with $\Gal (F_\varphi/\Q)$, such that
$$ \theta(\sigma) \circ \varphi = \varphi \circ \sigma
\qqq \sigma \in S \,. $$

\smallskip

\n A first step towards this result is to show that the
arithmeticity condition is equivalent to a covariance property
under left multiplication of $\rho$ by elements $\alpha\in
\GL_2(\RR)$, in terms of Galois automorphisms. The condition is
always satisfied for $\alpha\in \SL_2(\RR)$.

\smallskip

\n For each $g \in \GL_2(\A_f)$, we let ${\rm Gal}(g)\in {\rm
Aut}(F)$ be its natural action on $F$, written in a covariant way
so that
$$
{\rm Gal}(g_1\, g_2)=\,{\rm Gal}(g_1)\circ {\rm Gal}(g_2).
$$
With the standard contravariant notation 
$f\mapsto f^g$ (\cf \eg \cite{Lang}) we let, for all
$f\in F$,
\begin{equation}\label{galoisact}
{\rm Gal}(g)(f):= f^{\tilde{g}} \,,\quad \tilde{g}={\rm Det}(g)\; g^{-1}\,.
\end{equation}

\begin{lem}\label{arith-lem}
For any $\a \in \SL_2(\RR)$ one has
$$ f_{g,\a\,m}=\, {\rm Gal}(\a) f_{g',m}, $$
where $g\,\a= \a'\,g'$ is the decomposition of $g\,\a$ as a
product in $\GL_2(\RR).\GL_2(\Q)$.
\end{lem}

\n {\em Proof.} Notice that the decomposition $\a'\,g'$ is not
unique, but the left invariance $$f(\g\, g',\rho,\tau)=
f(g',\rho,\tau) \qqq \g \in \G$$ shows that the above condition is
well defined. Let $p_N: M_2(\RR)\to M_2(\Z/N\Z)$ be the
projection. Then $f_{g,\a\,m}=f_{g,p_N(\a)p_N(m)}$, for $f$ of
level $N$. Let $\gamma \in \G$ be such that $p_N(\gamma)=p_N(\a)$.
Then
$$ f_{g,\a m}(\tau)= f(g,\gamma m,\tau)= f(g\gamma, m,
\gamma^{-1}(\tau)). $$ Thus, for $g'=g\gamma$, one obtains the
required condition. $\Box$

\smallskip

\begin{lem} \label{arith-lem2}
A function $f$ is in $\cA_\Q$ iff condition \eqref{arithm1} is
satisfied and
\begin{equation}\label{galois}
f_{g,\a\,m}=\, {\rm Gal}(\a) f_{g',m}, \ \ \ \ \forall \a \in
\GL_2(\RR),
\end{equation}
where $g\,\a= \a'\,g'$ is the decomposition of $g\,\a$ as a
product in $\GL_2(\RR).\GL_2(\Q)$.
\end{lem}

\n {\em Proof.} By Lemma \ref{arith-lem}, the only nontrivial part
of the covariance condition \eqref{galois} is the case of diagonal
matrices $\delta=\left[
\begin{array}{cc} u & 0 \\ 0 & 1
\end{array} \right]$ with 
$u\in \GL_1(\RR)$.

\n To prove \eqref{galois} we can assume that $g=\,g_0\,\g$
with $g_0 $ diagonal as in (\ref{galpha})
and $\g\in \G$. Let then $\g\,\a=\,\delta\,\a_1$
with $\d$ as above and $\a_1 \in \SL_2(\RR)$.
One has 
$$
f_{g_0,\delta \a_1 m}= {\rm Gal}(\delta)\,f_{g_0, \a_1 m}
$$
 by Definition \ref{arith} since  ${\rm Gal}(\delta)$
 is given by raising the roots of unity that appear as coefficients of
the $q$-expansion to the power $k$ where 
$u$ is the residue of $k$ modulo $N$ (\cf \cite{Shimura} (6.2.1)
p.141). One then has
$$
f_{g,\a\,m}=\, {\rm Gal}(\g^{-1}) f_{g_0,\g\a m}={\rm Gal}(\g^{-1}\delta)\,f_{g_0, \a_1 m}
$$
and by Lemma \ref{arith-lem}, with $g_0 \a_1=\a_1'g_0'$ we get 
$$
f_{g,\a\,m}=\,{\rm Gal}(\g^{-1}\delta\a_1)\,f_{g_0',  m}=\,{\rm Gal}(\a)\,f_{g_0',  m}
$$
Moreover
$$
g\,\a=\,g_0\,\g\,\a=\,g_0\,\delta\,\a_1=\,\delta\,g_0\,\a_1=\,\delta\,\a_1'\,g_0'
$$
which shows that $g'=\,g_0'$

\n One checks similarly that the converse holds. $\Box$

\smallskip

\begin{prop}\label{subarit}
$\Ac_\Q$ is a subalgebra of the algebra of unbounded multipliers
of $A$, globally invariant under the action of the symmetry group
$S$.
\end{prop}

\noindent {\it Proof.} For each {\em generic} value of $\tau \in
\H$ the evaluation map
$$
h \in F \mapsto I_\tau(h)=h(\tau) \in \C
$$
gives an isomorphism of $F$ with a subfield $F_\tau \subset \C$
and a corresponding action ${\rm Gal}_\tau$ of $\GL_2(\A_f)$ by
automorphisms of $F_\tau$, such that
\begin{equation}\label{Galtau}
{\rm Gal}_\tau(g) (I_\tau(h))= I_\tau({\rm Gal}(g)(h)).
\end{equation}
We first rewrite the product as
$$
(f_1 * f_2)(g,\rho, \tau) \,= \, \sum_{ g_1 \in \G \backslash
{\rm GL}_2^+ ({\Q}),g_1\rho \in M_2(\RR)}
\, f_1(g g_1^{-1},g_1 \rho,g_1(\tau)) \;f_2(g_1,\rho, \tau) \,.
$$
\n The proof  that $(f_1 * f_2)_{g,m}\in F$ is the same as in
Proposition 2 of (\cite{CM10}). It remains to be shown that
condition \eqref{galois} is stable under convolution. Thus we let
$\a \in \GL_2(\RR)$ and we want to show that $f_1 * f_2$ fulfills
\eqref{galois}. We let $g'\in \GL_2(\Q)$ and $\b \in  \GL_2(\RR)$
with $g\,\a=\, \b\,g'$. By definition, one has
$$
(f_1 * f_2)(g,\a \,\rho, \tau) =  \sum_{ g_1 \in \G \backslash
{\rm GL}_2^+ ({\Q}), g_1\a\rho \in M_2(\RR)} \,f_1(g g_1^{-1},g_1
\,\a \,\rho,g_1(\tau))  \, f_2(g_1,\a \,\rho, \tau).
$$
We let $g_1 \a= \a' g'_1$ be the decomposition of $g_1 \a$,  and
use \ref{galois} to write the r.h.s. as
$$
(f_1 * f_2)(g,\a \,\rho, \tau) \,= \, \sum \, f_1(g g_1^{-1},g_1
\,\a \,\rho,g_1(\tau))\; {\rm Gal}_\tau(\a)( f_2(g'_1,\,\rho,
\tau) )\,,
$$
with $\Gal_\tau$ as in \eqref{Galtau}. The result then follows
from the equality
\begin{equation}\label{pass}
 f_1(g g_1^{-1},g_1 \,\a \,\rho,g_1(\tau))\,=
 {\rm Gal}_\tau(\a)( f_1(g^{'} g_1^{'\,-1},g'_1
 \,\rho,g'_1(\tau))),
\end{equation}
which we now prove. The equality $g_1 \a= \a' g'_1$ together with
\eqref{galois} shows that
$$
f_1(g g_1^{-1},g_1 \,\a \,\rho,g_1(\tau))\,=
 {\rm Gal}_{g_1(\tau)}(\a')(f_1(g^{'} g_1^{'\,-1},g'_1  \,\rho,g_1(\tau)))\,,
$$
using $\quad g g_1^{-1}\,\a'=\,g g_1^{-1}(g_1 \a\,g_1^{'\,-1} )=\,g \a \,g_1^{'\,-1}
= \b \,g'\,\,g_1^{'\,-1}$.

\smallskip
\n For any $h\in F$, one has
$$
I_{g_1(\tau)}({\rm Gal}(\a')h)=\,I_{g_1(\tau)}({\rm Gal}(g_1)
{\rm Gal}(\a){\rm Gal}(g_1^{'\,-1})(h))
$$
and, by construction of the Galois action \cite{Shimura},
$$
I_{g_1(\tau)}\circ{\rm Gal}(g_1)=\, I_\tau,
$$
so that in fact
$$
 {\rm Gal}_{g_1(\tau)}(\a')I_{g_1(\tau)}(h)=\,
{\rm Gal}_\tau(\a)I_{g'_1(\tau)}(h))\,.
$$
This proves \eqref{pass} and it shows that $\Ac_\Q$ is a
subalgebra of the algebra of unbounded multipliers of $A$. To
prove the invariance under $S$ is straightforward, since the
endomorphisms are all acting on the $\rho$ variable by right
multiplication, which does not interfere with condition
\eqref{galois}. $\Box$

\smallskip

\n In fact, modulo the nuance between ``forms" and functions, the
above algebra $\Ac_\Q$ is intimately related to the modular Hecke
algebra of \cite{CM10}.

\smallskip

\n We can now state the main result extending
Theorem \ref{BCthm}
to the two dimensional case.

\begin{thm}\label{main}
Let $l=(\rho,\tau)$ be a generic invertible
$\Q$-lattice and $\varphi_l\in \mathcal{E}_\infty$ be
the corresponding KMS$_\infty$
state. The image $\varphi_l(\Ac_\Q)\subset \C$ generates
the specialization $F_\tau \subset \C$
of the modular field $F$ obtained for
the modulus $\tau$.
The action of the symmetry group
$S$ of the dynamical system $(A,\sigma_t)$
is intertwined by $\varphi$
with the Galois group  of the modular field
$F_\tau$ by the formula
$$
\varphi \circ \a=\, {\rm Gal}_\tau(\rho\,\a\,\rho^{-1})\circ
\varphi.
$$
\end{thm}

\smallskip

\noindent {\it Proof.} We first need to exhibit
enough elements of $\Ac_\Q$. Let us first deal
with functions $f(g,\rho,\tau)$ which
vanish except when $g\in \G$. By construction these
are functions on the space
 $X$
of $\Q$-lattices
\begin{equation}
\label{X-spaceQlat}
X=\,(\hbox{Space of ${\mathbb Q}$-lattices in
${\mathbb C}$}) / {\mathbb C}^* \,\sim
\Gamma \backslash (M_2 (R) \times \H) \, .
\end{equation}
To obtain such elements of $\Ac_\Q$ we start with Eisenstein
series and view them as functions on the space of $\Q$-lattices.
Recall that to a pair $(\rho , \tau) \in Y$ we associate the
${\mathbb Q}$-lattice $(\L,\phi) =\t(\rho,\tau)$ by
\begin{equation}
\label{qlat} \Lambda = {\mathbb Z} + \, \tau \,{\mathbb Z} \, , \
\phi (a) = \rho_1 (a) - \tau \rho_2 (a) \in {\mathbb Q}\Lambda /
\Lambda,
\end{equation}
where $\rho_j (a) = \sum \rho_{jk} (a_k) \in {\mathbb Q} /
{\mathbb Z}$, for $a = (a_1 , a_2) \in ({\mathbb Q} / {\mathbb
Z})^2$. The Eisenstein series are given by
\begin{equation}
\label{eisen} E_{2k,a} (\rho , \tau) = \pi^{-2k} \sum_{y \in
\Lambda + \phi (a)} y^{-2k}.
\end{equation}
This is undefined when $\phi (a)\in \L$, but we shall easily deal
with that point below. For $k=1$ we let
\begin{equation}
\label{Xa-rhotau} X_a (\rho , \tau) = \pi^{-2} \left( \sum_{y \in
\Lambda + \phi (a)} y^{-2} - \underset{y \in
\Lambda}{\sum\nolimits'} \ y^{-2} \right)
\end{equation}
when $\phi (a) \notin \Lambda$ and $X_a (\rho , \tau) = 0$ if
$\phi (a) \in \Lambda$. This is just the evaluation of the
Weierstrass $\wp$-function on $\phi (a)$.

\smallskip

\n For $\gamma \in \Gamma = {\rm SL}_2 ({\mathbb Z})$ we have
$\quad\displaystyle X_a (\gamma \,\rho , \gamma \,\tau) = (c\tau +
d)^2 \, X_a (\rho , \tau), \quad$ which shows that the function
$c(\tau) X_a$ is $\Gamma$-invariant on $Y$, where
\begin{equation}
\label{fricke}
c(\tau) = -2^7 \, 3^5 \, \frac{g_2 g_3}{\Delta}
\end{equation}
 has weight $-2$ and no pole in $\H$.
We use $c$  as we used the covolume in the 1-dimensional case, to
pass to modular functions. This corresponds in weight $2$ to
passing from division values of the  Weierstrass $\wp$-function to
the Fricke functions (\cf \cite{Lang} \S 6.2)
\begin{equation}\label{Fricke}
f_v(\tau)=- 2^7 3^5 \;\frac{g_2\,g_3}{\Delta} \wp (\lambda(v,\tau)),
\end{equation}
where $v=(v_1,v_2)\in (\mathbb{Q}\slash\mathbb{Z})^2$
and $\lambda(v,\tau):= v_1 \,\tau + v_2$.
Here
 $\;g_2$, $g_3\;$ are the
coefficients giving the elliptic curve
$E_\tau = \C/\Lambda$
 in Weierstrass form,
$$
 y^2 = 4x^3 - g_2 x - g_3,
$$
with discriminant $\Delta= g_2^3 - 27 g_3^2$.
One has (up to powers of $\pi$)
$$
g_2= 60 \,e_4\,,\quad g_3= 140 \,e_6,
$$
where one defines the standard modular forms of even weight $k\in
2\,\Nb$ as
$$
e_k(\L):= \pi^{-k}\,\sum_{y\in  \Lambda \backslash \{0\}}\,y^{-k}
$$
with $q$-expansion ($q=e^{2\pi i \tau}$)
$$
e_k = \frac{2^k}{k!} B_{\frac{k}{2}} + (-1)^{k/2}
\frac{2^{k+1}}{(k-1)!}\,\sum_1^{\infty} \sigma_{k-1}(N)\, q^N,
$$
where the $B_n$ are the Bernoulli numbers and $\sigma_{n}(N)$ is
the sum of $d^n$ over the divisors $d$ of $N$. The $e_{2n}$ for
$n\geq 2$ are in the ring $\Qb[e_4,e_6]$ (cf \cite{WeilEll})
thanks to the relation
\begin{equation}\label{wealg0}
\frac{1}{3}(m-3)(4m^2-1)e_{2m}=\sum_2^{m-2}\,(2r-1)(2m-2r-1)e_{2r}e_{2m-2r}.
\end{equation}

\smallskip

\n Notice that $\Xc_a:=\,c\,X_a \in C (M_2 (R) \times \H) = C(Y)$
 is a continuous function
on $Y$. The continuity of $X_a$ as a function of $\rho$ comes from the fact
that it only involves the restriction
$\rho_N \in M_2 ({\mathbb Z} / N {\mathbb Z})$ of $\rho$
to $N$-torsion elements $a$
with $Na = 0$.

\n We view $\Xc_a$ as a function on $Z$ by
$$
\Xc_a(\g,\rho, \tau):= \Xc_a(\rho, \tau) \qqq \g \in \G,
$$
while it vanishes for $\g \notin \G$. Let us show that $\Xc_a \in
\Ac_\Q$. Since the Fricke functions belong to the modular field
$F$ we only need to check \eqref{galois}. For $\a \in \GL_2(\RR)$
and generic $\tau$ we want to show that
$$
\Xc_a(\a\,\rho,\tau)=\, {\rm Gal}_\tau(\a) \Xc_a(\rho,\tau).
$$
If $\rho(a)=0$ both sides vanish, otherwise they are both given by
Fricke functions $f_v,\,f_{v'}$, corresponding respectively to the
labels (using \eqref{qlat})
$$
v = s\, \a\,\rho(a)\,,\quad v'= s\,\rho(a)\,,\quad s=
\begin{bmatrix} 0 &-1 \\ 1 &0 \end{bmatrix}.
$$
Thus, $v'=s\,\a \,s^{-1}(v)=(s^{-1}\,\a^t\, s)^t(v)$ and the
result follows from \eqref{galoisact} and the equality
$$
\tilde{\a}=\,{\rm Det}(\a)\; \a^{-1}=\,s^{-1}\,\a^t\, s,
$$
with the Galois group $\GL_2(\Z/n\Z)/\pm 1$  of the modular field
$F_n$ over $\Q(j)$ acting on the Fricke functions by permutation
of their labels:
$$
f_v^{\sigma(u)}= f_{u^tv} \qqq u\in \GL_2(\Z/n\Z).
$$
This shows that $\Xc_a\in \Ac_\Q$ and it suffices to show that,
with the notation of the theorem, $\varphi_l(\Ac_\Q)$ generates
$F_\tau$, since the modular field $F$ is the field generated over
$\Q$ by all the Fricke functions. It already contains $\Q(j)$ at
level $2$ and it contains in fact $\Q^{ab}(j)$.

\smallskip
\n Let us now display elements $T_{r_1,r_2}\in\Ac_\Q$,
$r_j\in\Q_+^*, r_1|r_2$, associated to the classical
Hecke correspondences. We let $C_{r_1,r_2}
\subset \G\backslash
\GL_2(\Q)^+$ be the finite subset given by the double class
of
$
\begin{bmatrix} r_1 &0 \\ 0 &r_2 \end{bmatrix} $ in $\G\backslash
\GL_2(\Q)^+/\G$. We then define
$$
T_{r_1,r_2}(g,\rho,\tau)= 1 \quad {\rm if}\quad
 g\in C_{r_1,r_2}\,,\; g\rho\in M_2(\RR)\,,\quad
T_{r_1,r_2}(g,\rho,\tau)= 0 \quad {\rm otherwise}\,.
$$
One needs to check (\ref{galois}), but if $g\,\a= \a'\,g'$ is the
decomposition of $g\,\a$ as a product in $\GL_2(\RR).\GL_2(\Q)^+$,
then $g'$ belongs to the double coset of $g \in \G\backslash
\GL_2(\Q)^+/\G$, which gives the required invariance. It is not
true that the $T_{r_1,r_2}\in\Ac_\Q$ fulfill the relations of the
Hecke algebra $\Hc(\GL_2(\Q)^+,\G)$ of double cosets, but this
holds when $r_1,r_2$ are restricted to vary among positive {\em
integers}. To see this one checks that the map
$$
\tau(f)(g,y):=\,f(g)\quad {\rm if}\quad
 g\in M_2(\Z)^+\,,\quad
\tau(f)(g,y)= 0 \quad {\rm otherwise}
$$
defines an isomorphism
\begin{equation}\label{heckeinclusion}
\tau: \Hc(M_2(\Z)^+,\G)\to \Ac_\Q
\end{equation}
of the standard Hecke algebra $\Hc(M_2(\Z)^+,\G)$ of
$\G$-biinvariant functions (with $\G$-finite support) on
$M_2(\Z)^+$ with a subalgebra $\Hc\subset  \Ac_\Q $. Notice that
it is only because the condition $h\,y \in Y$ of definition
(\ref{convpro}) is now automatically satisfied that $\tau$ is a
homomorphism.
\smallskip

\n Let us now show the intertwining equality
\begin{equation}\label{intertwin}
\varphi_l \circ \a=\, {\rm Gal}_\tau(\rho\,\a\,\rho^{-1})\circ \varphi_l\qqq
\a \in S\,.
\end{equation}
One has
$$
\varphi_l(f)=f(1,\rho,\tau)\qqq f \in \Ac_\Q
$$
It is enough to prove \eqref{intertwin} for $\a \in \GL_2(\RR)$
and for $\a \in \GL_2(\Q)$.

\smallskip
\n For $\a \in \GL_2(\RR)$, the state $\varphi_l \circ \a$ is
given simply by
$$
(\varphi_l \circ \a)(f)=f(1,\rho\,\a,\tau)\qqq f \in \Ac_\Q\,,
$$
and using \eqref{galois} one gets \eqref{intertwin} in that case.

\smallskip
\n Let $m\in M_2^+(\Z)$. Then the state $\varphi_l \circ m$ is
more tricky to obtain, since it is not the straight composition
but the $0$-temperature limit of the states obtained by
composition of the KMS$_\beta$ state $\varphi_{l,\beta}$ with the
endomorphism $\t_m$ defined in (\ref{theta}). Indeed, the range of
$\theta_m$ is the reduced algebra by the projection $e_L$, with
$L=m(\Z^2)$, on which any of the zero temperature states vanishes
identically.

\smallskip

\n Let us first show that for finite $\beta$ we have
\begin{equation}\label{comp}
\varphi_{l,\beta}\circ
\t_m=\,\varphi_{l,\beta}(e_L)\,\varphi_{l',\beta},
\end{equation}
where $L=m(\Z^2)$ and $l'$ is given by
\begin{equation}\label{posbeta}
l'=(\rho',\,m^{'\,-1}(\tau))\,,\quad \rho\,
m = m'\, \rho'\in M_2^+(\Z).\GL_2(\RR)\,.
\end{equation}
By \eqref{kmsexplicit} we have
$$
\varphi_{\beta,l}(\t_m(f))= Z^{-1}\,
\sum_{\Gamma \backslash M_2(\Z)^+
} \,f(1, \mu \,\rho \,\tilde{m}^{-1}, \mu(\tau))\,{\rm Det}(\mu)^{-\beta}\,,
$$
where  $\mu\in M_2(\Z)^+$ is subject to the condition
$\mu \,\rho \,\tilde{m}^{-1} \in M_2(\RR)$.
The other values of $\mu$ a priori involved in the summation
\eqref{kmsexplicit}
do not contribute, since they correspond to the orthogonal of the
support of $\t_m(f)$.

\smallskip

\n One has ${\rm Det}(m)={\rm Det}(m')$ by construction, hence
$$
\rho \,\tilde{m}^{-1}=\,\rho\,m\; {\rm Det}(m)^{-1}=\, {\rm Det}(m')^{-1}
\,m'\,\rho'=\, \tilde{m'}^{-1}\,\rho'\,.
$$
Therefore the condition $\mu \,\rho \,\tilde{m}^{-1} \in M_2(\RR)$
holds iff $\mu = \nu \, \tilde{m'}$ for some $\nu \in M_2(\Z)^+$.
Thus, since ${\rm Det}(\mu)=\,{\rm Det}(\nu)\cdot{\rm Det}(\tilde
m')$, we can rewrite, up to multiplication by a scalar,
$$
\varphi_{\beta,l}(\t_m(f))= Z^{'\,-1}\,
\sum_{\Gamma \backslash M_2(\Z)^+
} \,f(1, \nu \,\rho' , \nu\,\tilde{m'}(\tau))\,{\rm Det}(\nu)^{-\beta}\,.
$$
This proves (\ref{comp}). It remains to show that on $\Ac_\Q$ we
have
$$
\varphi_{l'}(f)=\, {\rm Gal}_\tau(\rho\,m\,\rho^{-1})\circ
\varphi_l(f)\qqq f\in \Ac_\Q\,.
$$
Both sides only involve the values of $f$ on invertible
$\Q$-lattices, and there, by \eqref{galois} one has
$$
f(1,\,\a,\tau)= \, {\rm Gal}_\tau(\a) f(1,\,1,\tau) \qqq \a \in\GL(2,\RR)\,.
$$
Thus, we obtain
$$
\varphi_{l'}(f)=\,f(1,\rho',\,m^{'\,-1}(\tau))=
I_{m^{'\,-1}(\tau)}({\rm Gal}(\rho')f)=I_\tau({\rm Gal}(m'\,\rho')
f)\,.
$$
Since $ m'\, \rho'=\rho\, m $, this gives ${\rm
Gal}_\tau(\rho\,m\,\rho^{-1})\circ \varphi_l(f)$ as required.
$\Box$\\

\smallskip

\n We shall now  work out the algebraic relations
fulfilled by the $\Xc_a$ as extensions
of the division formulas of elliptic functions.

\smallskip

\n We first work with lattice functions of some weight $k$, or
equivalently with forms $f(g,y)\,dy^{k/2}$, and then multiply them
by a suitable factor to make them homogeneous of weight $0$ under
scaling. The functions of weight 2 are the generators, the higher
weight ones will be obtained from them by universal formulas with
modular forms as coefficients.

\smallskip

\n The powers $X_a^m$ of the function $X_a$ are then expressed as
universal polynomials with coefficients in the ring ${\mathbb Q}
\, [e_4 , e_6]$ in the following weight $2k$ functions ($k>1$):
\begin{equation}
\label{eq31} E_{2k,a} (\rho , \tau) = \pi^{-2k} \sum_{y \in
\Lambda + \phi (a)} y^{-2k}.
\end{equation}
These fulfill by construction (\cite{WeilEll}) the relations
\begin{eqnarray}
\label{weilpol}
E_{2m,a} &= &X_a (E_{2m-2,a} - e_{2m-2}) + \left( 1 - \begin{pmatrix} 2m \\
2 \end{pmatrix} \right) e_{2m} \\
&- &\sum_1^{m-2} \begin{pmatrix} 2k+1 \\ 2k\end{pmatrix} e_{2k+2}
(E_{2m-2k-2,a} - e_{2m-2k-2}) \, . \nonumber
\end{eqnarray}
These relations dictate the value of $E_{2k} (\rho , \tau)$ when
$\varphi (a) \in \Lambda$: one gets
\begin{equation}
\label{eq33} E_{2k} (\rho , \tau) = \nu_{2k} (\tau) \qquad
\hbox{if} \ \varphi (a) \in \Lambda,
\end{equation}
where $\nu_{2k}$ is a modular form of weight $2k$ obtained by
induction from (\ref{weilpol}), with $X_a$ replaced by $0$ and
$E_{2m}$ by $\nu_{2m}$. One has $\nu_{2k} \in {\mathbb Q} \, [e_4
, e_6]$ and the first values are
\begin{equation}
\label{eq34}
\nu_4 = -5 \, e_4 \, , \quad \nu_6 = -14 \, e_6 \, , \quad \nu_8 =
\frac{45}{7} \, e_4^2 \, , \ldots
\end{equation}

\smallskip

\n We shall now write the important algebraic relations between
the functions $X_a$, which extend the division relations of
elliptic functions from invertible ${\mathbb Q}$-lattices to
arbitrary ones.

\smallskip

\n In order to work out the division formulas for the Eisenstein series
$E_{2m,a}$ we need to
control the image of $\left( \frac{1}{N} \, {\mathbb Z} \right)^2
=\frac{1}{N} \Z^2$ under an
arbitrary element $\rho \in M_2 (R)$. This is done as follows
using the projections $\pi_L$ defined in (\ref{latticelabel0}).
\medskip

\begin{lem} \label{modN}
Let $N \in {\mathbb N}^*$, and $\rho \in M_2 (R)$. There exists a
smallest lattice $L \subset {\mathbb Z}^2$ with $L \supset N
{\mathbb Z}^2$, such that $\pi_L (\rho) = 1$. One has
$$
\rho ( \frac{1}{N} \, {\mathbb Z}^2) = \frac{1}{N}
\, L \, .
$$
\end{lem}

\noindent {\it Proof.} There are finitely many lattices $L$ with
$N {\mathbb Z}^2 \subset L \subset {\mathbb Z}^2$. Thus, the
intersection of those $L$ for which $\pi_L (\rho) = 1$ is still a
lattice and fulfills $\pi_L (\rho) = 1$ by (\ref{eq11}). Let $L$
be this lattice, and let us show that $\rho( \frac{1}{N} \Z^2)
\subset \frac{1}{N} \, L$. Let $m \in M_2 ({\mathbb Z})^+$ be such
that $m({\mathbb Z}^2) = L$. Then $\pi_L (\rho) = 1$ implies that
$\rho = m\mu$ for some $\mu \in M_2 (R)$. Thus $\rho (\frac{1}{N}
\Z^2) \subset m (\frac{1}{N} \Z^2) = \frac{1}{N} \, L$. Conversely
let $L' \subset {\mathbb Z}^2$ be defined by $\rho (\frac{1}{N}
\Z^2)
 = \frac{1}{N} \, L'$. We need to show that
$\pi_{L'} (\rho) = 1$, {\it i.e.} that there exists $m' \in M_2
({\mathbb Z})^+$ such that $L' = m' ({\mathbb Z}^2)$ and
${m'}^{-1} \, \rho \in M_2 (R)$. Replacing $\rho$ by $\gamma_1
\rho \, \gamma_2$ for $\gamma_j \in {\rm SL}_2 ({\mathbb Z})$ does
not change the problem, hence we can use this freedom to assume
that the restriction of $\rho$ to $\left( \frac{1}{N} \, {\mathbb
Z} \right)^2$ is of the form
\begin{equation}
\label{rhoNd}
\rho_N = \begin{bmatrix} d_1 &0 \\ 0 &d_2 \end{bmatrix} \, , \quad d_j \mid
N \, , \ d_1 \mid d_2 \, .
\end{equation}
One then takes $m' = \begin{bmatrix} d_1 &0 \\ 0 &d_2 \end{bmatrix} \in M_2
({\mathbb Z})^+$ and checks that $L' = m' ({\mathbb Z}^2)$ while ${m'}^{-1}
\rho$ belongs to $M_2 (R)$. $\Box$

\smallskip

\n Given an integer $N > 1$, we let $S_N$ be the set of lattices
\begin{equation}
\label{lattices} N \, {\mathbb Z}^2 \subset L \subset {\mathbb
Z}^2,
\end{equation}
which is the same as the set of subgroups of $({\mathbb Z} / N
{\mathbb Z})^2$. For each $L \in S_N$ we define a projection $\pi
(N,L)$ by
\begin{equation}
\label{eq36}
\pi (N,L) = \pi_L \prod_{L' \in S_N , L' \subsetneqq L} (1- \pi_{L'}) \, .
\end{equation}
By Lemma \ref{modN} the range of $\pi (N,L)$ is exactly the set of
$\rho \in M_2 (R)$ such that
\begin{equation}
\label{range}
\rho (\frac{1}{N} \, {\mathbb Z}^2) = \frac{1}{N}
\, L \, .
\end{equation}
The general form of the division relations is as follows.

\medskip

\begin{prop}  There exists canonical modular forms
$\omega_{N,L,k}$ of level $N$ and weight $2k$, such that for all
$k$ and $(\rho , \tau) \in Y $ they satisfy
$$
\sum_{Na = 0} X_a^k (\rho , \tau) = \sum_{L \in S_N} \pi (N,L) (\rho) \,
\omega_{N,L,k} (\tau) \, .
$$
\end{prop}

\medskip

\n In fact, we shall give explicit formulas for the
$\omega_{N,L,k}$ and show in particular that
$$
\omega_{N,L,k}(\gamma \, \tau)= (c\tau
+d)^{2k}\,\omega_{N,\gamma^{-1}\,L,k}(\tau),
$$
which implies that $\omega_{N,L,k}$ is of level $N$.

\smallskip

\n We prove it for $k=1$ and then proceed by
induction on $k$.
 The division formulas in weight $2$ involve the
$1$-cocycle on the group $ \displaystyle
 \GL^{+}_2(\mathbb{Q})  \,
$ with values in Eisenstein series of weight $2$ given in terms of
the Dedekind $\eta$-function by (cf. \cite{CM10})
\begin{equation} \label{dmu}
\mu_{\g} \, (\tau) \, = \, \frac{1}{12  \pi i} \, \frac{d}{d\tau}
\log \frac{\D | \g}{\D} \, = \, \frac{1}{2  \pi i} \,
\frac{d}{d\tau} \log \frac{\eta^4 | \g}{\eta^4} \, ,
\end{equation}
where we used the standard `slash operator' notation for the
action of $\GL^{+}_2 (\mathbb{R}) $ on functions on the upper half
plane:
\begin{equation} \label{act1}
f|_{k}\, \a \, (z)  =  {\rm Det}\, (\a)^{k/2} \, f (\a \cdot z) \,
j(\a , z)^{-k} \, ,
\end{equation}
$$ \a = \begin{bmatrix}a  &b\\
c   &d  \end{bmatrix} \in  \GL^{+}_2(\mathbb{R})
 , \quad
 \a \cdot z = \frac{az+b}{cz+d}
 \quad \hbox{and} \quad
j (\a, z) \, = \, cz+d \, .$$
Since $\mu_{\g} =0$ for $\g \in \Gamma$, the cocycle property
\begin{equation} \label{gmu}
\mu_{\g_{1}\cdot \g_{2}} \, = \, \mu_{\g_{1}}| \g_{2} \, + \,
\mu_{\g_{2}}
\end{equation}
shows that, for $m\in M_2(\Z)^+$, the value of $\mu_{m^{-1}} $
only depends upon the lattice $L= m(\Z^2)$. We shall denote it by
$\mu_L$.

\smallskip

\begin{lem}\label{div11} For any integer $N$, the $X_{a}$, $a \in
\mathbb{Q}\slash\mathbb{Z}$, fulfill the relation
$$
\sum_{N\,a = 0}\, X_{a}=\,N^2\,\sum_{L \in S_N} \pi (N,L)\;\,\mu_L
$$
\end{lem}

\n By construction the projections $\pi (N,L)$, $L \in S_N$ form a
partition of unity,
$$
\sum_{L \in S_N} \pi (N,L)=\,1\,.
$$
Thus to prove the lemma it is enough to evaluate
both sides on $\rho \in \pi (N,L)$.

\smallskip

\n We can moreover use the equality
$$
\mu_{\g^{-1}L}= \mu_L|\g \qqq \g \in \G
$$
to assume that $L$ and $\rho_N$ are of the form
$$
L = \begin{bmatrix} d_1 &0 \\ 0 &d_2 \end{bmatrix} \,\Z^2\,,
\qquad\rho_N = \begin{bmatrix} d_1 &0 \\ 0 &d_2 \end{bmatrix} \, ,
\quad d_j \mid
N \, , \ d_1 \mid d_2 \, .
$$
Let $d_2=n\,d_1$. The order of the kernel of  $\rho_N$
is $d_1\,d_2$ and the computation of
$\sum_{N\,a = 0}\, X_{a}(\rho,\tau)$ gives
$$
-N^2 e_2(\tau)\,+\,d_1\,d_2 \,N^2\,
\sum_{(a,b)\in\Z^2\backslash \{0\}}\,(a\, d_1-b\,d_2\tau)^{-2}
$$
which gives $N^2 (n \,e_2(n\,\tau)-e_2(\tau))=N^2 \mu_L$. \\

\n This proves the proposition for $k=1$. Let us proceed by
induction using (\ref{weilpol}) to express $X_a^k$ as $E_{2k,a}$
plus a polynomial of degree $<k$ in $X_a$ with coefficients in
$\Qb[e_4,e_6]$. Thus, we only need to prove the equality
$$
\sum_{N\,a = 0}\, E_{2k,a}=\,\sum_{L \in S_N} \pi
(N,L)\;\,\a_{N,L,k},
$$
where the modular forms $\a_{N,L,k}$ are given explicitly as
$$
\a_{N,L,k}=\, N^{2k}d^2n^{k+1} e_{2k}|m^{-1}\,-{\rm
Det}(m)\;(e_{2k}-\nu_{2k}),
$$
with $m\in M_2(\Z)^+$, $m(\Z^2)=L$, and $(d,dn)$ the elementary
divisor of $L$.

\n The proof is obtained as above by evaluating both sides on
arbitrary $\rho \in \pi (N,L)$.  $\Box$

\smallskip

\n One can rewrite all the above relations in terms of the
weight $0$ elements
$$
\Xc_a:=\,c\,X_a \,,\quad \Ec_{2k,a}:= \,c^k\, E_{2k,a}\in \Ac_\Q.
$$
In particular, the two basic modular functions $c^2\, e_4$ and
$c^3\, e_6$ are replaced by
$$
c^2\, e_4=\, \frac{1}{5}\,j\,(j- 1728)\,,\quad c^3\, e_6=\,
-\frac{2}{35}\,j\,(j- 1728)^2.
$$
We can now rewrite the relations (\ref{weilpol}) in terms of
universal polynomials
$$
P_n \in \Q(j)[X],
$$
which express the generators $\Ec_{2k,a}$ in terms of $\Xc_a$ by
$$
\Ec_{2k,a}=P_k(\Xc_a).
$$
In fact, from (\ref{weilpol}) we see that the coefficients of
$P_k$ are themselves polynomials in $j$ rather than rational
fractions, so that
$$
P_n \in \Q[j,X].
$$
The first ones are given by
$$
P_2=X^2-j(j-1728)\,,\quad P_3=X^3-\frac{9}{5}\,X
j(j-1728)+\frac{4}{5}j(j-1728)^2 \,,\quad \cdots
$$

\section{The noncommutative boundary of modular curves}\label{tower}

We shall explain in this section how to
combine the dual of the
$\GL_2$-system described above with the
idea, originally developed in the work of Connes--Douglas--Schwarz
\cite{CDS} and Manin--Marcolli \cite{ManMar}, of enlarging the
boundary of modular curves with a noncommutative space that
accounts for the degeneration of elliptic curves to noncommutative
tori.

\n The $\GL_2$-system described in the previous sections
admits a ``dual" system obtained by considering $\Q$-lattices
up to commensurability but no longer up to scaling.
Equivalently this corresponds to taking the cross product of the
$\GL_2$-system  by the action of the Pontrjagin dual of $\C^*$, which
combines the time evolution $\sigma_t$ with an action by the group $\Z$
of integral weights of modular forms. The resulting space is
the total space of the natural $\C^*$-bundle.

\smallskip

\n In ad\'elic terms this ``dual" noncommutative space $\Lc_2$
is described as follows,
\begin{prop}
There is a canonical bijection from the space of
$\GL_2(\Q)$-orbits of the left action of
$\GL_2(\Q)$ on $M_2(\A_f)\times \GL_2(\R)$ to the space
$\Lc_2$  of
commensurability classes of two-dimensional $\Q$-lattices.
\end{prop}

\noindent {\bf Proof.} The space of $\GL_2(\Q)$ orbits on
$M_2(\A_f)\times \GL_2(\R)$ is the same as the space of $\GL_2^+(\Q)$
orbits on $M_2(\RR)\times \GL_2^+(\R)$.  $\Box$

\smallskip

\n By the results of the previous sections, 
the classical space
obtained by considering the zero temperature
limit of  the quantum statistical mechanical
system describing commensurability classes of 2-dimensional
$\Q$-lattices up to scaling is the Shimura variety
that represents the projective limit of all the modular curves
\begin{equation}\label{GL2T0}
\GL_2(\Q) \backslash \GL_2(\A)/ \C^*\,.
\end{equation}

\smallskip

\n Usually, the Shimura variety is constructed as the projective limit
of the 
$$ \Gamma'\backslash \GL_2(\R)^+/\C^\ast $$ 
over congruence subgroups $\Gamma'\subset \Gamma$.
This gives a connected component in \eqref{GL2T0}. The other
components play a crucial role in the present context, in that
the 
existence of several connected components allows for
non-constant solutions of the equation $\zeta^n=1$.
Moreover all
the components are permuted by the Galois covariance property of
the arithmetic elements of the $\GL_2$ system.

\smallskip

\n  The  total space of
the natural $\C^*$-bundle, \ie the quotient
\begin{equation}\label{totspace}
\GL_2(\Q) \backslash \GL_2(\A),
\end{equation}
is the space of {\em invertible} 2-dimensional
$\Q$-lattices (not up to scaling).

\smallskip

\n  In the $\GL_1$ case, the analog of \eqref{totspace}, \ie
the space of id\`ele classes
$$\GL_1(\Q)\backslash \GL_1(\A),$$ is compactified by first
considering the noncommutative space of commensurability classes
of $\Q$-lattices not up to scaling
$$\Lc=\GL_1(\Q)\backslash
\A^\cdot ,$$
where $\A^\cdot$ is the space of ad\`eles with nonzero
archimedean component.
The next step, which is crucial in obtaining
the geometric
space underlying the spectral realization of the zeros of the
Riemann zeta function, is to add  an additional ``stratum'' that gives
the noncommutative
space of ad\`ele classes $$\overline{\Lc}=\GL_1(\Q)\backslash \A\,,$$
which will be analyzed in the next Chapter.

\smallskip

\n Similarly, in the $\GL_2$ case, the
classical space given by the Shimura variety \eqref{GL2T0} is
first ``compactified'' by adding noncommutative ``boundary
strata'' obtained by replacing $\GL_2(\A_f)$ in
$\GL_2(\A)=\GL_2(\A_f)\times \GL_2(\R)$ by all matrices $M_2(\A_f)$. 
As boundary stratum of the Shimura variety it corresponds to
degenerating the invertible $\Q$-structure $\phi$ on the lattice
to a non-invertible one and yields the notion of $\Q$-lattice.
 The corresponding
space of commensurability classes of $\Q$-lattices up
to scaling played a  central role in this whole chapter.

\smallskip

\n The space of commensurability classes of
2-dimensional $\Q$-lattices (not up to scaling) is
\begin{equation}\label{ltwo}
\Lc_2 = \GL_2(\Q)\backslash (M_2(\A_f)\times \GL_2(\R))\,.
\end{equation}
On $\Lc_2$ 
 we can consider
not just modular functions but all modular forms as functions. 
One obtains in this way an antihomomorphism of the
modular Hecke algebra of level one of \cite{CM10} (with variable
$\a \in\GL_2(\Q)^+$ restricted to $M_2(\Z)^+$) to the algebra of
coordinates on $\Lc_2$.

\smallskip

\n The further compactification at the archimedean place,
corresponding to $\Lc \hookrightarrow \overline{\Lc}$ in the $\GL_1$ case,
now consists of replacing $\GL_2(\R)$ by matrices $M_2(\R)$. This
corresponds to degenerating the lattices to pseudo-lattices (in
the sense of \cite{Man1}) or in more geometric terms, to a
degeneration of elliptic curves to noncommutative tori. It is this
part of the ``noncommutative compactification'' that was
considered in \cite{CDS} and \cite{ManMar}.

\smallskip

\n A $\Q$-pseudolattice in $\C$ is a pair $(\Lambda,\phi)$, with
$\Lambda=j(\Z^2)$ the image of a homomorphism $j:\Z^2 \to \ell$, with
$\ell\subset \R^2\cong \C$ a real 1-dimensional subspace, and with
a group homomorphism
$$ \phi: \Q^2/\Z^2 \to \Q \L/\L. $$
The $\Q$-pseudolattice is nondegenerate if $j$ is
injective and is invertible if $\phi$ is invertible.

\smallskip

\begin{prop}\label{pseudolatY}
Let $\partial Y:= M_2(\RR)\times \bP^1(\R)$. The map
\begin{equation}\label{thetaPsLat}
(\rho,\theta) \mapsto (\Lambda,\phi), \ \ \ \  \L=\Z+\theta \Z, \ \ \
\phi(x)=\rho_1(x) -\theta \rho_2(x)
\end{equation}
gives an identification
\begin{equation}\label{PsLat-space}
\Gamma \backslash \partial Y \simeq (\text{Space of $\Q$-pseudolattices
in $\C$})/ \C^*.
\end{equation}
This space parameterizes the degenerations of 2-dimensional $\Q$-lattices
\begin{equation}\label{lambdamap}
\lambda (y) = (\Lambda, \phi) \ \text{ where } \ \Lambda =
\tilde h({\mathbb Z} + i\, {\mathbb Z} ) \  \text{ and } \
\phi = \tilde h \circ \rho,
\end{equation}
for $y = (\rho , h) \in M_2(\RR)\times \GL_2(\R)$ and $\tilde h =
h^{-1} {\rm Det}\,(h)$, when $h\in \GL_2(\R)$ degenerates to a non-invertible
matrix in $M_2(\R)$.
\end{prop}

\n {\em Proof.} $\Q$-pseudolattices in $\C$ are of the form
\begin{equation}\label{pseudolattice}
\Lambda= \lambda( \Z + \theta \Z), \ \ \ \phi(a) = \lambda
\rho_1(a) - \lambda \theta \rho_2(a),
\end{equation}
for $\lambda \in \C^*$ and $\theta \in \bP^1(\R)$ and $\rho\in
M_2(\RR)$. The action of $\C^*$ multiplies $\lambda$, while
leaving $\theta$ unchanged. This corresponds to changing the
1-dimensional linear subspace of $\C$ containing the
pseudo-lattice and rescaling it. The action of $\SL_2(\Z)$ on
$\bP^1(\R)$ by fractional linear transformations changes $\theta$.
The nondegenerate pseudolattices correspond to the values
$\theta\in \bP^1(\R)\smallsetminus \bP^1(\Q)$ and the degenerate
pseudolattices to the cusps $\bP^1(\Q)$.

\smallskip

\n For $y = (\rho , h) \in M_2(\RR)\times
\GL_2(\R)$ consider the $\Q$-lattice \eqref{lambdamap},
for $\tilde h = h^{-1} {\rm Det}\,(h)$. Here we use the basis
$\{e_1=1,e_2=-i\}$ of the $\R$-vector space $\C$ to let ${\rm
GL}_2^+ ({\R})$ act on $\C$ as $\R$-linear transformations.
These formulas continue to make sense when
$h\in  M_2(\R)$ and the image
$\Lambda =
\tilde h({\mathbb Z} + i\, {\mathbb Z} ) $
is a pseudolattice when the matrix $h$ is  no longer invertible.

\smallskip

\n To see this more explicitly, consider the right action
\begin{equation}\label{rightCM2}
 m \mapsto m\cdot z
\end{equation}
of $\C^*$ on $M_2(\R)$ determined by the inclusion $\C^* \subset
\GL_2(\R)$ as in \eqref{inclGL2R},
The action of $\C^*$ on $M_2(\R) \smallsetminus \{ 0 \}$ is
free and proper. The map
\begin{equation}\label{rhomapM2}
 \rho (\alpha) = \left\{
\begin{array}{ll}
 \alpha(i) & (c,d)\neq (0,0) \\[2mm]
\infty & (c,d)=(0,0) \end{array} \right. \ \ \ \text{ with } \
\alpha = \left[
\begin{array}{cc} a & b \\ c & d \end{array}\right]
\end{equation}
defines an isomorphism
\begin{equation}\label{rhoisoM2}
 \rho: (M_2(\R)\smallsetminus \{ 0 \}) /\C^* \to \bP^1(\C),
\end{equation}
equivariant with respect to the left action of $\GL_2(\R)$ on
$M_2(\R)$ and the action of $\GL_2(\R)$ on $\bP^1(\C)$ by
fractional linear transformations. Moreover, this maps $M_2(\R)^+$
to the closure of the upper half plane
\begin{equation}\label{H2compact}
 \overline{\H}= \H \cup \bP^1(\R).
\end{equation}
The rank one matrices in $M_2(\R)$ map to $\bP^1(\R)\subset
\bP^1(\C)$.
In fact, the isotropy group of $m\in M_2(\R)$ is trivial if $m\neq
0$, since $m\cdot z =m$ only has nontrivial solutions for $m=0$,
since $z-1$ is invertible when nonzero. This shows that $M_2(\R)
\backslash \{ 0 \}$ is the total space of a principal
$\C^*$-bundle. $\Box$

\smallskip

\n Notice that, unlike the case of $\Q$-lattices of
\eqref{maptheta}, where the quotient $\Gamma\backslash Y$ can be
considered as a classical quotient, here the space $\Gamma\backslash
\partial Y$ should be regarded as a noncommutative space
with function algebra
$$ C( \partial Y ) \rtimes \Gamma. $$

\medskip

\n The usual algebro--geometric compactification of a modular
curve $Y_{\Gamma'} = \Gamma' \backslash \Hb$, for $\Gamma'$ a finite
index subgroup of $\Gamma$, is obtained by adding the cusp
points $\Gamma \backslash \bP^1(\Q)$,
\begin{equation}\label{cusps}
X_{\Gamma'} = Y_{\Gamma'} \cup \{ \text{ cusps } \} = \Gamma' \backslash
(\Hb \cup \bP^1(\Q))\, .
\end{equation}
Replacing $\GL_2(\R)$ by $M_2(\R)$
in \eqref{ltwo} corresponds to replacing the cusp points
$\bP^1(\Q)$ by the full boundary $\bP^1(\R)$ of $\Hb$. Since
$\Gamma$ does not act discretely on $\bP^1(\R)$, the quotient is
best described by noncommutative geometry, as the cross product
$C^*$-algebra $C(\bP^1(\R))\rtimes \Gamma'$ or, up to Morita
equivalence,
\begin{equation}\label{NCmodcurve}
C(\bP^1(\R)\times \bP) \rtimes \Gamma,
\end{equation}
with $\bP$ the coset space $\bP = \Gamma/\Gamma'$.

\smallskip

\n The noncommutative boundary of modular curves defined this way
retains a lot of the arithmetic information of the classical
modular curves. Various results of \cite{ManMar} show,
from the number theoretic point of view,
why the irrational points of $\bP^1(\R)$ in the boundary of $\Hb$
should be considered as part of the compactification of modular
curves.

\smallskip

\n The first such result is that the classical definition of
modular symbols (\cf \cite{Man3}), as homology classes on modular
curves defined by geodesics connecting cusp points, can be
generalized to ``limiting modular symbols'', which are asymptotic
cycles determined by geodesics ending at irrational points. The
properties of limiting modular symbols are determined by the
spectral theory of the Ruelle transfer operator of a dynamical
system, which generalizes the Gauss shift of the continued
fraction expansion by taking into account the extra datum of the
coset space $\bP$.

\smallskip

\n Manin's modular complex (\cf \cite{Man3}) gives a combinatorial
presentation of the first homology of modular curves, useful in
the explicit computation of the intersection numbers obtained by
pairing modular symbols to cusp forms. It is shown in
\cite{ManMar} that the modular complex can be recovered
canonically from the $K$-theory of the $C^*$-algebra
\eqref{NCmodcurve}.

\smallskip

\n Moreover, Mellin transforms of cusp forms of weight two for the
congruence subgroups $\Gamma_0(p)$,  with $p$ prime, can be
obtained by integrating along the boundary $\bP^1(\R)$ certain
``automorphic series'' defined in terms of the continued fraction
expansion and of modular symbols.

\smallskip

\n These extensions of the theory of modular symbols to the
noncommutative boundary appear to be interesting also in relation
to the results of \cite{CM10}, where the pairing with modular
symbols is used to give a formal analog of the Godbillon--Vey
cocycle and to obtain a rational representative for the Euler
class in the group cohomology $H^2(\SL_2(\Q),\Q)$.

\smallskip

\n The fact that the arithmetic information on modular curves is
stored in their noncommutative boundary \eqref{NCmodcurve} is
interpreted in \cite{ManMar2} as an instance of the physical
principle of holography. Noncommutative spaces arising at the
boundary of Shimura varieties have been further investigated by
Paugam \cite{Pau} from the point of view of Hodge structures.

\smallskip

\n This noncommutative boundary stratum of modular curves
representing degenerations of lattices to pseudo-lattices has been
proposed by Manin (\cite{Man1} \cite{Man2}) as a geometric space
underlying the explicit class field theory problem for real
quadratic fields. In fact, this is the first unsolved case of the
Hilbert 12th problem. Manin developed in \cite{Man1} a theory of
real multiplication, where noncommutative tori and pseudolattices
should play for real quadratic fields a role parallel to the one
that lattices and elliptic curves play in the construction of
generators of the maximal abelian extensions of imaginary
quadratic fields. The picture that emerges from this ``real
multiplication program'' is that the cases of $\Q$
(Kronecker--Weber) and of both imaginary and real quadratic fields
should all have the same underlying geometry, related to different
specializations of the $\GL_2$ system. The relation of the
$\GL_2$ system and explicit class field theory for imaginary
quadratic fields is analyzed in \cite{CMR}.

\bigskip
\section{The BC algebra and optical coherence}

\n It is very natural to look for concrete physical realizations
of the phase transition exhibited by the BC system. An
attempt in this direction has been proposed in \cite{planat}, in
the context of the physical phenomenon of quantum phase locking in
lasers.

\smallskip

\n This interpretation relates the additive generators $e(r)$ of
the BC algebra (\cf Proposition \ref{presentation}) with the
quantum phase states, which are a standard tool in the theory of
optical coherence (\cf \eg \cite{Lou}), but it leaves open the
interpretation of the generators $\mu_n$. Since on a finite
dimensional Hilbert space isometries are automatically unitary,
this rules out nontrivial representations of the $\mu_n$ in a
fixed finite dimensional space.

\smallskip

\n After recalling the basic framework of phase states and optical
coherence, we interpret the action of the $\mu_n$ as a
``renormalization'' procedure, relating the quantum phase states
at different scales.

\smallskip

\n There is a well known analogy (\cf \cite{laser} \S 21-3)
between the quantum statistical mechanics of systems with phase
transitions, such as the ferromagnet or the Bose condensation of
superconducting liquid Helium, and the physics of lasers, with the
transition to single mode radiation being the analog of
``condensation''. The role of the inverse temperature $\beta$ is
played in laser physics by the ``population inversion'' parameter,
with critical value at the inversion threshold. The injected
signal of the laser acts like the external field responsible for
the symmetry breaking mechanism. Given these identifications, one
in fact obtains similar forms in the two systems for both
thermodynamic potential and statistical distribution. The phase
locking phenomenon is also analogous in systems with phase
transitions and lasers, with the modes in the laser assuming same
phase and amplitude above threshold being the analog of Cooper
pairs of electrons acquiring the same energy and phase below
critical temperature in superconductors.

\smallskip

\n In a laser cavity typically many longitudinal modes of the
radiation are oscillating simultaneously. For a linewidth
$\Delta\nu$ around a frequency $\nu_0$ for the active medium in
the cavity of length $L$ and frequency spacing $\delta\nu = c/2L$,
the number of oscillating modes is $N=\lceil \Delta\nu/\delta\nu
\rceil$ and the field output of the laser is
\begin{equation}\label{fieldE}
E(x,t)= \sum_{n=-N/2}^{N/2} A_n \exp(-2\pi i\nu_n (t-x/c) + 2\pi i
\theta_n),
\end{equation}
with all the beat frequencies between adjacent modes
$\nu_n-\nu_{n-1}= \delta\nu$. Due to noise in the cavity all these
modes are uncorrelated, with a random distribution of amplitudes
$A_n$ and phases $\theta_n$.

\begin{figure}
\begin{center}
\epsfig{file=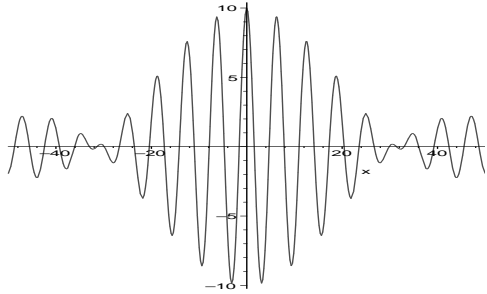}
\end{center}
\caption{Output pulse train in lasers above threshold.
\label{wave}}
\end{figure}

\n A mode locking phenomenon induced by the excited lasing atoms
is responsible for the fact that, above the threshold of
population inversion, the phases and amplitudes of the frequency
modes become locked together. The resulting field
\begin{equation}\label{fieldElock}
E(x,t)= A e^{2\pi i\theta} \exp(-2\pi i\nu_0 (t-x/c)) \left(
\frac{\sin (\pi \delta\nu (N+1)(t-x/c))}{\sin (\pi \delta\nu
(t-x/c))} \right),
\end{equation}
shows many locked modes behaving like a single longitudinal mode
oscillating inside the cavity (\cf Figure \ref{wave}). This
phenomenon accounts for the typical narrowness of the laser
linewidth and monocromaticity of laser radiation.

\smallskip

\smallskip

\n Since the interaction of radiation and matter in lasers is
essentially a quantum mechanical phenomenon, the mode locking
should be modeled by quantum mechanical phase operators
corresponding to the resonant interaction of many different
oscillators. In the quantum theory of radiation one usually
describes a single mode by the Hilbert space spanned by the
occupation number states $|n\rangle$, with creation and
annihilation operators $a^*$ and $a$ that raise and lower the
occupation numbers and satisfy the relation $[a, a^*]=1$. The
polar decomposition $a=S \sqrt{N}$ of the annihilation operator is
used to define a quantum mechanical phase operator, which is
conjugate to the occupation number operator $N=a^* a$. Similar
phase operators are used in the modeling of Cooper pairs. This
approach to the definition of a quantum phase has the drawback
that the one sided shift $S$ is not a unitary operator. This can
also be seen in the fact that the inverse Cayley transform of $S$,
which gives the cotangent of the phase, is a symmetric non
self-adjoint operator.

\smallskip

\n The emission of lasers above threshold can be described in
terms of coherent state excitations,
$$ |\alpha\rangle = \exp(-|\alpha|^2) \sum_n
\frac{\alpha^n}{(n!)^{1/2}} \, |n\rangle, $$ which are
eigenfunctions of the annihilation operators,
$$ a\, |\alpha\rangle = \alpha\, | \alpha \rangle. $$
These are quantum mechanical analogs of classical electromagnetic
waves as in \eqref{fieldE} \eqref{fieldElock}. One can show (\cf
\eg \cite{Lou} \S 7.4) that the field excitation in a laser
approaches a coherent state as the pumping increases to values
above the population inversion threshold, with the phase diffusion
governed by the equation of motion for quantum mechanical phase
states.

\smallskip

\n The problem in defining a proper quantum phase operator, due to
lack of self-adjointness, has been overcome by the following
approximation of the basic quantum operators on the Fock space
$\Hc$. One selects a scale, given by a positive integer $N \in \N$
and cuts down $\Hc$ to a finite dimensional subspace by the phase
state projector
$$ P_N = \sum_m | \theta_{m,N} \rangle \, \langle \theta_{m,N} |,
$$
where the orthonormal vectors $| \theta_{m,N} \rangle$ in $\Hc$
are given by
\begin{equation}\label{phasestate}
| \theta_{m,N} \rangle := \frac{1}{(N+1)^{1/2}} \sum_{n=0}^N
\exp\left( 2\pi i \frac{m\,n}{N+1} \right) \, | n \rangle .
\end{equation}
These are eigenvectors for the phase operators, that affect
discrete values given by roots of unity, replacing a continuously
varying phase.

\smallskip

\n This way, phase and occupation number behave like positions and
momenta. An occupation number state has randomly distributed phase
and, conversely, a phase state has a uniform distribution of
occupation numbers.

\smallskip

\n We now realize the ground states of the BC system as
representations of the algebra $\Ac$ in the Fock space $\Hc$ of
the physical system described above. Given an embedding $\rho:
\Q^{ab}\to \C$, which determines the choice of a ground state, the
generators $e(r)$ and $\mu_n$ (\cf Proposition \ref{presentation})
act as
$$ e(a/b)\,  | n \rangle = \rho(\zeta_{a/b}^n) \, | n \rangle, $$
$$ \mu_k \, | n \rangle = | kn \rangle. $$

\smallskip

\n In the physical system, the choice of the ground state is
determined by the primitive $N+1$-st root of unity
$$ \rho(\zeta_{N+1})= \exp (2 \pi i/ (N+1)). $$
One can then write \eqref{phasestate} in the form
\begin{equation}\label{phaseop}
| \theta_{m,N}\rangle = e\left( \frac{m}{N+1}\right) \cdot v_N .
\end{equation}
where we write $v_N$ for the superposition of the first $N+1$
occupation states
$$ v_N :=\frac{1}{(N+1)^{1/2}}\sum_{n=0}^N | n
\rangle. $$

\smallskip

\n Any choice of a primitive $N+1$-st root of unity would
correspond to another ground state, and can be used to define
analogous phase states. This construction of phase states brings in a
new hidden group of symmetry, which is different from the standard
rotation of the phase, and is the Galois group $\Gal(\Q^{ab}/\Q)$.
This raises the question of whether such symmetries are an artifact of
the approximation, or if they truly represent a property of the
physical system.

\smallskip

\n In the BC algebra, the operators $\mu_n$ act on the
algebra generated by the $e(r)$ by endomorphisms given by
\begin{equation}\label{muPmu}
\mu_n \, P(e(r_1),\ldots, e(r_k)) \, \mu_n^* =
\frac{\pi_n}{n^k}\sum_{ns =r}P(e(s_1),\ldots, e(s_k)) ,
\end{equation}
for an arbitrary polynomial $P$ in $k$-variables, with
$\pi_n=\mu_n\mu_n^*$ and $s=(s_1,\ldots,s_k)$,
$r=(r_1,\ldots,r_k)$ in $(\Q/\Z)^k$. In particular, this action
has the effect of averaging over different choices of the
primitive roots.

\smallskip

\n The averaging on the right hand side of \eqref{muPmu},
involving arbitrary phase observables $P(e(r_1),\ldots, e(r_k))$,
has physical meaning as statistical average over the choices of
primitive roots. The left hand side implements this averaging as a
renormalization group action.

\smallskip

\n Passing to the limit $N\to \infty$ for the phase states is a
delicate process. It is known in the theory of optical coherence
(\cf \eg \cite{mandel} \S 10) that one can take such limit only
after expectation values have been calculated. The analogy between
the laser and the ferromagnet suggests that this limiting
procedure should be treated as a case of statistical limit, in the
sense of \cite{wilson}. In fact, when analyzing correlations near
a phase transition, one needs a mechanism that handles changes of
scale. In statistical mechanics, such mechanism exists in the form
of a renormalization group, which expresses the fact that
different length or energy scales are locally coupled. This is
taken care here by the action of \eqref{muPmu}.

\bigskip
\bigskip
\bigskip

\n {\bf A.~Connes}

\n Coll\`ege de France, 3, rue Ulm, F-75005 Paris, France

\n I.H.E.S. 35 route de Chartres F-91440 Bures-sur-Yvette, France

\n Department of Mathematics, Vanderbilt University, TN-37240, USA

\n email: connes\@@ihes.fr, alain\@@connes.org

\medskip

\n {\bf M.~Marcolli}

\n Max--Planck Institut f\"ur Mathematik, Vivatsgasse 7,
D-53111 Bonn, Germany

\n email: marcolli\@@mpim-bonn.mpg.de

\end{document}